\documentclass[11pt]{amsart}

\usepackage{amsmath}
\usepackage{amssymb}
\usepackage{graphicx}
\usepackage[ps,dvips,xdvi,arrow,frame,matrix,graph,curve,cmtip]{xy}
\usepackage{enumerate}
\usepackage{float}

\newcommand{\N}{\mathbb{N}}

\newcommand{\signa}{{\rm sig}}

\newcommand{\gam}{\gamma}
\newcommand{\eps}{\varepsilon}
\newcommand{\Gam}{\Gamma}
\newcommand{\Gamz}{\Gamma_0}

\newcommand{\lam}{\lambda}

\newcommand{\R}{\mathbb{R}}
\newcommand{\C}{\mathbb{C}}
\newcommand{\Z}{\mathbb{Z}}
\newcommand{\D}{\mathbb{D}}
\newcommand{\s}{\mathbb{S}}

\DeclareMathOperator{\fix}{Fix}
\DeclareMathOperator{\stab}{Stab}
\DeclareMathOperator{\pstab}{pStab}
\DeclareMathOperator{\aut}{Aut}
\DeclareMathOperator{\diag}{diag}
\DeclareMathOperator{\spn}{span}
\DeclareMathOperator{\sym}{Sym}
\DeclareMathOperator{\MI}{MI}
\DeclareMathOperator{\GL}{GL}
\newcommand{\bigcupdot}{%
  \mathop{%
    \vphantom{\bigcup}%
    \mathpalette\setbigcupdot\cdot}\displaylimits} 
\newcommand{\setbigcupdot}[2]{\ooalign{\hfil$#1\bigcup$\hfil\cr\hfil$#2$\hfil\crcr}}

\newcommand{\iso}{\cong}

\newcommand{\dr}{{\, dr}}

\theoremstyle{plain} 
\newtheorem{thm}{Theorem}[section]

\newtheorem{prop}[thm]{Proposition}

\theoremstyle{definition}
\newtheorem{definition}[thm]{Definition}

\theoremstyle{remark}
\newtheorem{rem}[thm]{Remark}

\floatstyle{boxed}
\newfloat{Algorithm}{tbp}{lop} 
\floatname{Algorithm}{\sc{Algorithm}}

\begin{document}

\thanks{Partially supported by NSF Grant DMS-0074326}
\thanks{\today}

\title[Automated Bifurcation Analysis on Graphs]
{Automated Bifurcation Analysis for Nonlinear\\Elliptic Partial Difference Equations on Graphs}

\author{John M. Neuberger}
\author{N\'andor Sieben}
\author{James W. Swift}

\email{
John.Neuberger@nau.edu,
Nandor.Sieben@nau.edu,
Jim.Swift@nau.edu}

\address{
Department of Mathematics and Statistics,
Northern Arizona University PO Box 5717,
Flagstaff, AZ 86011-5717, USA
}


\begin{abstract}
We seek solutions $u\in\R^n$ to the semilinear elliptic partial
difference equation $-Lu + f_s(u) = 0$,
where $L$ is the matrix corresponding to the Laplacian operator
on a graph $G$ and $f_s$ is a one-parameter family of nonlinear functions.
This article combines the ideas introduced by the authors in two papers:
a) {\it Nonlinear Elliptic Partial Difference Equations on Graphs}
(J. Experimental Mathematics, 2006),
which introduces analytical and numerical techniques for solving such equations,
and
b) {\it Symmetry and Automated Branch Following for a Semilinear Elliptic PDE on a Fractal Region} (SIAM J. of Dynamical Systems, 2006),
wherein we present some of our recent advances concerning symmetry, bifurcation, and automation for PDE.

We apply the symmetry analysis
found in the SIAM paper to arbitrary graphs
in order to obtain better initial guesses for Newton's method, create informative graphics, and better understand the role of symmetry
in the underlying variational structure.
We use two modified implementations of the
gradient Newton-Galerkin algorithm (GNGA, Neuberger and Swift) to 
follow bifurcation branches in a robust way.
By handling difficulties that arise when encountering accidental degeneracies and higher-dimensional critical eigenspaces,
we can find many solutions of many symmetry types to the discrete 
nonlinear system.
We present a selection of experimental results which demonstrate our algorithm's capability
to automatically generate bifurcation diagrams and solution graphics starting with only an edgelist
of a graph.
We highlight interesting symmetry and variational phenomena.
\end{abstract}

\maketitle

\tolerance=10000

\begin{section}{Introduction.}

This paper considers {\it nonlinear partial difference equations} (PdE) on graphs.
In particular, we automate the bifurcation analysis
and branch following required for finding solutions
$u\in\R^n$ to the discrete nonlinear system
\begin{eqnarray}
\label{pde}
-Lu + f_s(u) = 0 .
\end{eqnarray}
Here, $L$ is the matrix corresponding to the Laplacian operator 
on a simple connected graph $G$
and $f_s:\R\to\R$ satisfies $f_s(0)=0$ and $f_s'(0) =s$.
The nonlinear term $f_s(u)\in\R^n$ is defined as a composition,
that is, $(f_s(u))_i=f_s(u_i)$.
The real number $s$ is treated as a bifurcation parameter.
The existence of the trivial solution $u = 0 \in \R^n$ is clear for all $s\in\R$, since $f_s(0)=0$.
By finding and following new, bifurcating branches of (generally) lesser symmetry we are able to find, within reason,
any solution that is connected by branches to the trivial branch.
Our code works for a wide range of nonlinearities.
It does not require that $f_s$ is odd, that the nonlinearity is superlinear
\cite{AR, CCN}, nor that $f_s$ has the form $f_s(t) = s t + H(t)$.
That being said, in this paper we choose $f_s$ to be the family of odd and superlinear functions defined by
$f_s(t)=st+t^3$ except when otherwise specified.
Our ultimate goal is to automate the process of accurately approximating
all solutions to Equation~(\ref{pde}) given only the edgelist for a given graph $G$, 
and then to sensibly present information about those solutions.

We first applied Newton's method to solve semilinear elliptic boundary value problems in \cite{NS},
where we sought solutions as critical points of an appropriate action
functional on a suitable function space.
This article combines the new ideas introduced in
\cite{Neu}, concerning nonlinear PdE on graphs, 
with the recent advances concerning symmetry, bifurcation, and automation
presented in our paper \cite{NSS2}.
By automating the symmetry analysis
and corresponding isotypic decompositions found in 
\cite{NSS2} 
for arbitrary graphs,
we are able to apply two modified implementations of the
gradient Newton-Galerkin algorithm (GNGA, see \cite{NS}) in order to completely automate
bifurcation branch following.
We are able to handle most difficulties that arise when encountering accidental degeneracies
(see Definition~\ref{accidental})
and high-dimensional critical eigenspaces.

In order to catalog experimentally found solutions according to symmetry and symmetry type
and understand the type of bifurcation that lead to the successful computation of each solution,
we make use of the automatic generation of the information described in  the {\it bifurcation digraph} 
(see \cite{NSS2}) corresponding to that experiment's underlying symmetry group.
A visual display of the bifurcation digraph can also be automatically generated for human use
in understanding the underlying variational structure and expected proliferation of solutions.

In our bifurcation diagrams, we show plots of $\|u\|_1$ versus $s$
for solutions $u$ to Equation~(\ref{pde}) with parameter $s$.
These diagrams can indicate by color or line type the symmetry of solutions,
which is invariant on each branch, or the Morse Index (MI) of solutions, which
typically changes at bifurcation and turning points.

In a generalization of the notion of a contour plot,
we have developed several visual representations of solutions to the discrete nonlinear problem (\ref{pde}). 
To a high degree, these plots too are automatically generated,
chosen where possible to
make the symmetry of solutions visible
and to yield a graphic that is informative and pleasing.

Our study of the finite dimensional semilinear elliptic PdE (\ref{pde})
closely follows the related works concerning the PDE
\begin{eqnarray}
\label{cpde}
	\left\{ \begin{array}{rl}
			\Delta u + f(u) = 0 & \hbox{in } \Omega \\
			\frac{\partial u}{\partial\eta}=0 & \hbox{on } {\partial \Omega},
			\end{array}
			\right .
\end{eqnarray}
as well as the similar zero-Dirichlet problem; see
\cite{CCN, Neu2}
and references therein.
The graph Laplacian $L$ corresponds to the negative Laplacian $-\Delta$ from PDE theory.
Both $L$ and $-\Delta$ have non-negative eigenvalues,
and both have zero eigenvalues corresponding to constant eigenvectors and eigenfunctions, respectively.
Note the sign difference between Equations~(\ref{pde}) and (\ref{cpde}).


Much is known about the spectrum of the graph Laplacian.
See, for example, \cite{Bap, Big, Chu}.
Most of the PdE literature concerns linear problems and/or positive solutions,
whereas we are interested in the existence and symmetry of all solutions,
in particular sign-changing ones, to nonlinear PdE.
Our first paper in this subject area \cite{Neu} contains a fairly thorough list of citations
relevant to the study of solutions to linear and nonlinear PdE,
e.g., works by A. Ashyralyev, S. S. Cheng, P. G. Kevrekidis et al, S. T. Liu, M. Lapidus, 
G. I. Marchuk, C. V. Pao, Yu. V. Pokorny\u\i, V. L. Pryadiev, P. Sobolevskii, J. C. Strikwerda, 
and G. Zhang.
Applications of limiting cases where we increase the number of vertices 
and use a scaling factor 
to approximate solutions to nonlinear PDE on fractals 
closely follow the linear results of 
R. S. Strichartz and A. Teplyaev.
Our survey article \cite{Neu2} summarizes some of our most relevant PDE results and provides
a list of open problems in that area. 
While PdE are generated whenever PDE are discretized via finite differences on a grid,
the PdE we study in this paper have few vertices and do not approximate PDE.
Our numerical techniques, symmetry analysis, and existence theorems \cite{CCN, Neu},
apply to both PdE and PDE.
By focusing on PdE, we can study large symmetry groups which would only arise from PDE
on domains with dimension 3 or greater.

In \cite{NS}, the Gradient Newton Galerkin Algorithm (GNGA) was developed to investigate
PDE (\ref{cpde}) using a basis of eigenfunctions of the corresponding (continuous) linear problem to
span a
suitably large finite dimensional subspace.
In \cite{NSS2} we adapted this algorithm in order to find many solutions 
of a PDE on a region with fractal
boundary, while in \cite{Neu} we used the entire basis since $n$ was 
small.
In the current work, we modify the GNGA slightly in two different ways;
the cylinder-augmented GNGA (cGNGA) is used to find initial solution points on new branches
near bifurcation points, and the tangent-augmented GNGA (tGNGA) is used to more effectively
follow solution branches.
A closely related approach to finding solutions to PDE with symmetry is found in
\cite{WangZhou1, WangZhou2}.


The paper is organized as follows.
In Section~\ref{prelims} we handle the preliminaries, stating definitions, theory, and notation
for graphs, graph Laplacians, symmetry, the variational method, and isotypic decompositions.
Section~\ref{preprocess} enumerates the various tasks we do prior to invoking the
continuation solver which
implements the GNGA to find solutions to Equation~(\ref{pde}).  In particular, we discuss
graph creation and layout from an edgelist, 
the computation of $\aut(G)$, bifurcation digraphs, the orthonormal basis of eigenfunctions of $L$,
and isotypic decomposition of symmetry-invariant fixed point spaces for bifurcation analysis.
We provide some details in Section~\ref{newton_w_c} concerning the 
implementations of the tGNGA, secant method, and cGNGA.
These three algorithms are used for following branches, 
finding bifurcation points, and finding new solutions
on bifurcating branches, respectively.
Section~\ref{solver} outlines our algorithms and heuristics for 
controlling the repeated application of the Newton and secant code to find representative branches 
from every conjugacy class of branches.  
Some postprocessing details for generating contour plots and bifurcation diagrams 
are given in Section~\ref{postprocess}.
Our main examples and numerical results are found in Section~\ref{examples}.
The concluding Section~\ref{conclusion} contains observations and ideas
for future refinements and applications of our methods.

\end{section}

\begin{section}{Preliminaries}
\label{prelims}

In this section we review background and notation for graph theory, symmetry, and the GNGA.

\subsection{Graphs}

Let $G=(V_G,E_G)$ be a simple connected graph with vertex set $V_G=\{v_1,\ldots,v_n\}$ 
and edge set $E_G$. The \emph{degree} of a vertex $v_i$ is denoted by $d(v_i)$.
An \emph{automorphism} of $G$ is a bijection $\alpha:V_G\to V_G$ such that 
$\{\alpha(v_i),\alpha(v_j)\}\in E_G$ if and only if
$\{v_i,v_j\}\in E_G$.
The \emph{symmetry group} of $G$ is the group $\aut(G)$ of automorphisms.
If $\pi$ is a permutation in $\s_n$ then we define $\alpha_{\pi}:V_G\to V_G$ by $\alpha_{\pi}(v_i)=v_{\pi(i)}$.
Not every permutation defines an automorphism of $G$ but 
every automorphism of $G$ is determined uniquely by a permutation.

\subsection{Cayley graphs}

We are often interested in graphs with prescribed symmetry. We construct these graphs as
decorated Cayley graphs \cite{White}.
Given a group $\Gamma$ and a set of generators $\Delta$, the \emph{Cayley color digraph} 
$\text{Cay}_{\Delta}\Gamma$ is a directed labeled graph $(G,c)$  with vertex set 
$V_G=\Gamma$ and edge set 
\[
E_G=\{(g,gd)\mid\ g\in \Gamma, d\in \Delta\}.
\] 
Edge $(g,gd)$ is labeled with the color $c(g,gd)=d$. 
%
%
The group $G$ acts on
$\text{Cay}_{\Delta}\Gamma$ by left multiplication;
in fact,
$\aut( \text{Cay}_{\Delta}\Gamma) \cong \Gamma$. 
To create a simple undirected graph whose automorphism group is $\Gamma$, we replace the directed colored
edges of the Cayley color graph with undirected decorated edges. The decoration adds extra vertices along
the edges.
The resulting simple graph is called a \emph{decorated Cayley graph}.

\subsection{Graph Laplacian}
\label{graphLap}

The {\em Laplacian} of $G$ is determined by the matrix $L$ defined be
letting $L_{ii}=d(v_i)$, $L_{ij}=-1$ if $\{v_i,v_j\}\in E_G$, 
and $L_{ij}=0$ if $i\not =j$ but $\{v_i,v_j\}\not\in E_G$.
This Laplacian can be viewed as enforcing the zero Neumann boundary condition \cite{Big}.
For example, if we solve the appropriately scaled Equation~(\ref{pde}) on the path $P_n$
for large $n$, we get approximate solutions to Equation~(\ref{cpde}).
Consideration of other boundary conditions is possible and interesting, but
is the subject of other and future reports.
The incidence (first difference) matrix $D$ of an
arbitrary orientation of $G$ satisfies $L=D^TD$;
we do not use this fact but observe that the variational equations for PdE most closely resemble
those for PDE when expressions like $Lu\cdot v$ are replaced with $Du\cdot Dv$ (see~\cite{Neu}).
The eigenvalues and corresponding eigenvectors of $L$ are denoted by 
$0=\lambda_1<\lambda_2\leq\cdots\leq\lambda_n$ and $\{\psi_j\}_{j=1}^n$, respectively.

Let
$X=\{(u,s)\in \R^n\times \R \mid -Lu+f_s(u)=0\}$
be the
\emph{solution set} of PdE~(\ref{pde}).
We write 
\[
u=\sum_{j=1}^m a_j\psi_j\in U_m:=\spn\Psi_m\subset \R^n,
\]
where $\Psi_m=\{\psi_1,\ldots,\psi_m\}$ is an orthonormal set of eigenvectors of $L$, 
and use the notation
$ ([u]_{\Psi_m},s) := (a,s) \in\R^m\times\R$. 
Since we are working with modest sized graphs, in this paper we take $m$ to be $n$.

Where possible, when multiple eigenvalues are encountered, choices of the associated eigenvectors in
$\Psi_m$ are made to respect symmetry in a similar fashion as was done in \cite{NSS}.  
For details see Section~\ref{basisgen}.

\subsection{Symmetry of functions}
\label{functSym}

To study the symmetry of solutions to Equation~(\ref{pde}) we consider
$\Gamma_0=\aut(G)\times \Z_2$, where $\Z_2=\{1,-1\}$ is written multiplicatively. 
The natural action of $\Gamma_0$ on $\R^n$ is defined by
\begin{equation}
\label{Gamma_action}
(\gamma\cdot u)_i=\beta u_{\pi^{-1}(i)}, 
\end{equation}
where $\gamma=(\alpha_\pi,\beta)\in \Gamma_0$ and $u\in \R^n$.
We usually write $\alpha$ for $(\alpha,1)$ and $-\alpha$ for $(\alpha,-1)$. 
The \emph{symmetry} of $u$ is the isotropy subgroup $\sym(u):=\stab(u,\Gamma_0)=\{\gamma\in\Gamma_0\mid \gamma\cdot u=u\}$.
Two subgroups $\Gamma_i$ and $\Gamma_j$ of $\Gamma_0$ are called \emph{conjugate} if $\Gamma_i=\gamma\Gamma_j \gamma^{-1}$ for some 
$\gamma\in\Gamma_0$. The \emph{symmetry type} of $u$ is the conjugacy class $[\sym(u)]$ of the symmetry of $u$.
We use the notation ${\mathcal G}:=\{\Gamma_0,\ldots,\Gamma_q \}$ for the set of symmetries and 
${\mathcal S}:=\{S_0=[\Gamma_0],\ldots ,S_r\}$ for the set of symmetry types.

In general it is difficult to compute $\mathcal G$,
but the following definition helps for some graphs.
A \emph{generic vertex} of a graph is a vertex $v$ such that $\{ \alpha \in \aut(G) \mid \alpha (v) = v\}$
contains only the identity map.
The $\aut(G)$ orbit of a generic vertex has the same size as $\aut(G)$.
The proof of the following proposition follows \cite{NSS2}.
\begin{prop}
\label{all_sym_prop}
If the graph $G$ has a generic vertex, then 
${\mathcal G} = \{ \Gam \leq \Gam_0 \mid \Gam = \Gam_0 \text{ or } -1 \not \in \Gam \}$.
\end{prop}
\begin{proof}
Assume that $G$ has a generic vertex, which we label as $v_1$.
Consider the function $u$ such that $u_1 = 1$ and $u_i = 0$ for $i > 1$.
Then for any subgroup $\Gam \leq \Gam_0$ the function
$$
\sum_{\gam \in \Gam} \gam \cdot u
$$
has symmetry $\Gam$ if $-1 \not \in \Gam$, and symmetry $\Gam_0$ otherwise.
On the other hand, only $u = 0$ satisfies $-u = u$.
So, if $\Gam \leq \Gam_0$ is an isotropy subgroup containing $-1$ then $\Gam = \Gam_0$.
\end{proof}
\begin{rem}
We have a counterexample which shows that the converse of Proposition~\ref{all_sym_prop} is false.
Our counterexample $G$ is the union of
a decorated Cayley graph of $\Z_3$ and a decorated Cayley graph of $\Z_5$,
with 15 additional edges joining each element in $\Z_3$ with each element in $\Z_5$.
The symmetry group of $G$ is $\aut(G) \cong \Z_3 \times \Z_5$.
The set of symmetries $\mathcal G$
of $G$
consists of $\Gam_0$ and all the subgroups
of $\Z_3 \times \Z_5$,
but $G$ has no generic vertex.
\end{rem}
A decorated Cayley graph of any group automatically has a generic vertex, namely 
any of the vertices corresponding to elements of the group.
Graphs with generic vertices are good models for PDE where the domain $\Omega$ has
a particular symmetry.
If we can find a graph $G$ such that $\aut (G) = \aut (\Omega)$
and $G$ has a generic vertex, then $\mathcal G$ for $G$ is the same as the set of
possible symmetries of solutions to PDE (\ref{cpde}) on $\Omega$.

We define a \emph{branch of solutions} to be a maximal subset of the solution 
space $X$ that is
a $C^1$ manifold with constant symmetry.
The \emph{trivial branch} $\{ (0, s) \mid s \in \R\}$ contains the \emph{trivial solution} $u = 0$,
which has symmetry $\Gamma_0$ 
if $f_s$ is odd, and symmetry $\aut(G)$ otherwise.
The positive \emph{constant branch} is $\{ ((c, \ldots , c), s) \mid f_s(c) = 0, c > 0, s \in \R\}$.
The negative constant branch is similarly defined.
A \emph{bifurcation point} is a solution in the closure of at least two different solution branches.
We call the branch containing the bifurcation point the {\em mother}, and the other branches the {\em daughters}.
For example, the (positive and negative) constant branches are daughters of the trivial branch, which contains
the bifurcation point $(0, 0) \in \R^{n}\times \R$.

Equation~(\ref{pde}) can be interpreted as $\nabla J_s (u) = 0$, where $\nabla J_s: \R^n \rightarrow \R^n$
is defined by $-\nabla J_s(u) = -L u + f_s(u)$.
The operator $\nabla J_s$ is $\aut(G)$-equivariant, i.e., $\nabla J_s(\alpha u) = \alpha \cdot \nabla J_s(u)$ for all $\alpha\in\aut(G)$.
Furthermore, if $f_s$ is odd, then  $\nabla J_s$ is $\Gam_0$-equivariant.
If $u$ is a solution to Equation~(\ref{pde}) with $f_s$ odd, then
$\gamma \cdot u$ is also a solution to Equation~(\ref{pde}) for all $\gamma \in \Gamma_0$.
Following the standard treatment~\cite{GSS, NSS2},
for each $\Gam_i \leq \Gam_0$ we define
the {\em fixed point subspace} of the $\Gamma_0$ action on $V = \R^n$ to be
$$
\fix(\Gamma_i, V) = \{ u \in \R^n \mid \gamma \cdot u = u \text{ for all } \gamma \in \Gamma_i\}.
$$
If $f_s$ is odd,
these fixed point subspaces are $\nabla J_s$-invariant.
Otherwise, $\fix(\Gamma_i, V)$ is $\nabla J_s$-invariant for all $\Gam_i \leq \aut(G)$.
Recall that a subspace $W \subset \R^n$ is $\nabla J_s$-invariant if
$\nabla J_s(W) \subset W$.
We say that a subspace ${\mathcal A} \subset \R^n$
is an {\it anomalous invariant subspace} (AIS) if it is $\nabla J_s$-invariant but is not a fixed point subspace. 
If $\mathcal A$ is an AIS, we sometimes say $u \in {\mathcal A}$ is {\em anomalous}.
Note that the {\em constant subspace} $W_c := \{ (c, \ldots, c ) \mid c \in \R \}$ is the fixed point
subspace $\fix(\aut G, \R^n)$ if $G$ is vertex transitive.
Otherwise, the constant subspace is an AIS, which we denote by ${\mathcal A}_c$.
The book~\cite{GS} is a good reference on invariant spaces of nonlinear operators, although our definition
of anomalous invariant subspaces appears to be new.

We define an {\em anomaly-breaking bifurcation} to be a bifurcation where the
daughters have the same symmetry
as the mother, and
the mother is in an AIS that does not contain at least one of the daughters.
We have not been able to describe a general theory for anomaly-breaking bifurcations.

\subsection{Isotypic Decomposition}
\label{isoSub}
To analyze the bifurcations of a branch of solutions with symmetry $\Gamma_i$,
we need to understand the isotypic decomposition of the action of $\Gamma_i$ on
$\R^n$. 

Suppose a finite group $\Gamma$ acts on $V=\R^n$ according to the representation 
$g\mapsto \alpha_g:\Gamma\to \aut(V)\cong {\GL}_{n}(\R)$. 
In our applications we choose $\Gamma\in{\mathcal G}$ and the 
group action is the one in Equation~(\ref{Gamma_action}).
Let $\{\alpha^{(k)}_{\Gamma}:\Gamma\to {\GL}_{d_\Gamma^{(k)}}(\R)\mid k \in K_\Gamma\}$ be the set of 
irreducible representations of $\Gamma$ over $\R$. 
We write $\alpha^{(k)}$ and $K$ when the subscript $\Gam$ is understood.
It is a standard result of representation theory that
there is an orthonormal basis $B_{\Gamma}=\bigcup_{k\in K} B_{\Gamma}^{(k)}$ 
for $V$ such that
$B_\Gamma^{(k)}=\bigcupdot_{l=1}^{\,L_k} B_\Gamma^{(k,l)}$ and 
$[\alpha_g|_{V_{\Gamma}^{(k,l)}}]_{B_\Gamma^{(k,l)}}=\alpha^{(k)}(g)$ for all $g\in\Gamma$,  
where $V_{\Gamma}^{(k,l)}:=\spn(B_\Gamma^{(k,l)})$.
Each $V_{\Gamma}^{(k,l)}$ is an irreducible subspace of $V$.
Note that $B_\Gamma^{(k)}$ might be empty for some $k$, corresponding to $V_\Gamma^{(k)}=\{0\}$.
The {\em isotypic decomposition of } $V$ under the action of $\Gamma$ is
$$
V = \bigoplus_{k\in K} V_{\Gamma}^{(k)}, 
$$
where $V_{\Gamma}^{(k)}=\bigoplus_{l=1}^{L_k} V_{\Gamma}^{(k,l)}$ are the {\em isotypic components}.

The isotypic decomposition of $V$ under the action of each $\Gamma_i$ is required by
our algorithm.  The decomposition under the action of $\aut (G)$ is the same
as the decomposition under the action of  $\Gamma_0$.
While there are twice as many irreducible representations of  $\Gamma_0 = \aut (G) \times \Z_2$
as there are of $\aut (G)$, if $\alpha^{(k)}_{\Gamma_0}(-1) = I$ then $V_{\Gamma_0}^{(k)} = \{ 0 \}$.
The other half of the irreducible representations have $\alpha^{(k)}_{\Gamma_0}(-1) = -I$.
The irreducible representations of $\Gamma_0$ and of $\aut (G)$ can be labeled so that
$V_{\Gamma_0}^{(k)} = V_{\aut (G)}^{(k)}$ for $k \in K_{\aut (G)}$.

The isotypic components are uniquely determined, but the decomposition into 
irreducible spaces is not.
Our goal is to find $B_\Gamma^{(k)}$ for all $k$ by finding the projection 
$P_\Gamma^{(k)}:V\to V_\Gamma^{(k)}$. 
To do this, we first need to
introduce representations over the complex numbers $\C$ for two reasons.
First,
irreducible representations over $\C$ are better understood than those over $\R$.
Second, our GAP program uses the field $\C$ since
irreducible representations over $\R$ are not
readily obtainable by GAP.

There is a natural action of $\Gamma$ on $W:=\C^n$ given by the representation 
$g\mapsto \beta_g:\Gamma\to \aut(W)$ such that $\beta_g$ and $\alpha_g$ have the 
same matrix representation. The isotypic decomposition 
$W=\bigoplus_{k\in \tilde K} W_{\Gamma}^{(k)}$ is defined as above using the set
$\{\beta^{(k)}:\Gamma\to {\GL}_{\tilde d_\Gamma^{(k)}}(\C)\mid k \in \tilde K_\Gam\}$ of 
irreducible representations of $\Gamma$ over $\C$.

The {\em characters} of the irreducible representation $\beta^{(k)}$ are 
$\chi^{(k)}(g) := {\rm Tr} \, \beta^{(k)}(g)$. The projection $Q_\Gamma^{(k)}:W\to W_\Gamma^{(k)}$
is known to be
\begin{equation}
\label{projQ}
Q_{\Gamma}^{(k)} ={\frac{\tilde d_\Gamma^{(k)}}{|\Gamma|}}\sum_{g\in \Gamma} \chi^{(k)}(g) \beta_g .
\end{equation}
We are going to get the $P_\Gamma^{(k)}$'s in terms of $Q_\Gamma^{(k)}$'s with the help
of the Frobenius-Schur indicator
$$
\nu^{(k)} := 
\frac{1}{|G|}
\sum_{g \in G} \chi^{(k)}(g^2) \in \{ -1, 0, 1\}.
$$
Recall~\cite{Dornhoff} that
\begin{enumerate}[\upshape (i)]
\item $\nu^{(k)}=1$  implies $\chi^{(k)} = \overline{ \chi^{(k)} }$, 
    in which case we say $ \beta^{(k)} $ is a \emph{real} irreducible representation;
\item $\nu^{(k)}=0$  implies $ \chi^{(k)} \neq \overline{ \chi^{(k)} }$, 
    in which case we say $ \beta^{(k)}$ is \emph{complex};
\item $\nu^{(k)}=-1$ implies $ \chi^{(k)} = \overline{ \chi^{(k)}}$, 
    in which case we say $ \beta^{(k)}$ is \emph{quaternionic}. 
\end{enumerate}
Sometimes the term \emph{quasi-real} is used in place of quaternionic. 
For $k,k'\in\tilde K_\Gam$, 
we say $k\sim k'$ if $\chi^{(k)}=\overline{\chi^{(k')}}$. 
Complex representations come in complex conjugate pairs, so $\sim$ is an equivalence relation. 
Then $K_\Gam$ can be chosen to be any
complete set of representatives of the quotient set $\tilde K_\Gam/\sim$. 
We calculate the projection operators in $\R^n$ from the projection operators in $\C^n$ using the following formulas:
\begin{enumerate}[\upshape (i)]
\item if $\nu^{(k)}=1$ then $P_\Gam^{(k)}=Q_\Gam^{(k)}\mid_V$ and 
	$d_\Gamma^{(k)}=\tilde d_\Gamma^{(k)}$;
\item if $\nu^{(k)}=0$ then $P_\Gam^{(k)}=\left(Q_\Gam^{(k)}+\overline{Q_\Gam^{(k)}}\right)\mid_V$ and 
	$d_\Gamma^{(k)}=2\tilde d_\Gamma^{(k)}$;
\item if $\nu^{(k)}=-1$ then $P_\Gam^{(k)}=Q_\Gam^{(k)}\mid_V$ and 
	$d_\Gamma^{(k)}=2\tilde d_\Gamma^{(k)}$,
\end{enumerate}
for all $k\in K_\Gam$.

\subsection{GNGA}

We now review the GNGA for PdE \cite{Neu}.
Let $F_s:\R\to\R$ be the primitive defined by $F_s(t)=\int_0^t f_s(r) \dr$,
e.g., $F_s(t)= \frac12s t^2 + \frac14 t^4$.
The \emph{action} functional 
$J_s:\R^n\to\R$ is defined by
$$
J_s(u) = {\textstyle \frac12} Lu\cdot u - \sum_{i=1}^n F_s(u_i).
$$
For $u,v\in\R^n$ it is easy to see that
$$
J_s'(u)(v) = -(-Lu + f_s(u))\cdot v,
$$
so that $u$ is a critical point of $J_s$ if and only if $(u,s)\in X$,
i.e., $u$ is a solution to Equation~(\ref{pde}) for parameter $s$.

For $u \in U_m$, we compute the the coefficients of the gradient vector $g_s(u)\in\R^m$ by
\begin{eqnarray}
\label{grad}
g_s(u)_j = Lu\cdot \psi_j - f_s(u)\cdot\psi_j =
(L\sum_{k=1}^m a_k \psi_k)\cdot\psi_j - f_s(u)\cdot\psi_j =  a_j \lambda_j - f_s(u)\cdot\psi_j.
\end{eqnarray}
If $m = n$, then $g_s(u) = 0$ if and only if $\nabla J_s(u) = 0$.
When $m < n$, the solutions to  $g_s(u) = 0$ are approximate solutions to Equation~(\ref{pde}).
In this paper we assume $m =n$,
but the formulas use $m$ where appropriate, keeping in mind applications to large graphs,
e.g., those arising from a PDE.

Similarly, the Hessian matrix $h_s(u) = (J_s''(u)(\psi_j,\psi_k))_{j,k=1}^m$ can be computed as
\begin{eqnarray}
\label{hess}
h_s(u)_{jk} = L\psi_j\cdot \psi_k - \diag(f_s'(u))\psi_j\cdot\psi_k =
\lambda_j\delta_{jk} - \diag(f_s'(u))\psi_j\cdot\psi_k,
\end{eqnarray}
where $\delta_{jk}$ is the Kronecker delta and $\diag(f_s'(u))$ is a diagonal matrix.
Using the coefficient vector $a$ and eigenvalues $\{\lambda_j\}_{j=1}^m$ to compute the difference terms in 
$Lu\cdot\psi_j$ and $L\psi_j\cdot\psi_k$ significantly
reduces the number of matrix and vector operations required to define the linear system for a
search direction $\chi$ satisfying $h_s(u)\chi=g_s(u)$.
Applying Newton's method to find zeroes of $(u,s)\mapsto g_s(u)$ is the basis of our gradient Newton-Galerkin algorithms (GNGA).

We define the \emph{signature} $\signa(u,s)$ to be the number of negative eigenvalues of the matrix $h_s(u)$
representing the self-adjoint bilinear operator $D^2J_s(u)$.
If $(u,s)$ is a nondegenerate solution to Equation~(\ref{pde}), then sig$(u,s)$ equals the \emph{Morse index} $\MI(u,s)$.
The MI can be thought of as the number of ``down'' directions of the critical point,
that is, $\MI(u,s) = 0$ for minima of $J_s$, 
$\MI(u,s) = n$ for maxima, and $\MI(u,s)\in\{1,\ldots,n-1\}$ for saddle points in between.
The search direction $\chi$ can be found using any number of linear solvers;
we use a minimum norm least squares solver to avoid problems with
noninvertible Hessians $h_s(u)$.
Noninvertible Hessians inevitably occur at bifurcation points, and \emph{fold points} (points
where the solution branch is not monotonic in $s$).
When the Hessian is singular, the eigenspace of the Hessian with eigenvalue 0 is
called the {\em critical eigenspace}, and is denoted by $E$. 

\end{section}

\begin{section}{Preprocessing}
\label{preprocess}

In this section we describe the various tasks that must be performed prior to 
approximating solutions to Equation~(\ref{pde}).
In particular, we must create the edgelist, visualize the graph, analyze the symmetry of the problem,
compute the possible bifurcations,
and generate suitably bases from the eigenfunctions of the Laplacian.
The data files generated during preprocessing are used by the continuation solver,
as well as in the postprocessing phase when creating graphics in order to visualize the results.
All of these files for a single graph, Example~\ref{cycle_sub},
can be found at the website
\noindent
\begin{verbatim} http://NAU.edu/Jim.Swift/PdE. \end{verbatim}

\subsection{Graph Creation}

A graph $G$ is determined by an edgelist file.
Each line contains a pair of integers $i$ and $j$, indicating that $\{v_i, v_j\} \in E_G$.
This file is usually created by a text editor. 
We also have the option to create the edgelist file automatically by GAP \cite{GAP} as a decorated Cayley
graph of a given group (see Section~\ref{Qexample}). This is the main human input for our process. 

\subsection{Graph Layout Code}

\label{layout}
To create a visualization of the graph,
we use a standard spring embedding algorithm to create an \emph{embedding} $\ell:V_G\to \R^2$  
of the graph. 
This is done by a C++ program.
The program starts with a random placement of the vertices, that is, $\ell$ is initialized with random values. 
We  then calculate the ``force'' $F_i=E_i+H_i$ on each vertex $v_i$, where $E_i$ is generated by
an equal ``electric charge'' $Q$  on the vertices and $H_i$ is generated by ``springs'' of natural length 
$\nu$ replacing the edges of the graph. Specifically, with $d_{ij}=\ell(v_i)-\ell(v_j)$, we have
\[
   \begin{aligned}
    E_i = \sum_{j\neq i} \frac{Q^2}{\eps+\|d_{ij}\|^D}\,d_{ij}, \qquad
    H_i = \sum_{\{v_i,v_j\}\in E_G} \frac{\nu-\|d_{ij}\|}{\|d_{ij}\|}\,d_{ij}.
   \end{aligned}
\]
Experiments show that the values $Q^2=1$, $\eps=0.001$, $\nu=1$ and $D=1.1$ work well.
Iteratively replacing $\ell(v_i)$ by $\ell(v_i)+\delta F_i$ using a stepsize of $\delta=0.1$, we
simulate the movement of this physical system with added damping until an equilibrium is reached.
This stable position usually shows some aspects of the symmetry of the graph.
The {\em complexity} of the layout is defined to be the number of distinct distances between vertices.
The program tries several initial positions and picks the layout that minimizes the complexity.
Layouts with higher complexity are also stored for possible use.
Finally, we rotate the optimal placement so that the most common edge slope is horizontal. 
The output of the program is a file containing the coordinates of the vertices. We create
figures automatically from this file using Gnuplot, \Xy-pic and Mathematica, together with solution data
generated by the continuation solver. 

\subsection{Automorphism Group Code}
To analyze the symmetry of the graph we need to find its automorphism group $\aut(G)$. This is done by Nauty
\cite{nauty}, which is a very efficient program that can handle fairly large graphs. It creates a file
containing all the permutations or only the generators of the automorphism group. This file is the
input of the GAP program that performs the full symmetry analysis.

\subsection{Symmetry Analysis Code}
\label{sym_analy}

In this subsection we give some details of the GAP computations done to analyze the symmetry of the problem.
Since GAP uses irreducible representations over the complex numbers, some care must be taken.

To compute the set of symmetries
$\mathcal G$, we need the following definition.
If $\Gamma$ acts on $V$ and $U$ is a subspace of 
$V$ then
\[
\pstab(U,\Gamma)=\{ \gamma\in\Gamma \mid \gamma\cdot u = u \text{ for all } u\in U \} .
\]
The isotropy subgroups of 
the $\Gamma_0$ action on $\R^n$ are precisely the subgroups
$\Gamma$ of the finite group $\Gamma_0$ which satisfy  
\[
\pstab(\fix(\Gamma,V),\Gamma)=\Gamma.
\]
This computation is easily performed by GAP.

Next, we determine the possible symmetries of the daughters of a bifurcation point with 
symmetry
$\Gamma_i\in\mathcal{G}$.
For each $i$,
we use GAP to find the irreducible representations
$\{\alpha_{\Gamma_i}^{(k)} \mid k\in K_{\Gam_i}  \}$ of $\Gam_i$, and the characters $\chi^{(k)}$.
The characters are used to produce the projection operators $Q_{\Gam_i}^{(k)}$ defined in
Equation~(\ref{projQ}).
The isotypic components
$V_{\Gamma_i}^{(k)}$ of the $\Gamma_i$ action on $W = \C^n$ are computed as the row spaces of the projection operators.
The kernel of the irreducible representation, denoted $\Gam_{i,k}'$, is also computed by GAP.
Then for each $k\in K_{\Gam_i}$ for which $V_{\Gamma_i}^{(k)}$ is nontrivial
we generate the set $\mathcal{H}_{i,k}$ of isotropy subgroups of the 
$\Gamma_i$ action on $V_{\Gamma_i}^{(k)}$.
The set $\mathcal{H}_{i,k}$ is partially ordered with
$\Gam_i$ at the top and $\Gam_{i,k}'$ at the bottom, and is often called the lattice of isotropy subgroups \cite{GSS, NSS2}.
If there are no subgroups in $\mathcal{H}_{i,k}$ properly between $\Gamma_i$ and $\Gamma_j$ for some 
$\Gamma_j\in\mathcal{H}_{i,k}$ then
$\Gamma_j$ is called a \emph{maximal isotropy subgroup}.
For each of these maximal isotropy subgroups there is a possible generic bifurcation from a mother
with symmetry $\Gam_i$ to a daughter with symmetry $\Gam_j$,
represented by the {\em bifurcation arrow}
\[
\xymatrix{\Gamma_i \ar[r]^k & \Gamma_j} .
\]
The bifurcation arrows always join isotropy subgroups in $\mathcal{G}$, since
$\mathcal{H}_{i,k}\subset\mathcal{G}$.
For each of these bifurcation arrows,
GAP computes the groups $\Gam_i/\Gam_{i,j}'$, $N_{\Gam_i}(\Gam_j)$,
and $N_{\Gam_i}(\Gam_j)/\Gam_j$
(the normalizer of  $\Gam_j$ in $\Gam_i$, denoted $N_{\Gam_i}(\Gam_j)$,
is the largest subgroup of $\Gam_i$ for which $\Gam_j$ is a normal subgroup).

At a nondegenerate bifurcation of a solution with symmetry $\Gam_i$, the critical eigenspace $E$ is 
an irreducible subspace lying in one of the $V_{\Gam_i}^{(k)}$ and we say that the mother undergoes
a bifurcation with $\Gam_i/\Gam_{i,k}'$ symmetry. 
Note that $\Gam_i/\Gam_{i,k}'$ acts freely on $E \subset V^{(k)}_{\Gamma_i}$.
The bifurcation arrows with a given label $k$
represent the symmetries $\Gam_j$ of the daughters that are expected to bifurcate from the mother
when $E \subset V^{(k)}_{\Gamma_i}$.
If the system (\ref{pde}) is restricted to $\fix(\Gam_j, \R^n)$ then the effective symmetry 
of the bifurcation is $N_{\Gam_i}(\Gam_j)/\Gam_j$.
For example, if $N_{\Gam_i}(\Gam_j)/\Gam_j \cong \Z_2$, then there is a pitchfork bifurcation 
creating two conjugate daughter branches.
For details, see \cite{GSS, NSS2}.

We say that two arrows $\xymatrix{\Gamma_i \ar[r]^k & \Gamma_{j_1}}$ and $\xymatrix{\Gamma_i \ar[r]^k & \Gamma_{j_2}}$
are {\em equivalent} if  $\Gam_{j_1}$ and $\Gam_{j_2}$ are conjugate.
Since we only seek non-conjugate solutions in our continuation solver,
one bifurcation arrow from each equivalence class,
together with its auxiliary information,
is written into a file by GAP.
This file is used by our continuation solver when a bifurcation is encountered,
as described below.

The amount of material contained in the bifurcation arrows is overwhelming.
To summarize the collected information about the possible bifurcations we draw a bifurcation
digraph (see \cite{NSS2}, which gives an equivalent definition).
This is an extension of the usual lattice of isotropy subgroups.
For an example, see Figure~\ref{P3_digraph_condd}.
Our continuation solver requires the label of the irreducible representation $k$,
but does not require information about the group $\Gam_i/\Gam_{i,k}'$;
on the other hand, humans find
the symmetry group of the bifurcation
more informative than $k$, so
it is included in the bifurcation digraph instead of $k$.

\begin{definition}
\label{digraphdef}
The {\em bifurcation digraph} of the $\Gamz$ action on a real vector space $V=\R^n$
is a directed graph with labeled arrows between the symmetry types.
We draw an arrow from $[\Gam_i]$ to $[\Gam_j]$ if and only if $\Gam_j$ is conjugate to 
a maximal isotropy subgroup 
of the $\Gam_i$ action on some isotypic component $V^{(k)}_{\Gam_i}$.
The label on this arrow is $\Gam_i/\Gam_{i,k}'$, where $\Gam_{i,k}'$ is the kernel
of the $\Gam_i$ action on $V^{(k)}_{\Gam_i}$.
We use the arrow types
\begin{align}
\mbox{solid } \xymatrix{\ar[r]&}
&\mbox{ if }
N_{\Gam_i}(\Gamma_j)/\Gamma_j \cong \Z_2, \nonumber \\
\mbox{dashed } \xymatrix{\ar@{-->}[r]&}
&\mbox{ if }
N_{\Gam_i}(\Gamma_j)/\Gamma_j 
                      \cong \Z_1, \mbox{ and} \nonumber\\
\mbox{dotted }\xymatrix{\ar@{.>}[r]&}
&\mbox{ otherwise, } \nonumber
\end{align}
to indicate the nature of the bifurcation.
\end{definition}

Assume that $E \subset V^{(k)}_{\Gam_i}$ is an irreducible subspace, and a critical eigenspace of a bifurcation
point with symmetry $\Gam_i$.
The group $N_{\Gam_i}(\Gamma_j)/\Gamma_j$ acts freely on $\fix(\Gam_j, E)$.
The size of this factor group, which determines the arrow type in the bifurcation digraph, is passed to the continuation solver.
The arrow type gives us information about the dimension of $\fix(\Gam_j, E)$.
For real vector spaces
$V = \R^n$, the solid and dashed arrows imply that
$\dim_\R \fix(\Gamma_j, E) = 1$, whereas the dotted arrows give
$\dim_\R \fix(\Gamma_j, E) > 1$.
The first two arrow types, with  1-dimensional fixed point spaces, are called EBL bifurcations
since the Equivariant Branching Lemma  \cite{GSS} guarantees (under certain conditions) that solutions
with symmetry $\Gam_j$ are born at the bifurcation.
These bifurcating branches are particularly easy to follow numerically, since there is only one critical
eigenvector (up to a scalar multiple) with the symmetry $\Gam_j$.
The EBL can be considered an extension of the classic bifurcation results from~\cite{Rabinowitz}.

If there is a dotted arrow to $\Gamma_j$
and certain nondegeneracy conditions hold, then there is a daughter with symmetry $\Gam_j$ born
for gradient systems \cite{GSS}.
No general theory predicts where the daughters lie when projected to $\fix(\Gamma_j, E)$;
our approach to following such branches numerically requires randomly choosing perturbations within $\fix(\Gamma_j, E)$.

The bifurcation digraph is often very complicated so we use
condensation classes instead of conjugacy classes to get a simpler {\em condensed bifurcation digraph}.
See Figure~\ref{P3_digraph_condd} for an example.
Let $\phi\in\aut(\Gamma_0)$. If $\phi(\Gamma_i)$ is a symmetry group for all symmetry groups $\Gamma_i$
then we say that $\phi$ is \emph{symmetry preserving}. The symmetry preserving automorphisms
form a subgroup $\aut_c(\Gamma_0)$ of $\aut(\Gamma_0)$. A symmetry preserving map induces a permutation
of symmetries so $\aut_c(\Gamma_0)$ acts on $\mathcal G$. Conjugate elements of ${\mathcal G}$ are on the 
same orbit since the inner automorphisms of $\Gamma_0$ are in $\aut_c(\Gamma_0)$. Hence $\aut_c(\Gamma_0)$
also acts on $\mathcal S$. The orbit equivalence classes of this latter action
are called \emph{condensation classes},
and are computed automatically by a GAP program.
The condensed bifurcation digraph is the quotient of the bifurcation digraph by the orbit equivalence
of the $\aut_c(\Gam_0)$ action on $\mathcal S$.
Hence, the vertices of the condensed bifurcation digraph are the condensation classes.

\subsubsection{Digraph Layout Code}
Using a C++ program similar to the graph layout code, we generate a layout for the bifurcation
digraph and the condensed bifurcation digraph
The difference is that the vertices can move 
horizontally but their vertical position is determined by the size of the group they represent.
Graphics of the layouts are then created by Gnuplot and \Xy-pic.

\subsection{Basis Generation Code}

\label{basisgen}


To generate $\Psi_m = \{ \psi_1, \ldots, \psi_m \}$, we use a somewhat complicated procedure
that increases the efficiency of the GNGA.
We first pick a single arrow 
$\xymatrix{ \Gamma_0 \ar[r]^k & \Gamma_{j_k} }$ for each $k \in K_{\Gam_0}$.
Using Mathematica, we then find a basis of 
eigenvectors $D^{(k)}$ of the symmetric operator $L$ restricted to the invariant subspace
$V_{\Gamma_0}^{(k)}\cap V_{\Gamma_{j_k}}^{(1)}$. 
We compute the orthonormal basis $\Psi^{(k)}$ of $V_{\Gamma_0}^{(k)}$ using the Gram-Schmidt process on 
$\{ \gamma\cdot v \mid \gamma \in \Gamma_0, v\in D^{(k)}  \}$. 
We start the Gram-Schmidt process with the already orthonormal set $D^{(k)}$ so that these elements
survive in the resulting basis.
The eigenvalues of $L$ and the corresponding eigenvectors in
$\Psi_m=\cup_k \Psi^{(k)}$ are written to a file.

The basis generation code also calculates the projection operator
$P_{\Gamma_i}^{(k)}$ for each $i$ and $k \in K_{\Gam_i}$
from the characters and permutations produced by our automorphism group and symmetry analysis codes.
The isotypic component $V^{(k)}_{\Gamma_i}$ is the range of $P_{\Gamma_i}^{(k)}$.
For each $i$ and $k$ the coordinates
in $\Psi_m$
of the elements of $B_{\Gamma_i}^{(k)}$ are written
to a file.

\end{section}

\begin{section}{Newton's Method with Constraints}
\label{newton_w_c}

To follow branches or find new bifurcating branches, we treat the parameter $s$ as the $(m+1)^{\rm st}$ unknown.
Thus, when we say $p=(a,s)\in\R^{m+1}$ is a solution, we mean that $u=\sum a_j\psi_j$ solves Equation~(\ref{pde}) with parameter $s$.
We restrict the search for a particular solution to some hypersurface in $\R^{m+1}$,
satisfying an $(m+1)^{\rm st}$ equation of the form $\kappa(a,s)=0$.
In this section we will describe two different choices of $\kappa$,
one for following branches and another for finding new branches emanating from bifurcation points.
In the first case, we take an {\em old} and {\em current} pair
of solutions  $p_o$ and $p_c$ along a symmetry invariant branch
to obtain
a reasonable initial guess $p_g$ for iterating to find a {\em new} solution $p_n$ satisfying the constraint
of lying on a hyperplane normal to the branch.
The constraint used at a bifurcation point $p^*$ instead forces the new solution $p_n$ to have a projection
onto the critical eigenspace of a specified norm.

In either case, the iteration we use is:

\vbox{
\begin{itemize}
\item compute the constraint $\kappa$, gradient vector $g:=g_s(u)$, and Hessian matrix $h:=h_s(u)$
\item solve $\left[\begin{array}{cc}
h & \frac{\partial g}{\partial s} \\
(\nabla_a\kappa)^T & \frac{\partial\kappa}{\partial s}
\end{array}\right]
\left[\begin{array}{c} \chi_a \\ \chi_s \end{array}\right]
=
\left[\begin{array}{c} g \\ \kappa \end{array}\right]
$
\item $(a,s)\leftarrow (a,s) - \chi$, \ \ $u=\sum a_j \psi_j$.
\end{itemize}
}
\noindent
Equations~(\ref{grad}) and (\ref{hess}) are used to compute $g$ and $h$.
The $(m+1)^{\rm st}$ row of the matrix is defined by
$(\nabla_a\kappa, \frac{\partial\kappa}{\partial s})=\nabla\kappa\in\R^{m+1}$;
the search direction is
$\chi=(\chi_a,\chi_s)\in\R^{m+1}$.
Since this is Newton's method on $(g,\kappa)\in\R^{m+1}$ instead of just $g\in\R^m$,
when the process converges we have not only that $g=0$
(hence $p=(a,s)$ is a solution to Equation~(\ref{pde})),
but also that $\kappa=0$.

\subsection{Tangent-augmented Newton's method (tGNGA)}
We use the tGNGA to follow branches.
Given two consecutive solutions $p_o$ and $p_c$ on a given branch, we compute the
(approximate) tangent vector $v=(p_c - p_o)/ \|p_c - p_o\|\in\R^{m+1}$.
The initial guess is then $p_g=p_c +cv$.
In our experiments the speed $c$ has a minimum and maximum range, for example from 0.01 to 0.4,
and is modified dynamically according to various heuristics
(see for example Figure \ref{points}).
For the tGNGA, the constraint is that each iterate $p=(a,s)$ must lie on
the hyperplane passing through the initial guess $p_g$, perpendicular to $v$.
That is, $\kappa(a,s) := (p-p_g)\cdot v$.
Easily, one sees that $(\nabla_a\kappa(a,s), \frac{\partial\kappa}{\partial s}(a,s)) = v$.
In general, if $f_s$ has the form $f_s(u) = s u + H(u)$, then $\frac{\partial g}{\partial s} = -a$.
For example, when $f_s(u) = su + u^3$, a calculation shows that
$g_s(a)_j=a_j(\lambda_j-s)-(\sum_{k=1}^m a_k\psi_k)^3\cdot\psi_j$,
hence $\frac{\partial g}{\partial s} = -a$.
Newton's method is invariant in this plane so that in fact $\chi\cdot v=0$ at each step.
Hence, the linear system to be solved each iteration can be described by:
$$
\left[\begin{array}{cc}
h & \frac{\partial g}{\partial s} \\
(v_a)^T & v_s
\end{array}\right]
\left[\begin{array}{c} \chi_a \\ \chi_s \end{array}\right]
=
\left[\begin{array}{c} g \\ 0 \end{array}\right],
$$
where $v=(v_a,v_s)\in\R^{m+1}$.
Our function {\tt tGNGA}$(p_g, v)$
returns, if successful,
a new solution $p_n$ satisfying the constraint.

\subsection{The secant method and processing bifurcation points}

\label{secant}

In brief,
when using the tGNGA to follow a solution branch and the
MI changes at consecutively found solutions,
say from $k$ at the solution $p_o$
to $k+d$ at the solution $p_c$,
we know by the continuity of $D^2J_s$ that there exists a third, nearby solution $p^*$ where
$h$ is not invertible and the $r^{\rm th}$ eigenvalue of $h$ is zero,
where
$r=k+ \lceil\frac{d}2\rceil$.
Let $p_0 = p_o$, $p_1 = p_c$, with $\beta_0$ and $\beta_1$
the $r^{\rm th}$ eigenvalues of $h$ at the points $p_0$ and $p_1$, respectively.

We effectively employ the vector secant method by iterating

\vbox{
\begin{itemize}
\item
$\displaystyle{ p_g = p_i - \frac{(p_i-p_{i-1})\beta_i}{(\beta_{i}-\beta_{i-1})} }$
\item
$ p_{i+1} = {\tt tGNGA}(p_g, v) $
\end{itemize}
}
\noindent
until the sequence $(p_i)$ converges.
The vector $v=(p_c - p_o)/\|p_c - p_o\|$ is held fixed throughout,
while the value $\beta_i$ is the newly computed $r^{\rm th}$ eigenvalue
of $h$ at $p_i$.
If our function {\tt secant}$(p_o, p_c)$ is successful, it returns
a solution point $p^*=(a^*,s^*)$,
lying between $p_o$ and $p_c$,
where $h$ has $r$ zero eigenvalues within some tolerance.
We take the critical eigenspace $E$ to be the span of the corresponding eigenvectors.
If $p^*$ is not a turning point, then it is a bifurcation point.

\subsection{Cylinder-augmented Newton's method (cGNGA)}

The cGNGA is used to find initial solution points on new branches near bifurcation points
$p^*$.
After such a point has been detected and the corresponding critical eigenspace
has been computed, we search for a new solution bifurcating from the main branch by
running the cGNGA.
The first step is to choose a subspace $E$ of the critical eigenspace.

To ensure that we find the mother solution rather than a daughter,
we insist that the Newton iterates belong to the cylinder
$C :=\{(a, s) \in\R^{n+1}:\|P_E (a - a^*)\|=\eps\}$,
where $P_E$ is the orthogonal projection onto $E$
and the radius $\eps$ is a small fixed parameter.
The input file to the continuation solver sets the value of $\eps$.
At a symmetry breaking bifurcation the critical eigenspace is orthogonal to
the fixed point subspace of the mother,
so the mother branch does not intersect the cylinder.
We conjecture that
at anomaly breaking bifurcations
the critical eigenspace is orthogonal to
the AIS of the mother branch.
It is not true that the critical eigenspace is orthogonal to the mother in general,
although we observed this to be true for all the numerical results we include in this paper.
Consider
Equation~(\ref{pde}) with $f_s(u) = (s u + u^3)(u^2 -1)$,
which has a bifurcation where two solution branches in ${\mathcal A}_c$ cross
at $u = (1, 1, \ldots, 1)$, $s = -1$. 
At this bifurcation point $a^*$ and the critical eigenvector are parallel.

The constraint we use to put each Newton iterate on the cylinder is
$\kappa(a,s)=\frac12(\|P_E (a - a^*)\|^2-\eps^2)=0$.
The initial guess we use is $p_g:= (a^*, s^*) +\eps (e,0)$, where $e$ is a
randomly chosen unit vector in $E$.  Clearly, $p_g$ lies on the cylinder $C$.
A computation shows that 
$\nabla_a \kappa (a,s) = P_E(a-a^*)$, and $\frac{\partial\kappa}{\partial s}(a,s) = 0$.
Hence, the search direction $\chi$ is found by solving
$$
\left[\begin{array}{cc}
h & \frac{\partial g}{\partial s} \\
\left(P_E(a-a^*)\right)^T & 0
\end{array}\right]
\left[\begin{array}{c} \chi_a \\ \chi_s \end{array}\right]
=
\left[\begin{array}{c} g \\ {\kappa} \end{array}\right].
$$
Again, $\frac{\partial g}{\partial s} = -a$ when $f_s$ has the form $f_s(u) = s u + H(u)$.
When successful, {\tt cGNGA($p^*, p_g, E$)}
returns a new solution $p_n$ of Equation~(\ref{pde}) that lies on the cylinder $C$.

We take $E$ to be various low-dimensional subspaces of the critical eigenspace,
corresponding to the possible symmetries that bifurcations theory predicts must exist.
For example, at an EBL bifurcation $E$ is spanned by a single eigenvector.
When the dimension of $E$ is greater than one, we call {\tt cGNGA} repeatedly with
several random choices of the critical eigenvector $e$.
The details are given in Equation~(\ref{Ej}) and Algorithm~\ref{find_daughters}.
The theory we apply does not guarantee a complete prediction of all daughter solutions.
Therefore we also call {\tt cGNGA} with $E$ equal to the full critical eigenspace.
In this way, if the dimension of the critical eigenspace is not too big
we have a high degree of confidence that we are capturing all relevant solutions,
including those that arise due to accidental degeneracy and that are neither predicted nor ruled out by understood bifurcation theory.

\end{section}

\begin{section}{Continuation Solver}
\label{solver}

Our continuation solver is implemented in C++.
We start our search for solution branches
with the known solution $p_c=(0,s_0)\in\R^{m+1}$, which lies on the trivial branch,
together with the initial direction vector $v=(0,\pm1)\in\R^{m+1}$,
which points in the direction of another solution on the trivial branch.
The branch queue is initialized with the job $(p_c, v)$. 
Every job in the branch queue is fed to
{\tt follow\_branch}
until the branch exits some window in $\R^{m+1}$.
After every new point is computed,
{\tt find\_bifpoints} is called.
If a bifurcation point is found,
{\tt find\_daughters} is called and
a job is added to the branch queue for every solution found.
For efficiency, 
{\tt find\_daughters} only returns solutions on
distinct, non-conjugate bifurcating branches.
The process stops when the branch queue is empty.
Thus, our continuation solver finds
a representative branch from each conjugacy class of branches connected to the trivial branch, within the chosen window.

\subsection{Branch Following}

Once an instance of {\tt follow\_branch} has started, consecutive solutions $p_o$ and $p_c$ are used to generate 
the next direction vector 
$v=\frac{p_c-p_o}{\|p_c-p_o\|}$.
The \emph{speed} $c$ is modified according to the scale and complexity 
of features, e.g., severe turning points, proliferation of proximal bifurcation points, or failure of the algorithm to converge. 
When the algorithm converges especially quickly, for example, we use other heuristics to increase the speed as our 
guesses are somehow too good.  In that way, many solution points are found near trouble spots while 
much fewer are needed on long, 
featureless parts of a branch (see for example Figure~\ref{points}).
Each point along the branch together
with its corresponding data is written to a file, to be used later in generating bifurcation
diagrams and reports, as well as for diagnosing the occasional failure.
Algorithm~\ref{follow_branch} is executed repeatedly until the branch queue is empty.

\begin{Algorithm}
\begin{enumerate}
\item[] {\bf while} ($p_c\in {\tt window}$ and $c>\tau$)
\begin{enumerate}
\item[] set $p_g = p_c + cv$
\item[] set $p_n$ = {\tt tGNGA}($p_g$, $v$)
\item[] {\bf if} $p_n$ is unacceptable 
\begin{enumerate}
\item[] halve speed $c$
\end{enumerate}
\item[] {\bf else}
\begin{enumerate}
\item[] use heuristics to adjust speed $c$ 
\item[] set $p_o = p_c$
\item[] set $p_c = p_n$
\item[] set $v = \frac{p_c - p_o}{\|p_c - p_o \|}$
\item[] {\bf forall} $p_i^* \in $ {\tt find\_bifpoints}($p_o$, $p_c$)

\begin{enumerate}
\item[] compute critical eigenspace $E_i$ of $p_i^*$
\item[] {\bf forall} $q_i^j \in {\tt find\_daughters}(p_i^*, E_i)$
\item[] \hskip .4in set $v_i^j = \frac{p_i^* - q_i^j}{\|p_i^* - q_i^j\|}$
\item[] \hskip .4in add $(p_i^*,v_i^j)$ to {\tt branch\_queue}
\end{enumerate}
\end{enumerate}
\end{enumerate}
\end{enumerate}
\caption{{\tt follow\_branch}($p_c$, $v$)}
\label{follow_branch}
\end{Algorithm}

\begin{Algorithm}
\begin{enumerate}
\item[] {\bf if} $|\MI(p_o)-\MI(p_c)|=0$ 
  \begin{enumerate}
  \item[] {\bf return} $\{ \, \}$
  \end{enumerate}
\item[] set $p={\tt secant}(p_o, p_c)$
\item[] compute critical eigenspace $E$ of $p$
\item[] {\bf if} $\dim(E)=|\MI(p_o)-\MI(p_c)|$ 
\begin{enumerate}
  \item[] {\bf if} $p$ is not a fold point
  \begin{enumerate}
  \item[] {\bf return} $\{ p \}$
  \end{enumerate}	    
\end{enumerate}
\item[] {\bf else}
  \begin{enumerate}
  \item[] set $p_g=(p_o+p_c)/2$
  \item[] set $v=\frac{p_c-p_o}{\|p_c-p_o\|}$
	\item[] set $p={\tt tGNGA}(p_g,v$) 
	\item[] {\bf return} ${\tt find\_bifpoints}(p_o, p) \cup {\tt find\_bifpoints}(p, p_c)$
  \end{enumerate}
\end{enumerate}
\caption{{\tt find\_bifpoints}$(p_o, p_c)$}
\label{find_bifpts}
\end{Algorithm}

The implementations of the tGNGA and secant methods are straightforward following
Section~\ref{newton_w_c}.
Details concerning {\tt find\_bifpoints} and {\tt find\_daughters} can be found in
Algorithms~\ref{find_bifpts} and \ref{find_daughters}, respectively.

\begin{figure}
\begin{center}
\input{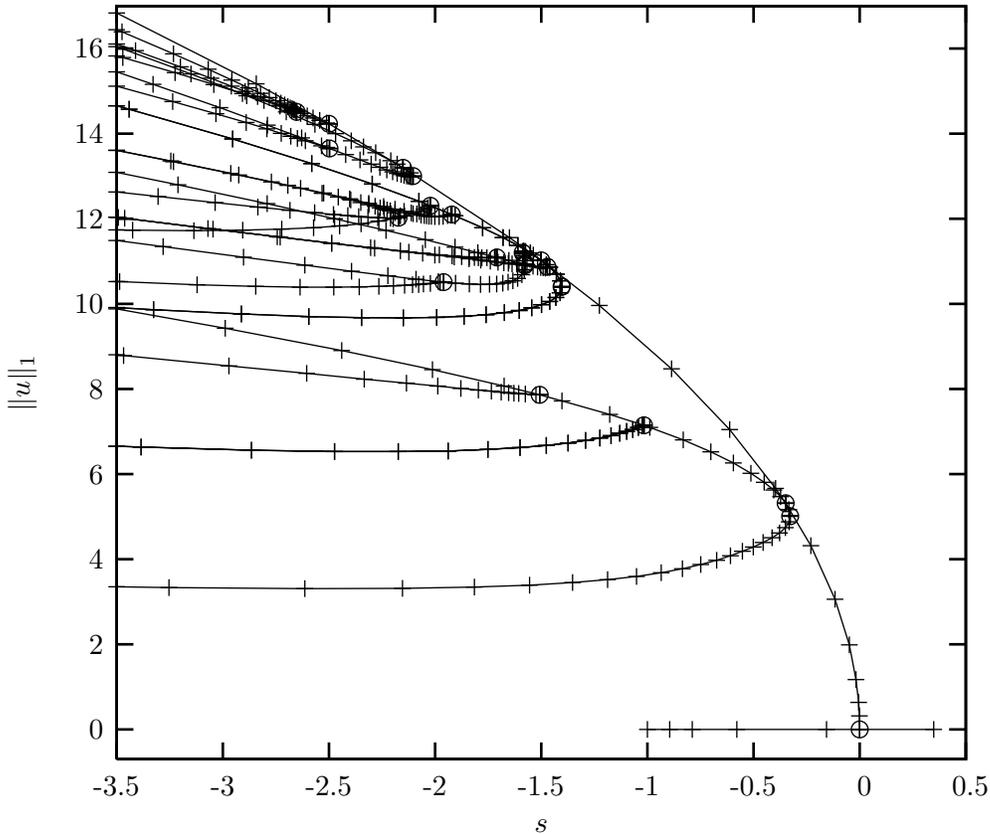}
\caption{Bifurcation diagram for the first primary branch of the Cayley graph of $\s_3$ 
(see Section~\ref{cayleys3}).
The graphic demonstrates how
the density of computed points (denoted by crosses) is increased near interesting features.
We use heuristics to adjust the speed.
For example, the speed $c$ is halved if tGNGA fails to converge in four iterations.
Further, 
the speed is multiplied by a factor in $(0, 2]$ based on the angle formed by the
last three points, 
where the factor is 1 if the angle is 0.1 radians.
}
\label{points}
\end{center}
\end{figure}

\subsection{Finding bifurcation points}

When a MI change is observed between two consecutive solutions $p_o$ and $p_c$  
on a branch, we know
that there exists at least one degenerate intermediate solution $p^*$ on the branch.
Since there can be multiple degenerate points on the branch in the branch segment,
our bifurcation point finding algorithm recursively calls the secant method of
Section~\ref{secant} (see Algorithm~\ref{find_bifpts}).

\subsection{Finding daughter branches.}
Given a bifurcation point $p^* = (u^*,s^*)$ with symmetry
$\Gamma_i$ and its critical eigenspace $E$, 
we want to find all of the bifurcating branches.
Our function {\tt find\_daughters},
described in Algorithm~\ref{find_daughters},
finds as many daughters as it can.

If $p^*$ is nondegenerate, then $E$ is an irreducible subspace so it is contained in exactly one isotypic component
$V_{\Gamma_i}^{(k)}$
of the $\Gam_i$ action on $V = \R^n$ and has dimension $d_{\Gamma_i}^{(k)}$.
The first step of {\tt find\_daughters}
is to find the intersection of $E$ with each isotypic component.
These determine the set
$\textsf{K}=\{ k \mid  E \cap V_{\Gamma_i}^{(k)} \neq \{ 0 \}\}$. 
Here we use the set of bases $\{B_{\Gam_i}^{(k)} \mid k \in K_{\Gam_i}\}$ computed in Section~\ref{basisgen}.

A bifurcation point is degenerate if $E$ is contained in the fixed
point subspace
$\fix(\Gam_i, V)$
of the mother.
When $u^*=0$, $E$ cannot be contained in the zero-dimensional fixed point subspace
$\fix(\Gam_0, V) = \{0\}$
of the mother.
In all other cases, we label
the irreducible representations of $\Gam_i$ so that the trivial representation is $k = 1$,
and $\fix(\Gam_i, V) = V_{\Gamma_i}^{(1)}$.

\begin{definition}
\label{accidental}
The bifurcation point $p^*$ with symmetry $\Gam_i$
has \emph{accidental degeneracy} if any of the following conditions hold:
\begin{enumerate}
\item \textsf{K} contains 1 and $\Gam_i \neq \Gam_0$;
\item \textsf{K} is not a singleton set;
\item $\dim(E \cap V_{\Gamma_i}^{(k)})> d_{\Gamma_i}^{(k)}$ for some $k \in \textsf{K}$.
\end{enumerate}
Then we say that $p^*$ has a \emph{degeneracy of Type 1, 2 or 3} respectively. 
\end{definition}

Condition (1) says that $E$ has a nontrivial intersection with the fixed point
subspace of the mother.
If $E$ is not an irreducible subspace then condition (2) or (3) holds.

\begin{Algorithm}
\begin{enumerate}
\item[] compute \textsf{K} and \textsf{J} for $\Gam_i$ and $E$
\item[] compute $E_j$ for $j \in  \textsf{J}$
\item[] set ${\tt list\_of\_ daughters} = \{ \, \}$
\item[] {\bf if} $\dim(E) > 1$
\begin{enumerate}
\item[] set $E_0 = E$
\item[] set $\textsf{J} = \textsf{J} \, \cup \{0\}$
\end{enumerate}
\item[] {\bf forall} $j\in \textsf{J}$
\begin{enumerate}
\item[] set {\tt num\_no\_changes} $ = f_{nc}(\dim(E_j))$
\item[] set {\tt no\_changes} $ = 0$
\item[] {\bf while} ${\tt no\_changes} < {\tt num\_no\_changes}$
\begin{enumerate}
\item[] choose a random $e \in E_j$ with $\| e \| = \eps$
\item[] set $p_g = (u^* + e, s^*)$
\item[] set $q={\tt cGNGA}(p^*, p_g ,E_j)$
\item[] {\bf if} the $\Gam_i$ orbit of $q$ and {\tt list\_of\_daughters} are disjoint
\begin{enumerate}
\item[] set {\tt no\_changes} $ = 0$
\item[] add $q$ to {\tt list\_of\_daughters}
\item[] set $p_g = (u^* - e, s^*)$
\item[] set $q={\tt cGNGA}(p^*, p_g ,E_j)$
\item[] {\bf if} the $\Gam_i$ orbit of $q$ and {\tt list\_of\_daughters} are disjoint
\item[] \hskip .4in add $q$ to {\tt list\_of\_daughters}
\end{enumerate}
\item[] {\bf else}
\begin{enumerate}
\item[] increment {\tt no\_changes}
\end{enumerate}
\end{enumerate}
\end{enumerate}
\item[] {\bf return} {\tt list\_of\_daughters}
\end{enumerate}
\caption{\tt find\_daughters($p^*, E$)}
\label{find_daughters}
\end{Algorithm}

The set of expected bifurcating symmetry indices of the daughter solutions is
\[
\textsf{J}:=\{ j \mid  \xymatrix{\Gamma_i \ar[r]^{k} & \Gamma_j} 
  \text{ is a bifurcation arrow for some } k\in \textsf{K} \text{ or }
(j = i \neq 0 \text{ and } 1 \in \textsf{K} )\}.
\]
For each $j \in \textsf{J}$, we define
\begin{equation}
\label{Ej} 
E_j = E \cap \fix(\Gamma_j, V). 
\end{equation}
For nondegenerate bifurcations,
$\sym(u^* + e) = \Gam_j$ for all nonzero $e \in E_j$,
since $\Gam_j$ is a maximal isotropy subgroup of the $\Gam_i$ action on $E$.
Our algorithm uses initial guesses $u^* + e$ with $\| e \| = \eps$ in the
cylinder-augmented Newton's method function $\tt{cGNGA}$ to
look for solutions with symmetry $\Gam_j$.

The daughters are found by repeated applications of $\tt{cGNGA}$.
A heuristic function $f_{nc}: \N \rightarrow \N$ with $f_{nc}(1) = 1$
takes as input the dimension of $E_j$ or $E$, and outputs the number
of  $\tt{cGNGA}$ consecutive calls allowed without finding a new daughter.
The default function is defined by $f_{nc}(d) = 1 + 20(d-1)^2$, but
this can be changed if one suspects that a daughter branch has not yet been found.

The {\tt find\_daughters} subroutine prints information such as the number of random choices
it took to find a new daughter, so that the user can modify the $f_{nc}$ function if desired.

The search for solutions in all of $E$ when $\dim(E) > 1$ is 
included to find possible daughter branches
with symmetry not predicted by any bifurcation arrows.
This is needed at bifurcation points with Type 2 degeneracy, as in Example~\ref{cayleys3}.
Daughters with submaximal symmetry exist for
bifurcations with certain symmetries \cite{FieldRichardson}
even when $E$ is an irreducible subspace,
although we did not encounter this in the examples we studied.
\end{section}

\begin{section}{Postprocessing}
\label{postprocess}

Since even a small graph can have a large symmetry group and other features which
lead to a proliferation of solutions via possibly complicated bifurcations,
we must artfully display select subsets of our results in a human understandable format.
This section briefly describes the methodology and tools we have developed which
process the output from the preprocessing and continuation solver phases in order to
generate graphics automatically, edit and annotate those graphics, 
and research new variational and symmetry phenomena.

Our heuristics for automatically changing speed and retrying Newton's method with
better initial guesses are sufficient to generate 
all the results presented in Section~\ref{examples}, and many more.
In a few instances some adjustment of the initial speed was required.  
This adjustment is facilitated by files which track every solution and present 
all associated information in human readable formats.  
Generally, these files also contain the actual data that the programs described in Sections~\ref{contoursub} and \ref{schematic} 
use to generate graphical output.

\subsection{Contour plots}
\label{contoursub}
Solutions in $X$ are displayed
with a contour plot program written in Mathematica.
The contour plot program uses
the embedding of the graph found by the layout program described in Section~\ref{layout}.
The vertex $v_i$ is colored white if $u_i > 0$, gray if $u_i = 0$, and black if $u_i < 0$.
Furthermore, the vertex is shown as a disk whose area is proportional to $|u_i|$.
If $|u_i|$ is below some cutoff, a small disk is drawn.
When a solution $u \in \R^n$ is passed to the contour plot program,
each of the solutions in the group orbit
$\{\gamma \cdot u \mid \gamma \in \aut(G) \}$ are tested by heuristics
that attempt to find which one is the best.  
We plot the solution $u$ that minimizes the size of the set
$\{ u_i u_j \|d_{ij}\| \mid 1 \leq i < j \leq n \}$.
To break a tie, the program chooses a solution which has
a horizontal or vertical line of reflection symmetry, if such a solution exists.

We say a symmetry of a solution is {\it visible} if it is also a symmetry of the contour plot.
By changing several parameters the layout program can easily be made to generate alternate layouts
which may make more symmetry of a given solution visible.
Layouts can also be entered by hand or copied from the
output of other programs.  Once we have viewed a solution's contour plot for a particular layout,
we can view and save any solution in the orbit of that solution. 
Saved graphics are most easily viewed in an automatically
created HTML file which annotates each representative solution with useful information such as
Morse index, symmetry and symmetry type, $J$ value, branch number, and bifurcation history.
In these ways, we greatly reduce the human effort needed to generate informative graphics
in a format suitable for publication.

\subsection{Bifurcation diagrams}
\label{schematic}

A \emph{(schematic) bifurcation diagram} is the graph of $\{ (s,y(u)) \mid (u,s)\in X \}$ 
where $y:\R^n\to \R$ is some \emph{schematic function} \cite{GSS}. 
The schematic function $y$ is needed to reduce the graphics to two dimensions.
A good choice such as the taxicab norm
$y(u) = \|u\|_1 = |u_1| + |u_2| + \cdots + |u_n|$
visually separates branches. 
The $\Gamma_0$-invariance of this choice ensures that only one curve
is shown for each equivalence class of solution branches. 
Thus, it avoids apparent discontinuities 
in the diagrams due to inconsistent choices of representatives of orbit classes 
for bifurcating branches. 

In \cite{NSS2} we used the value of a PDE solution $u$ at a generic point of the domain. 
We actually solved a PdE, and 
defined $y$
by $y(u)=u_i$ for a fixed $i$, where the vertex $v_i$ has trivial symmetry
as described in Proposition~\ref{all_sym_prop}.
Note that for graphs in general (for example cycles),
there may not be a generic vertex.
While this choice of $y$ is not a $\Gamma_0$-invariant function, 
we were able to get meaningful bifurcation
diagrams without redundant branches by exploiting the simplicity of the symmetry group $\D_6$ 
in a way that cannot be done for general groups \cite{NSS}.  

In the current project,
we also investigated schematic functions of the form $y(u)=y_w(u)=\sum_{i=1}^n w_i|u_i|$, 
for some choice of \emph{weight vector} $w$.
We only include results using 
$w=(1,\ldots,1)$, which makes $y(u)=\|u\|_1$,
but find this topic a interesting area
for future research.
\end{section}

\begin{section}{Examples}
\label{examples}

We considered many different graphs in our numerical experimentation, with an eye for
examples that revealed interesting phenomena in symmetry, bifurcation, or variational structure.
Among other things,
we want to know which of the possible symmetries are represented in the solution space $X$,
how the symmetry of solutions relate to the symmetries of eigenfunctions of the Laplacian,
and what are the relationships between Morse index and nodal structure.
These questions were first raised for PdE in \cite{Neu}.
We also investigate several examples of anomaly breaking bifurcations.
The chosen examples demonstrate capabilities such as our ability to
handle high multiplicity bifurcations and accidental degeneracies. 
In general, the output automatically generated during our experiments was sufficient for the 
creation of the graphics and tables included in this section.
The following is an index of the experiments that we have decided to include in this section.
\begin{enumerate}
\item[\ref{path_sub}:] The path $P_3$.
We demonstrate the continuation solver for five different nonlinearities, not all odd nor all superlinear. 
Our code works without modification when $f_s$ is not odd.
We discuss the branch of constant solutions present in all our experiments.
An accidental degeneracy of Type 1 is featured.

\item[\ref{cycle_sub}:]  The cycle $C_4$.
Branches connected to the trivial branch and the existence of solutions of every possible symmetry are discussed.

\item[\ref{nosym_sub}:] Graphs with no symmetry.
The smallest graphs with no symmetry have 6 vertices, and there are 9 such graphs.
We show results for the two graphs that have AIS other than ${\mathcal A}_c$.
These AIS are associated with integer eigenvalues of $L$
and lead to anomaly-breaking bifurcations.

\item[\ref{cayleys3}:] A decorated Cayley graph of the symmetric group $\s_3$.
We demonstrate that we can automatically generate a graph and solutions to Equation~(\ref{pde}) on that graph with
a predetermined symmetry group.  We highlight an accidental degeneracy of Type 2 at an integer eigenvalue.

\item[\ref{Cayley_Z5}:] A decorated Cayley graph of $\Z_5$.  This example has several non-EBL bifurcations.
The critical eigenspace for the 
bifurcations with $\Z_{10}$ symmetry can be in either of two different isotypic components;
the same is true for the bifurcations with $\Z_5$ symmetry.
We show contour plots of the bifurcating solutions that occur in the different cases.

\item[\ref{Qexample}:] A decorated Cayley graph of the quaternion group $Q$.
We construct an example with several occurrences of bifurcations with $Q$ symmetry.
At these non-EBL bifurcations,
all points in the 4-dimensional critical eigenspace (except for the origin) have the same symmetry.

\item[\ref{petersensub}:] The Petersen graph.
Some information concerning the number of Newton iterations and overall computing time is presented for
this fairly complicated example which has bifurcation points of high multiplicity.  We find a numerical
counterexample to a variational and nodal structure conjecture.

\item[\ref{dodeca_sub}:] The dodecahedron.  In this example we encounter
an accidental degeneracy of Type 3. 
The degeneracy is explained by an AIS that relates certain solutions on the dodecahedron
to solutions on the Petersen graph.

\item[\ref{soccerball_sub}:]  The truncated icosahedron (soccer ball).
This example has more vertices and a large number of
high multiplicity eigenvalues.
We choose a layout that visually shows the resemblance of several solutions
to spherical harmonics.

\end{enumerate}

\begin{subsection}{The path $P_3$}
\label{path_sub}

\begin{table}
\begin{center}
\begin{tabular}{|c|c|l|c|}
\hline
$\Gamma_i\cong$ & $[\Gamma_i]$ & $\Gamma_i$ & Contour Plot of Solution in $\fix(\Gamma_i)$\\
\hline  \hline
$\Z_2\times\Z_2$ & $S_0$ & $\Gamma_0 = \langle \alpha_{(13)}, -1 \rangle$ & 
 \scalebox{.3}{\includegraphics{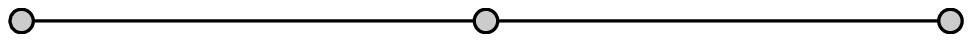}}   \\ \hline
$\Z_2$           & $S_1$ & $\Gamma_1=\langle \alpha_{(13)} \rangle$ &
 \scalebox{.3}{\includegraphics{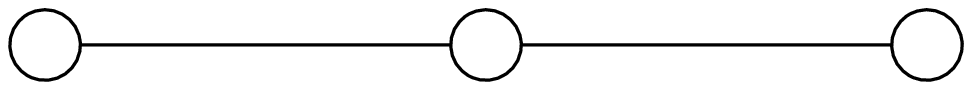}}   $\quad$
 \scalebox{.3}{\includegraphics{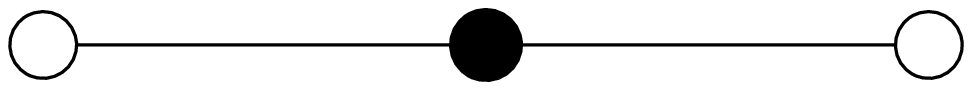}}
\\ \cline{2-4}
		 & $S_2$ & $\Gamma_2=\langle -\alpha_{(13)} \rangle$ &
\scalebox{.3}{\includegraphics{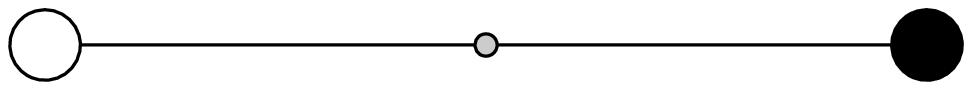}}  \\ \hline
$\Z_1$ 		 & $S_3$ & $\Gamma_3=\langle 1 \rangle$ &
\scalebox{.3}{\includegraphics{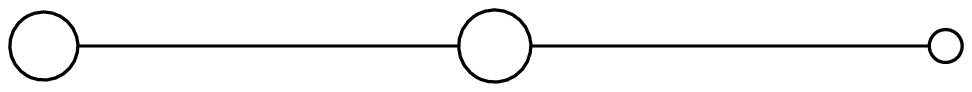}}  \\ \hline
\end{tabular}
\end{center}
\medskip

\caption
{
Symmetries for the path $P_3$ (Example~\ref{path_sub}).
The first column shows the isomorphism class of the elements in a
condensation class; the second and third columns give the symmetry type and the symmetry.
The fourth column shows contour plots for selected solutions with each symmetry type.
Two solutions with symmetry $\Gamma_1$ are shown, the constant solution on the left
is in the AIS ${\mathcal A}_c$.
}
\label{P3.symmetries}
\end{table}

In this very simple example one can easily see the entire bifurcation digraph and possible symmetries of solutions.
It is an easy exercise to increase the number of vertices in the path and consider scaling in order to approximate solutions to an ODE
with Neumann boundary conditions.

Let $G=P_3$ be the path with three vertices. Then $\aut(G)=\langle \alpha_{(1\,3)} \rangle\cong\Z_2$,
and so $\Gamma_0\cong\Z_2\times\Z_2$. There are four symmetries in $\mathcal G$, shown 
in Table~\ref{P3.symmetries}.
The symmetry types are all singletons, with $S_i = [\Gam_i] = \{ \Gamma_i \}$.

The bifurcation digraph and condensed bifurcation 
digraph (see Section~\ref{sym_analy}) are shown in Figure~\ref{P3_digraph_condd}. 
The automorphism group $\aut(\Gamma_0)$ is isomorphic to $\Z_2$, and is generated by $\phi$,
where $\phi(\alpha_{(1\,3)})=-\alpha_{(1\,3)}$ and
$\phi(-1)=-1$.
Thus, $\phi$ interchanges $\Gamma_1$ and $\Gamma_2$, while leaving $\Gamma_0$ and $\Gamma_3$ fixed.
There are three condensation classes, as seen in Figure~\ref{P3_digraph_condd}.

Figure~\ref{P3_bif_diagram} shows the bifurcation diagrams for several nonlinearities
that can be chosen by a flag in our continuation solver;
our implementation handles non-odd nonlinearities with no modification, provided $f_s(0) = 0$.
In addition to the functions shown in Figure~\ref{P3_bif_diagram}, our code works
for many other families of functions.
For example, we can use $f_s(u) = \sinh(s u)$, which is not of the form
$f_s(u) = s u + H(u)$ since it
has a nonlinear dependence on $s$. We also performed experiments with
asymptotically linear nonlinearities.
The solutions which bifurcate from $u = 0$, $s = 0$ in Figure~\ref{P3_bif_diagram}
are all constant solutions of the form $u=(c,\ldots,c)$.
These constant solution branches satisfy  $f_s(c) = 0$.
For example, in the first diagram, with  $f_s(u) = s u + u^3$, the constant solution has $c = \sqrt{-s}$.
The Hessian evaluated at a constant solution
is a diagonal matrix with $h_{ii} = \lam_i - f_s'(c)$.
It is an exercise to determine
the values of $s$ at which the Hessian is singular on this branch.

As noted in Section~\ref{functSym},
we find that
$W_c = {\mathcal A}_c$
is an AIS
since $P_3$ is not vertex transitive.  
In the first diagram in Figure~\ref{P3_bif_diagram}, the constant branch has a bifurcation
at $s = -1.5$ to a daughter $u \not \in {\mathcal A}_c$.  
This is {\em not} a symmetry-breaking
bifurcation since both mother and daughter have symmetry type $S_1$ 
(see Table~\ref{P3.symmetries}).
The critical eigenspace lies in the fixed-point subspace of the mother, hence
this bifurcation point has accidental degeneracy of Type 1 
(see Definition~\ref{accidental}).

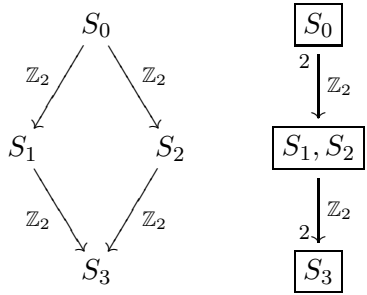
\begin{figure} 
\scalebox{1}{
\xymatrix@C=10pt{
&
S_0 \ar[dl]_<>(.5){\Z_2} \ar[dr]^(.5){\Z_2}
&
&
&
\fbox{$S_0$} \ar[d]^(.5){\Z_2}_<>(.1)2
\\
S_1 \ar[dr]_<>(.5){\Z_2}
&
&
S_2 \ar[dl]^(.5){\Z_2}
&
&
\fbox{$S_1, S_2$} \ar[d]^(.5){\Z_2}_<>(.85)2
\\
&
S_3
&
&
&
\fbox{$S_3$}
}
}
\caption{
The bifurcation digraph (left) and condensed bifurcation digraph (right)
for the $\Z_2 \times \Z_2$ action on $P_3$ (Example~\ref{path_sub}).
The elements in each condensation class are enclosed in a box.
The small numerals on the arrows tell the
number of connections emanating from each symmetry type in a box.
A missing small numeral means 1.
}
\label{P3_digraph_condd}
\end{figure}

\begin{figure}
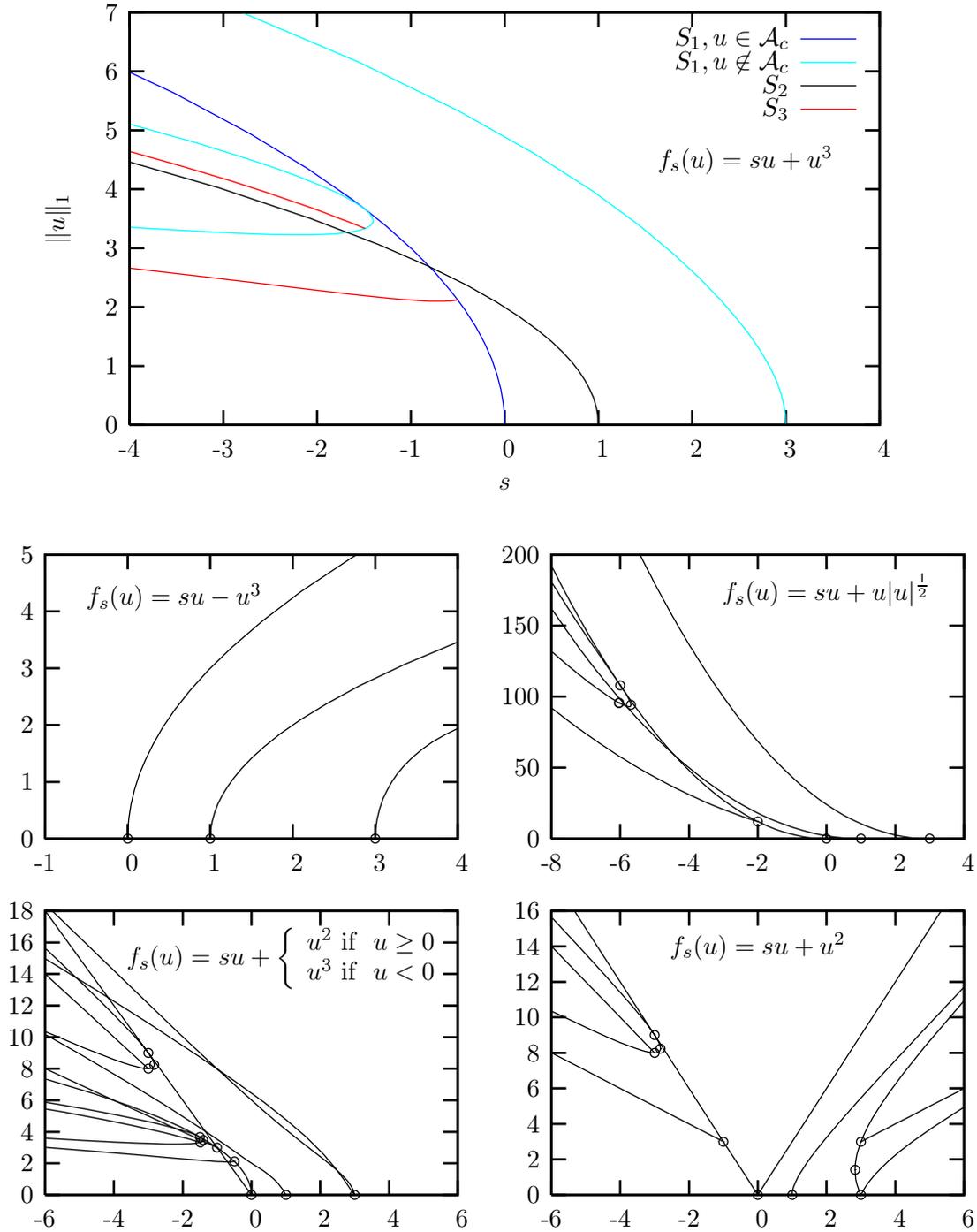

\begin{center}
\input P3_f0_bif.tex \\
\vspace{.2in}

\input P3_f1_bif.tex 
\input P3_f3_bif.tex \\

\input P3_f2_bif.tex 
\input P3_f4_bif.tex \\

\end{center}
\bigskip
\caption
{
Bifurcation diagrams for the graph $P_3$ with various nonlinearities (Example~\ref{path_sub}).  The first diagram is for our standard
odd, superlinear nonlinearity.  All of the diagrams use the taxicab norm $\|u\|_1$ plotted against $s$.
Note that the nonlinearities featured in the bottom row are not odd, but our procedures still work.
Extending results from \cite{Neu}, one computes that the secondary bifurcations
on the constant branch for the five cases respectively are at: 
$s=-\frac{\lambda_i}2$; nonexistent; $s=-2\lambda_i$; $s=-\frac{\lambda_i}2$ for $c <0$ and $s=-\lambda_i$ for $c > 0$;
$s=-\lambda_i$ for $c>0$ and nonexistent for $c < 0$.
}
\label{P3_bif_diagram}
\end{figure}

\end{subsection}

\begin{subsection}{The cycle $C_4$}
\label{cycle_sub}

We investigated far too many families of graphs to include them all, but make a brief mention of the cycle 
$C_4$ due to two interesting phenomena that we observed.
Figure~\ref{C4.bifurcation.diagram} shows every branch that is connected to the trivial branch.
Lee and Neuberger \cite{Lee} found one additional branch for $C_4$
that is not connected to the trivial branch, hence is missing from Figure~\ref{C4.bifurcation.diagram}.  
The cubic system~(\ref{pde}) with the default nonlinearity has at most $3^4 = 81$ real solutions.
Lee and Neuberger found exactly 81 real solutions for $s < s^* \approx -3$
by using their asymptotic form of solutions for large, negative $s$.

Also, notice that there is no branch of symmetry type $S_8$ in our figure.
With the vertices numbered cyclically,
functions of symmetry type $S_8$ are of the form $u=(a, b, -a, -b)$, with $a\not=b$ nonzero real numbers.
The additional branch found in \cite{Lee} does not have symmetry type $S_8$ either,
which provides strong evidence that some systems do not realize all possible symmetry types.

\begin{figure}
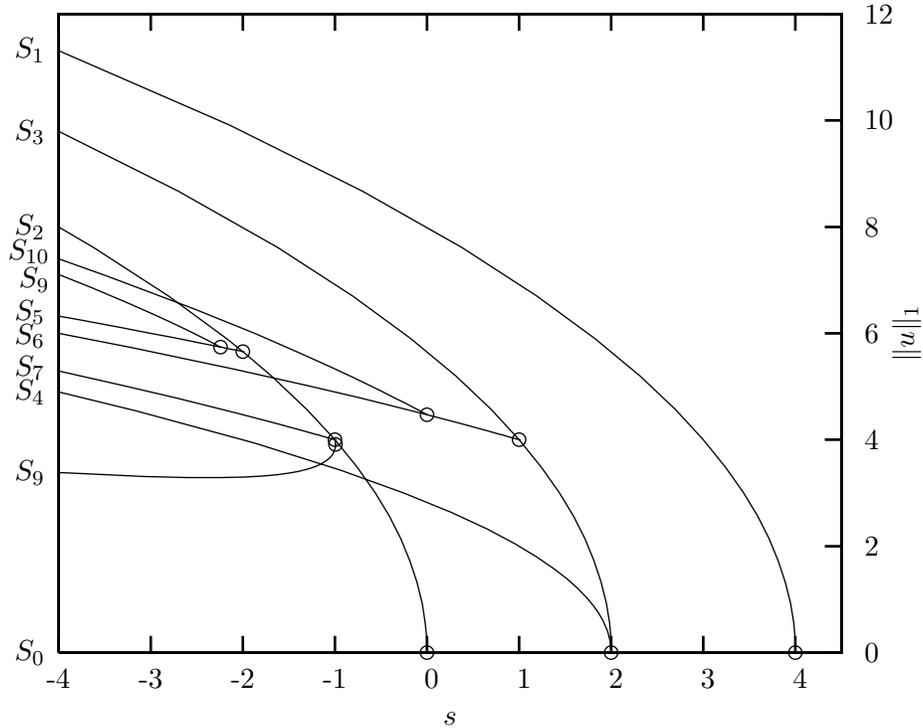

\begin{center}
\input C4.tex
\caption{Bifurcation diagram for $C_4$ (Example~\ref{cycle_sub}), showing all branches connected to the trivial branch.
The symmetry type of each branch is indicated.
There are 11 possible symmetry types for this system, but no solutions
with symmetry type $S_8$ are found.
One more branch (not connected to the trivial branch) is found 
in \cite{Lee}, but this branch does not have symmetry type $S_8$ either.
The output files automatically generated by our suite of programs for this graph, including this bifurcation diagram,
can be viewed at
{\tt http://NAU.edu/Jim.Swift/PdE}.
}
\label{C4.bifurcation.diagram}
\end{center}
\end{figure}

\end{subsection}

\begin{subsection}{Graphs with no symmetry}
\label{nosym_sub}

We considered graphs with no symmetry.  In this case, the set of possible symmetries of functions is
${\mathcal G} = \{ \Gam_0, \Gam_1\}$, where $\Gam_0 = \Z_2$ and $\Gam_1 = \{ 1 \}$.
The trivial solution has symmetry $\Gam_0$ and all other solutions have trivial symmetry.
Thus, nontrivial solutions cannot undergo symmetry-breaking bifurcations. 
The ``expected'' behavior based solely on symmetry theory is
to have $\Z_2$ bifurcations at the eigenvalues of $L$, which are typically simple, with no secondary bifurcations. 
As in Example~\ref{path_sub}, there are anomaly-breaking bifurcations of constant solutions in ${\mathcal A}_c$.
In this section we describe another AIS that is present for 
some graphs with trivial symmetry.

We have done automated experiments computing the automorphism groups of all graphs with 6 vertices or fewer.
Other than the graph with one vertex and no edges, all graphs $G$ with 5 or fewer vertices had nontrivial $\aut(G)$.
There are exactly 9 graphs with 6 vertices and $\aut(G) \cong \{1\}$.
Two of these graphs have AIS other than ${\mathcal A}_c$.
Non-constant anomalous solutions to Equation~(\ref{pde}) for these two graphs are shown in Figure~\ref{nosyms}.
For both of these graphs, the vertices can be numbered so that the AIS is
$$
{\mathcal A}_2 = \{ (a, a, b, b, b, b) \mid a \in \R, \ b \in \R \}.
$$
It is noteworthy that every AIS for these 9 graphs with no symmetry contains an eigenvector of
$L$ with an integer eigenvalue, and every eigenvector of $L$ with an integer eigenvalue is contained in
an AIS.  For example, $(1,1,1,1,1,1) \in {\mathcal A}_c$ is an eigenvector of $L$ with eigenvalue 0 for any graph,
and $(2, 2, -1, -1, -1, -1) \in {\mathcal A}_2$ is an eigenvector of $L$ with eigenvalue 3 for the two graphs shown in Figure~\ref{nosyms}.

Figure~\ref{nosyms} also shows the bifurcation diagram for one of the graphs that has a non-constant AIS.
The bifurcation diagram for the other graph is similar.
We observe secondary bifurcations on the two primary branches bifurcating at the integer eigenvalues 0 and 3.
The secondary branch born of the constant solution at $s = -3/2$ contains solutions that are in ${\mathcal A}_2$, and
these solutions have tertiary bifurcations to non-anomalous solutions.

Finally, in Figure~\ref{nosyms} one sees the phenomena of branch grouping by MI as $s\to-\infty$. 
We have observed this ``grouping by MI'' in all our experiments which use
the schematic $y(u)=\|u\|_1$, a superlinear $f$, and sufficiently negative $s$.
The reason for the grouping is largely explained in \cite{Lee}, where the asymptotic form of solutions in this realm 
takes on values in $\{ 0, c_s, -c_s\}$ at each vertex,
where $f_s(c_s) = 0$.
The MI is computed by counting the number of nonzero components in the solution vector $u$, which also directly accounts for 
the value $y(u)=\|u\|_1$.

\begin{figure}
\begin{center}
\begin{tabular}{cc}
\includegraphics[scale = 0.55]{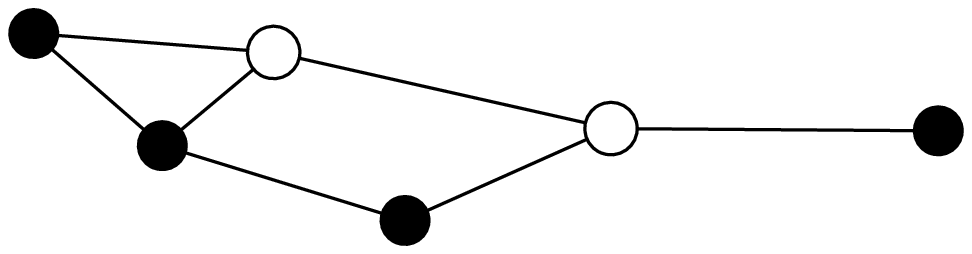} \quad & \quad
\includegraphics[scale = 0.53]{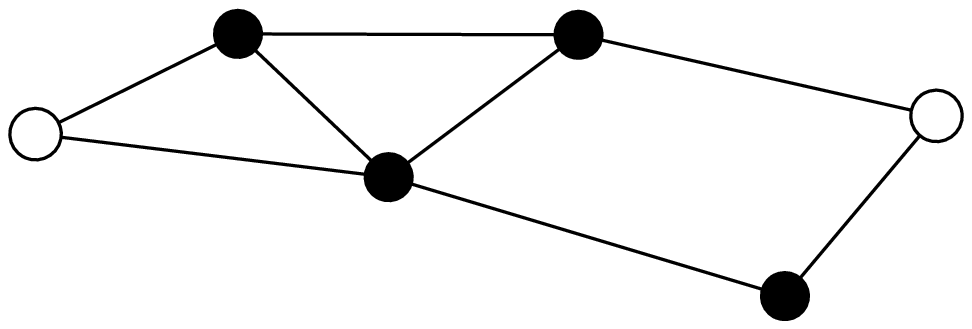} \\
\end{tabular}
\vspace{.15in}

\input nosym_bif.tex
\caption{Sign-changing anomalous solutions on 
two nonsymmetric graphs.
In both cases the solutions
lie on a primary branch bifurcating at $s = 3$.
All secondary and tertiary bifurcations are anomaly-breaking for these two graphs.
The bifurcation diagram for the graph on the left is shown.
}
\label{nosyms}
\end{center}
\end{figure}
\end{subsection}

\begin{subsection}{A decorated Cayley graph of the symmetric group $\s_3$}
\label{cayleys3}

Cayley graphs provide a way for us to generate a graph with a particular symmetry group.  
In Figure~\ref{S3_graphs} we show how 
the Cayley color digraph $\text{Cay}_{\{(1\,2),(2\,3)\}}\s_3$ is used to generate a decorated Cayley graph with $\D_3$ symmetry.
The symmetries and symmetry types for this graph are shown in Table~\ref{S3.symmetries}.
In Figure~\ref{S3_digraph} we show
uncondensed and condensed bifurcation digraphs containing all arrow types and nontrivial conjugacy classes.  
The uncondensed diagram has been annotated with contour plots to give visual cues as to the corresponding symmetries.  
The layouts and contour plots were all automatically generated by our suite of programs.

The matrix $L$ for this graph
has the triple eigenvalue $\lam_4 = \lam_5 = \lam_6 = 3$, whereas the irreducible representations of $\s_3$ are
one or two-dimensional.  Thus, the bifurcation point $(u,s ) = (0, 3)$ has an
accidental degeneracy of Type 2, since the
critical eigenspace $E$ is the direct sum of two irreducible spaces
(see Definition~\ref{accidental}).
Such accidental degeneracy is a common feature of our experiments,
since the matrix entries of the graph Laplacians are integers. 
The bifurcation of the constant branch at $s = -1.5$, seen in Figure~\ref{points},
also has an accidental degeneracy of Type 2.

\begin{figure}
\begin{center}
\scalebox{1}{
\xygraph{ !{0;<.35cm,0cm>:}
[]="O"
"O"
[r(-8.36773)] [u(10.8)]{(\,)}="v1"
"O"
[r(-3.75896)] [u(10.8)]{(1\,2)}="v2"
"O"
[r(-12.2978)] [u(4.5)]{(2\,3)}="v3"
"O"
[r(-10)] [u(0.5)]{(1\,3\,2)}="v5"
"O"
[r(0.146992)] [u(4.5)]{(1\,2\,3)}="v6"
"O"
[r(-2.16289)] [u(0.5)]{(1\,3)}="v7"
"v1":@<2pt>@{.>}"v2"
"v2":@<2pt>@{.>}"v1"
"v3":@<2pt>@{-->}"v1"
"v2":@<2pt>@{-->}"v6"
"v1":@<2pt>@{-->}"v3"
"v5":@<2pt>@{.>}"v3"
"v3":@<2pt>@{.>}"v5"
"v7":@<2pt>@{-->}"v5"
"v6":@<2pt>@{-->}"v2"
"v7":@<2pt>@{.>}"v6"
"v5":@<2pt>@{-->}"v7"
"v6":@<2pt>@{.>}"v7"
}
}
\qquad\qquad
\scalebox{.45}{
\includegraphics{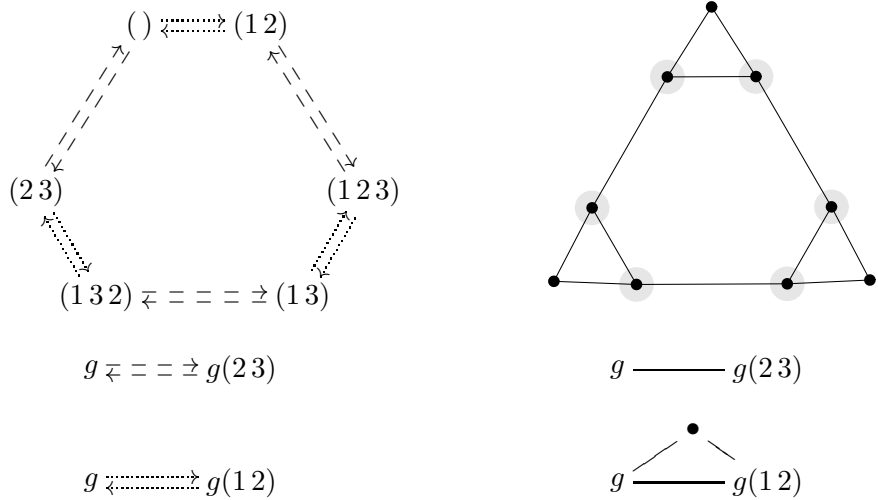}
}
\end{center}
\medskip
\begin{center}
\scalebox{1}{
\xygraph{ !{0;<1.0cm,0cm>:}
[]="O"
"O"
[r(0)] [u(0)]{g}="v0"
"O"
[r(2)] [u(0)]{g(1\,2)}="v1"
"O"
[r(7)] [u(0)]{g}="v2"
"O"
"O"
[r(9)] [u(0)]{g(1\,2)}="v4"
"O"
[r(8)] [u(0.7)]{\bullet}="v5"
"O"
[r(0)] [u(1.5)]{g}="v10"
"O"
[r(2)] [u(1.5)]{g(2\,3)}="v11"
"O"
[r(7)] [u(1.5)]{g}="v12"
"O"
[r(9)] [u(1.5)]{g(2\,3)}="v14"
"v0":@<2pt>@{.>}"v1"
"v1":@<2pt>@{.>}"v0"
"v2"-"v4"
"v2"-"v5"
"v4"-"v5"
"v10":@<2pt>@{-->}"v11"
"v11":@<2pt>@{-->}"v10"
"v12"-"v14"
}
}
\qquad
\end{center}
\caption
{
Cayley graphs of $\s_3$ (Example~\ref{cayleys3}). The graph on the top left is the Cayley color digraph
$\text{Cay}_{\{(1\,2),(2\,3)\}}\s_3$.  The graph on the right is a decorated Cayley graph with $\D_3$ symmetry.
The highlighted vertices of the decorated graph correspond to the vertices of the Cayley color digraph.
The bottom pictures show how the colored directed edges are replaced with decorated undirected edges. Since the
generators are involutions, a pair of directed edges can be replaced by a single edge whose decoration
does not encode the edge direction.
}
\label{S3_graphs}
\end{figure}

\def\psize{.13}
\begin{figure}
\scalebox{1}{
\xygraph{ !{0;<.9cm,0cm>:}
[]="O"
"O"
[r(4.7)] [u(-.1)]{  
	       \rlap{
	       \kern 0.4cm
	       \raise .4cm
	       \hbox{$\scriptstyle S_{6}$}
	       }
               \scalebox{\psize}{
               \includegraphics{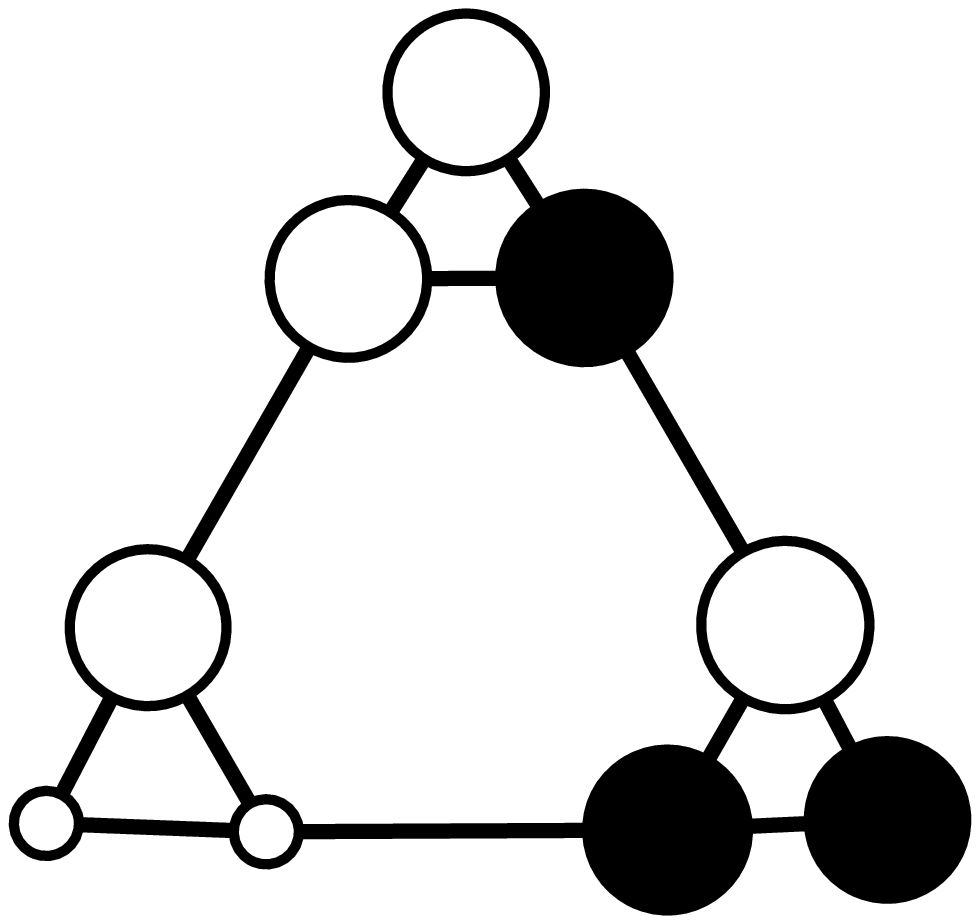} 
	       }
		}="v0"
"O"
[r(7.7)] [u(2)]{  
	       \rlap{
	       \kern 0.4cm
	       \raise .4cm
	       \hbox{$\scriptstyle S_{5}$}
	       }
               \scalebox{\psize}{
               \includegraphics{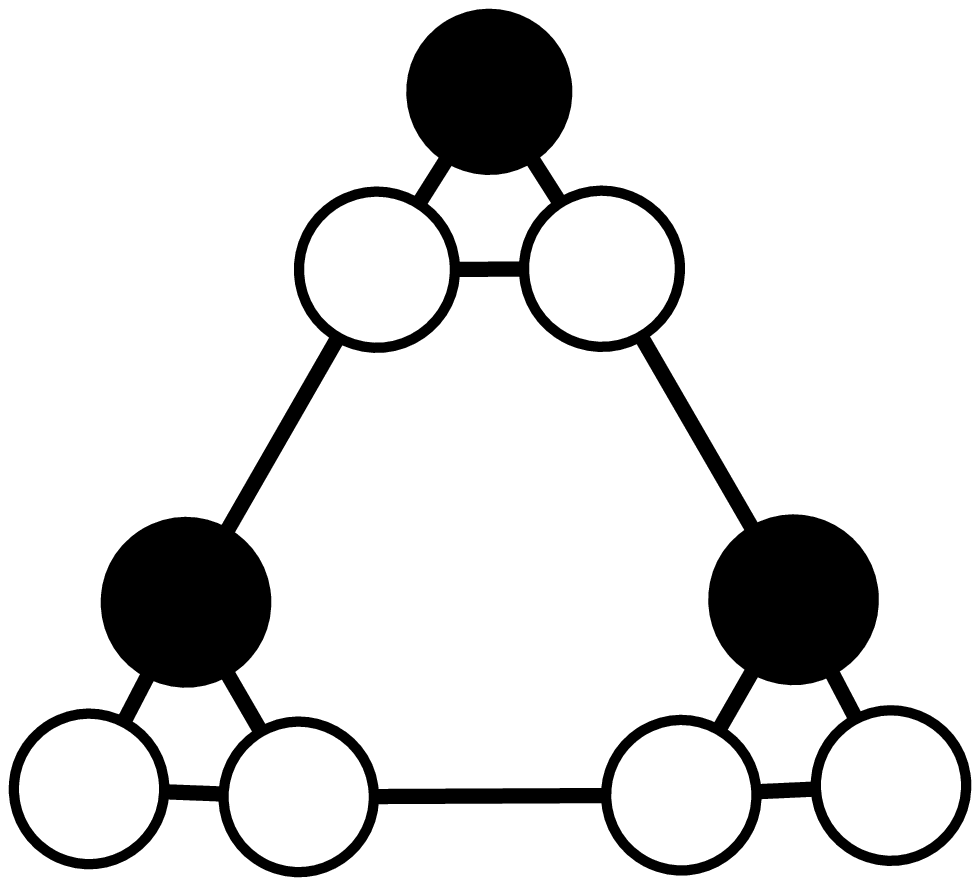} 
	       }
	       }="v1"
"O"
[r(1.7)] [u(2)]{  
	       \rlap{
	       \kern 0.4cm
	       \raise .4cm
	       \hbox{$\scriptstyle S_{4}$}
	       }
               \scalebox{\psize}{
               \includegraphics{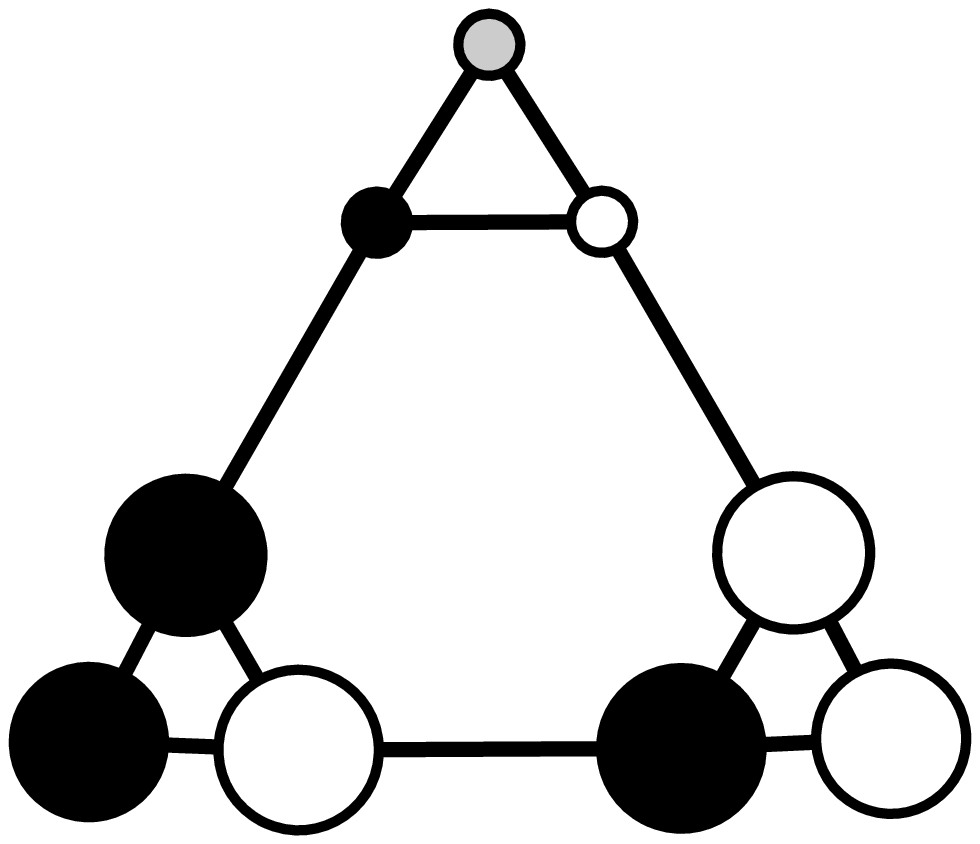} 
	       }
		}="v6"
"O"
[r(4.7)] [u(2)]{  
	       \rlap{
	       \kern 0.4cm
	       \raise .4cm
	       \hbox{$\scriptstyle S_{3}$}
	       }
               \scalebox{\psize}{
               \includegraphics{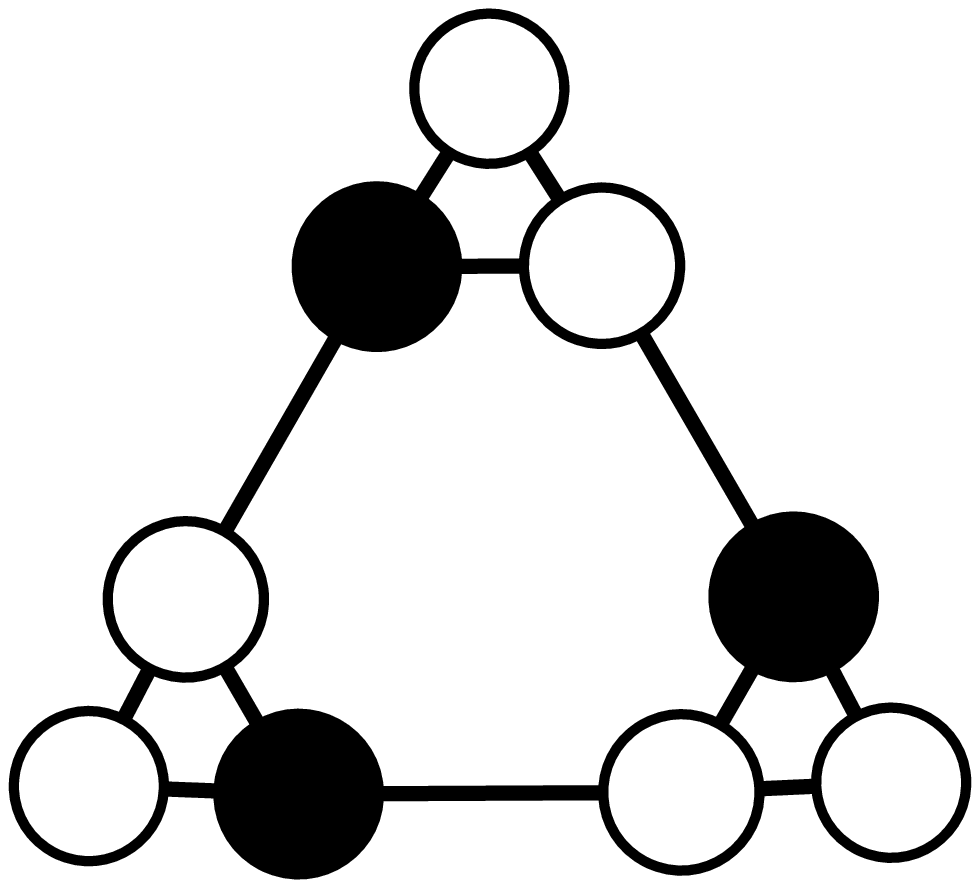} 
	       }
	       }="v2"
"O"
[r(7.7)] [u(4.4)]{ 
	       \rlap{
	       \kern 0.4cm
	       \raise .4cm
	       \hbox{$\scriptstyle S_{2}$}
	       }
               \scalebox{\psize}{
               \includegraphics{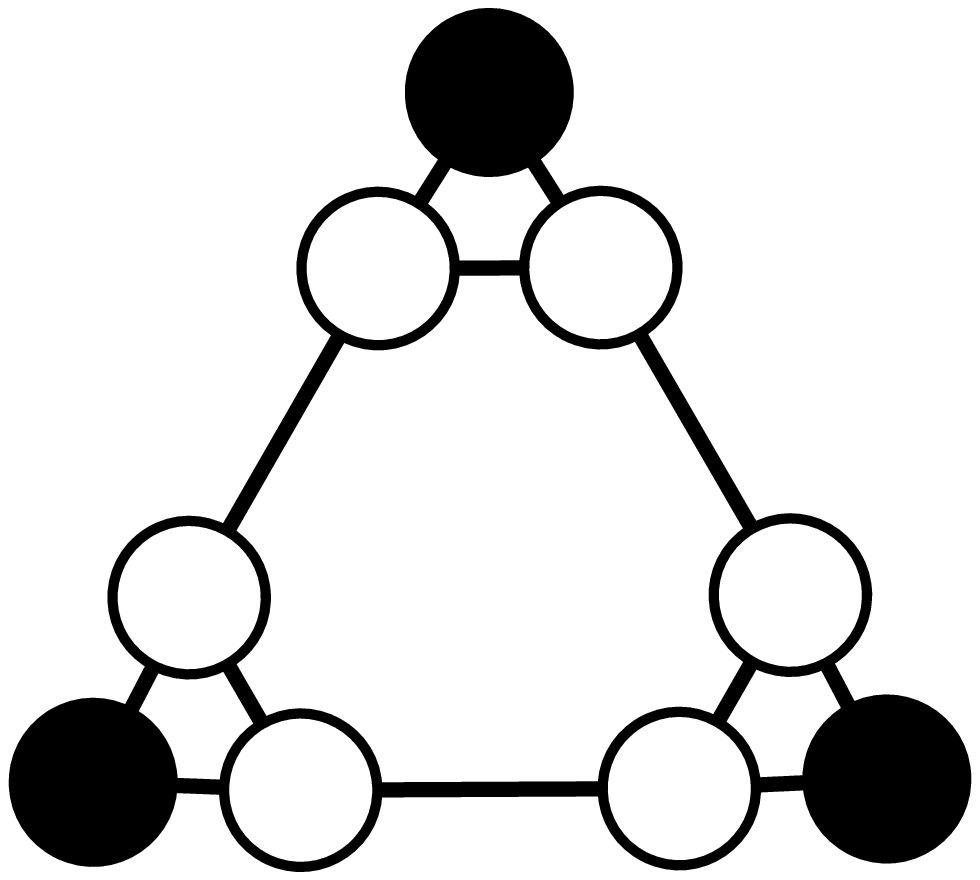} 
	       }
	       }="v3"
"O"
[r(1.7)] [u(4.4)]{
	       \rlap{
	       \kern 0.4cm
	       \raise .4cm
	       \hbox{$\scriptstyle S_{1}$}
	       }
               \scalebox{\psize}{
               \includegraphics{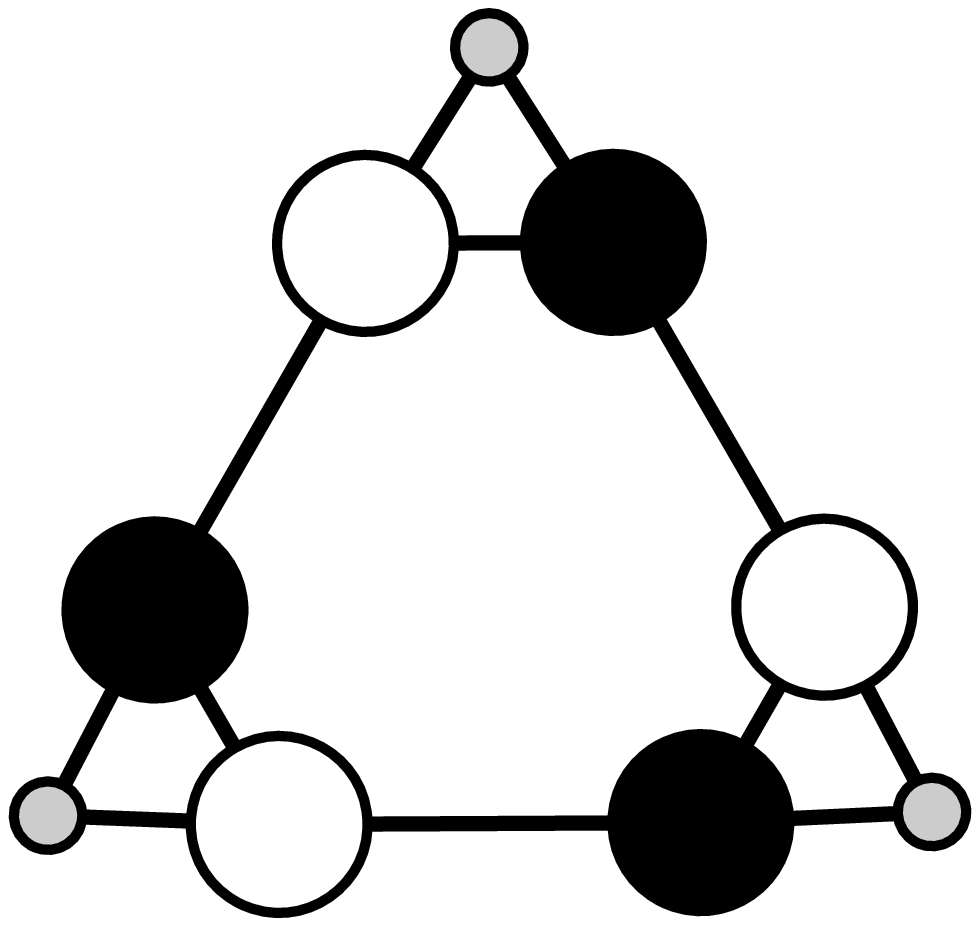} 
	       }
               }="v5"
"O"
[r(4.7)] [u(6)]{  
	       \rlap{
	       \kern 0.4cm
	       \raise .4cm
	       \hbox{$\scriptstyle S_{0}$}
	       }
               \scalebox{\psize}{
               \includegraphics{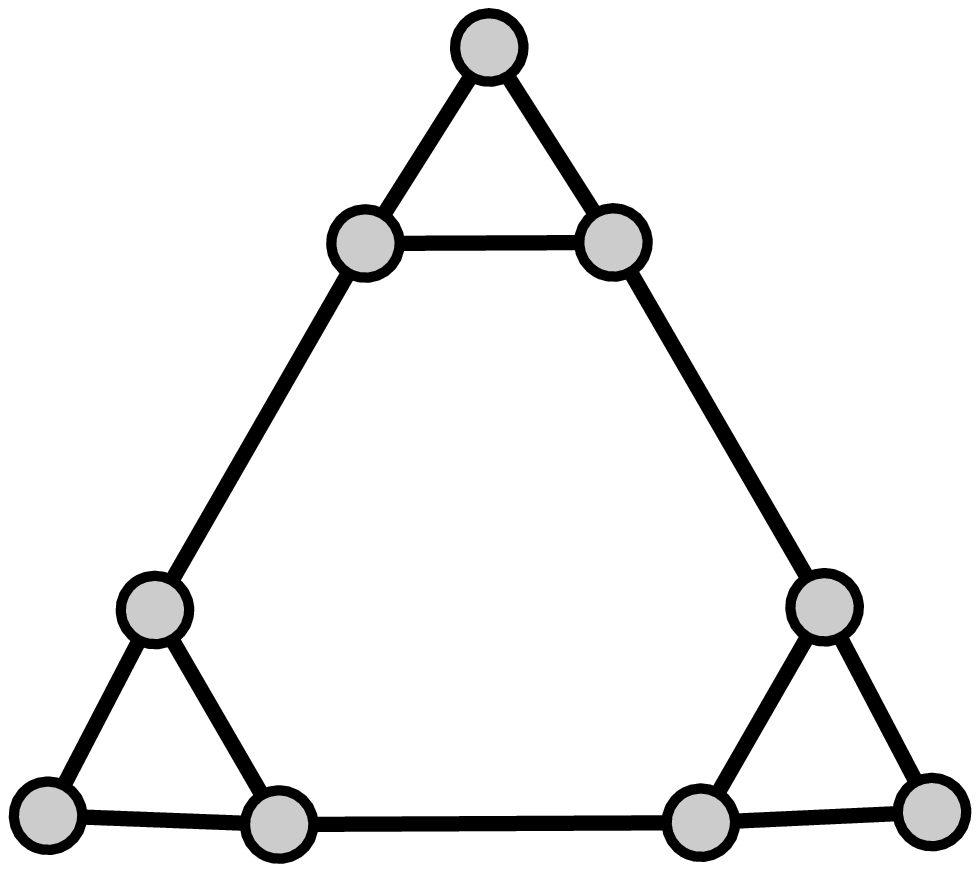} 
	       }
		}="v4"
"v1":^{\Z_2}"v0"
"v6":_{\Z_2}"v0"
"v2":@{.>}_{\Z_3}"v0"
"v3":@{-->}^{\s_3}"v1"
"v5":@{-->}_{\s_3}"v6"
"v3":^{\Z_2}"v2"
"v5":_{\Z_2}"v2"
"v4":_{\D_6}"v1"
"v4":^{\D_6}"v6"
"v4":^{\Z_2}"v3"
"v4":_{\Z_2}"v5"
}
\qquad\quad
}
\scalebox{1}{
\xygraph{ !{0;<.9cm,0cm>:}
[]="O"
"O"
[r(0.7)] [u(0)]{ \fbox{ $  S_{6} $} }="v0"
"O"
[r(3.1)] [u(2)]{  \fbox{ $ S_{4},S_5$} }="v1"
"O"
[r(0.7)] [u(2)]{  \fbox{$ S_{3} $}}="v2"
"O"
[r(0.7)] [u(4)]{ \fbox{$ S_{1},S_2 $}  }="v3"
"O"
[r(3.1)] [u(6)]{  \fbox{$ S_{0} $}  }="v4"
"v1":^{\Z_2}_<>(.9)2"v0"
"v2":@{.>}_{\Z_3}"v0"
"v3":@{-->}_{\s_3}"v1"
"v3":_{\Z_2}^<>(.8)2"v2"
"v4":_{\D_6}^<>(.1)2"v1"
"v4":_{\Z_2}^<>(.1)2"v3"
}
}
\quad \, 
\caption
{
Bifurcation digraphs for a decorated Cayley graph of $\s_3$ (Example~\ref{cayleys3}).
The digraph on the left is not condensed while the digraph on the right is condensed.
}
\label{S3_digraph}
\end{figure}

\begin{table}
\begin{center}
\begin{tabular}{|c|c|l|l|l|}
\hline
$\Gamma_i\cong$ & $[\Gamma_i]$ & \multicolumn{3}{l|}{$\Gamma_i$} \\ \hline  \hline
$\s_3\times\Z_2\cong \D_6$ & $S_0$ & \multicolumn{3}{l|}{$\Gamma_0$}   \\ \hline
$\s_3$           & $S_1$ & \multicolumn{3}{l|}{$\Gamma_1=\langle -(1\,2), -(2\,3) \rangle$}   \\ \cline{2-5}
		 & $S_2$ & \multicolumn{3}{l|}{$\Gamma_2=\langle (1\,2), (2\,3) \rangle$}   \\ \hline
$\Z_3$ 		 & $S_3$ & \multicolumn{3}{l|}{$\Gamma_3=\langle (1\,2\,3)\rangle$}   \\ \hline
$\Z_2$           & $S_4$ & $\Gamma_4=\langle  -(1\,2) \rangle$ &
                           $\Gamma_5=\langle  -(2\,3) \rangle$ & 
			   $\Gamma_6=\langle  -(3\,1) \rangle$   \\ \cline{2-5}
		 & $S_5$ & $\Gamma_7=\langle  (1\,2) \rangle$ &
                           $\Gamma_8=\langle  (2\,3) \rangle$ &  
			   $\Gamma_9=\langle  (3\,1) \rangle$   \\ \hline
$\Z_1$		 & $S_6$ & \multicolumn{3}{l|}{$\Gamma_{10}$ }   \\ \hline
\end{tabular}
\end{center}
\medskip

\caption
{
Symmetries for a decorated Cayley graph of $\s_3$.
The first column shows the isomorphism class of the elements in a
condensation class; the second and third columns give the symmetry type and the symmetry.
}
\label{S3.symmetries}
\end{table}

\end{subsection}

\begin{subsection}{The Cayley graph of $\Z_5$}
\label{Cayley_Z5}

We constructed a decorated Cayley graph of $\Z_5 = \langle a \mid a^5 = 1 \rangle$ with 15 vertices.
We conjecture that this is the smallest graph $G$ with $\aut(G) = \Z_5$.
This example is interesting for two reasons:  First, there are non-EBL
bifurcations with $\Z_5$ and $\Z_5 \times \Z_2 \cong \Z_{10}$ symmetry.
Secondly, there are two inequivalent 2-dimensional irreducible representations of $\Z_5$
with trivial kernels.

Figure~\ref{Z5_digraph} shows the irreducible representations of $\Z_{10}$ over $\R$.
All ten of the irreducible representations of $\Z_{10}$ over $\C$ are one dimensional:
two are real, and eight are complex.
Section~\ref{isoSub} describes how to construct the irreducible representations of $\Z_{10}$ over $\R$.
The three irreducible representations of $\Z_5$ over $\R$ can be obtained by restricting $\alpha^{(k)}$ to $\Z_5$.
As described for general groups in Section~\ref{isoSub},
the irreducible representations of $\Z_{10} \cong \Gam_0$ with $\alpha^{(k)}(-1) = -I$ are listed first.
We choose this ordering so that $V=\bigoplus_{k=1}^3 V_{\Gamma_0}^{(k)}$
and the isotypic components with $k = $ 4, 5, or 6 satisfy $V_{\Gamma_0}^{(k)} = \{0\}$.

It is illuminating to describe the functions in $V_{\Gamma_0}^{(k)}$.
In our layout of the graph, the generator $a$ of $\Z_5$ acts as a rotation by $2\pi/5$.
The vertices lie on three concentric circles, labeled by $j \in \{1, 2, 3\}$. 
Let $\theta_i$ be the angle of vertex $v_i$ in the layout.
Then $V_{\Gamma_0}^{(1)} = \{ u \in V \mid u_i = c_j \text{ when } v_i \text{ lies on circle } j\}$.
Similarly,
$V_{\Gamma_0}^{(2)} = \{ u \in V \mid u_i = c_j \cos(\theta_i) + d_j \sin(\theta_i) \text{ when } v_i \text{ lies on circle } j\}$
and
$V_{\Gamma_0}^{(3)} = \{ u \in V \mid u_i = c_j \cos(2 \theta_i) + d_j \sin(2 \theta_i) \text{ when } v_i \text{ lies on circle } j\}$.
The dimensions of these isotypic components satisfy $3 + 6 + 6 = 15$.

There are exactly three symmetries in $\mathcal G$: $\Gamma_0 \cong \langle a, -1 \rangle = \Z_5 \times \Z_2 \cong \Z_{10}$,
$\Gamma_1 \cong \langle a \rangle = \Z_5$, and $\Gamma_2 \cong \{ 1\} \cong \Z_1$.
The symmetry types are singletons: $S_i = \{ \Gamma_i \}$ for $i \in  \{0, 1, 2\}$.
The bifurcation digraph has exactly three arrows, a solid arrow $S_0 \rightarrow S_1$ with label $\Z_2$, 
a dotted arrow $S_0 \rightarrow S_2$ with label $\Z_{10}$, and 
a dotted arrow $S_1 \rightarrow S_2$ with label $\Z_5$.

While there are three arrows in the bifurcation digraph, there are five bifurcation arrows for this graph,
as shown in Figure~\ref{Z5_digraph}.
For example, the two bifurcation arrows 
$\Gam_0 \rightarrow \Gam_2$, with labels $k = 2$ and $k = 3$ both correspond to 
the single dotted arrow from $S_0 \rightarrow S_2$ with label $\Z_{10}$ in the bifurcation digraph.
The continuation solver needs to know if the critical eigenspace of the origin
is in $V_{\Gamma_0}^{(2)}$ or $V_{\Gamma_0}^{(3)}$,
but from a theoretical point of view there is a bifurcation with $\Z_{10}$ symmetry in either case.
The different nodal structures of the daughters bifurcating from the trivial solution
shown in Figure~\ref{Z5_digraph} are explained by the above descriptions of
$V_{\Gamma_0}^{(k)}$.
When $E \subset V_{\Gamma_0}^{(2)}$,
daughters change sign once, whereas 
when $E \subset V_{\Gamma_0}^{(3)}$,
daughters change sign twice.
Similarly, the perturbations of the constant solution in the bifurcation with $\Z_5$ symmetry 
change sign once and twice, respectively.

\def\psize{.16}
\begin{figure}
\scalebox{1}{
\raise .64 in
\hbox{
\begin{tabular}{c|c|c}
$k$  & $\alpha^{(k)}(a)$ & $\alpha^{(k)}(-1)$  \\
\hline
1 & 1 & $-1$ \\
2 & $R(2\pi/5)$ & $-I_2$ \\
3 & $R(4\pi/5)$ & $-I_2$ \\
\hline
4 & 1 & $1$ \\
5 & $R(2\pi/5)$ & $I_2$ \\
6 & $R(4\pi/5)$ & $I_2$ \\
\end{tabular}
}
\qquad
\scalebox{1}{
\raise 1.2 in
\hbox {
\xymatrix@C=10pt{
S_0 \ar[rd]^(.5){\Z_2} \ar@{.>}[dd]_(.5){\Z_{10}} \\
& S_1 \ar@{.>}[ld]^(.5){\Z_5} \\
S_2
}
}
}
\qquad
\xygraph{ !{0;<.8cm,0cm>:}
[]="O"
"O"
[r(0)] [u(0)]{  
               \rlap{
	       \kern 0.6cm
	       \raise .7cm
	       \hbox{$\scriptstyle \Gamma_{2}$}
	       }
               \scalebox{\psize}{
               \includegraphics{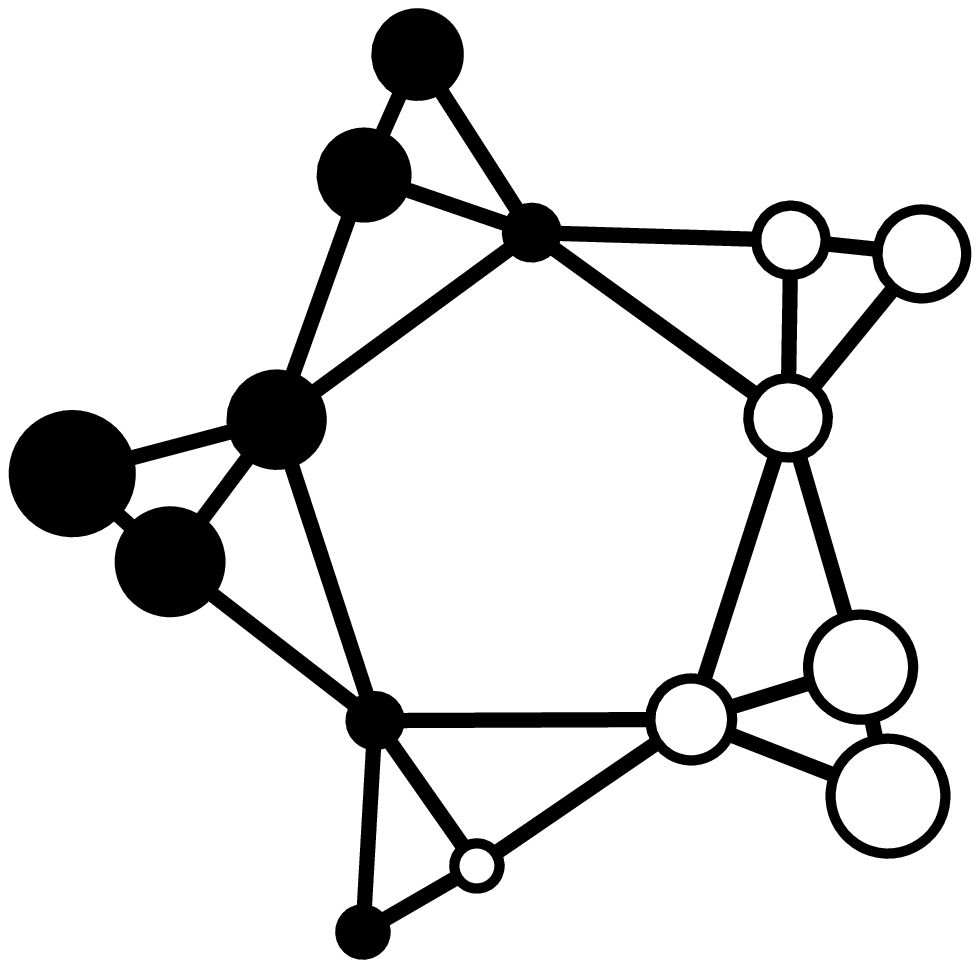} 
	       }
              }="v0a"
"O"
[r(2)] [u(0)]{  
               \rlap{
	       \kern 0.6cm
	       \raise .7cm
	       \hbox{$\scriptstyle \Gamma_{2}$}
	       }
               \scalebox{\psize}{
               \includegraphics{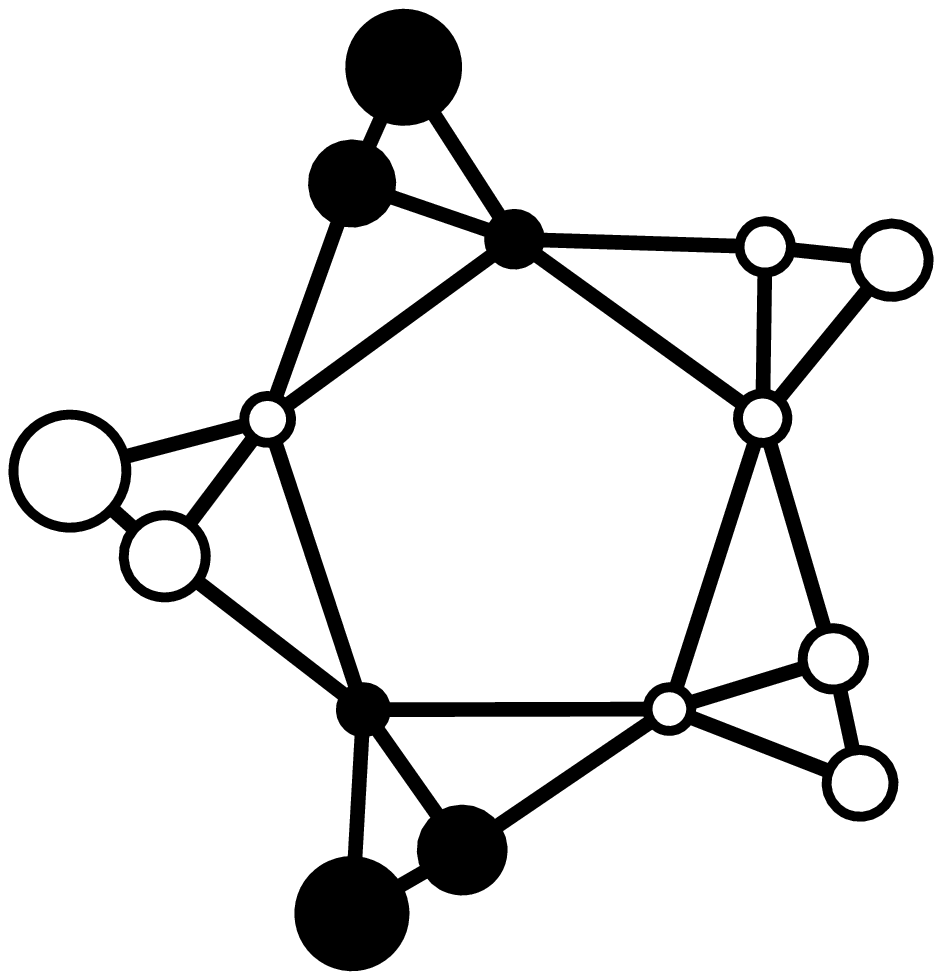} 
	       }
              }="v0b"
"O"
[r(4)] [u(0)]{  
               \rlap{
	       \kern 0.6cm
	       \raise .7cm
	       \hbox{$\scriptstyle \Gamma_{2}$}
	       }
               \scalebox{\psize}{
               \includegraphics{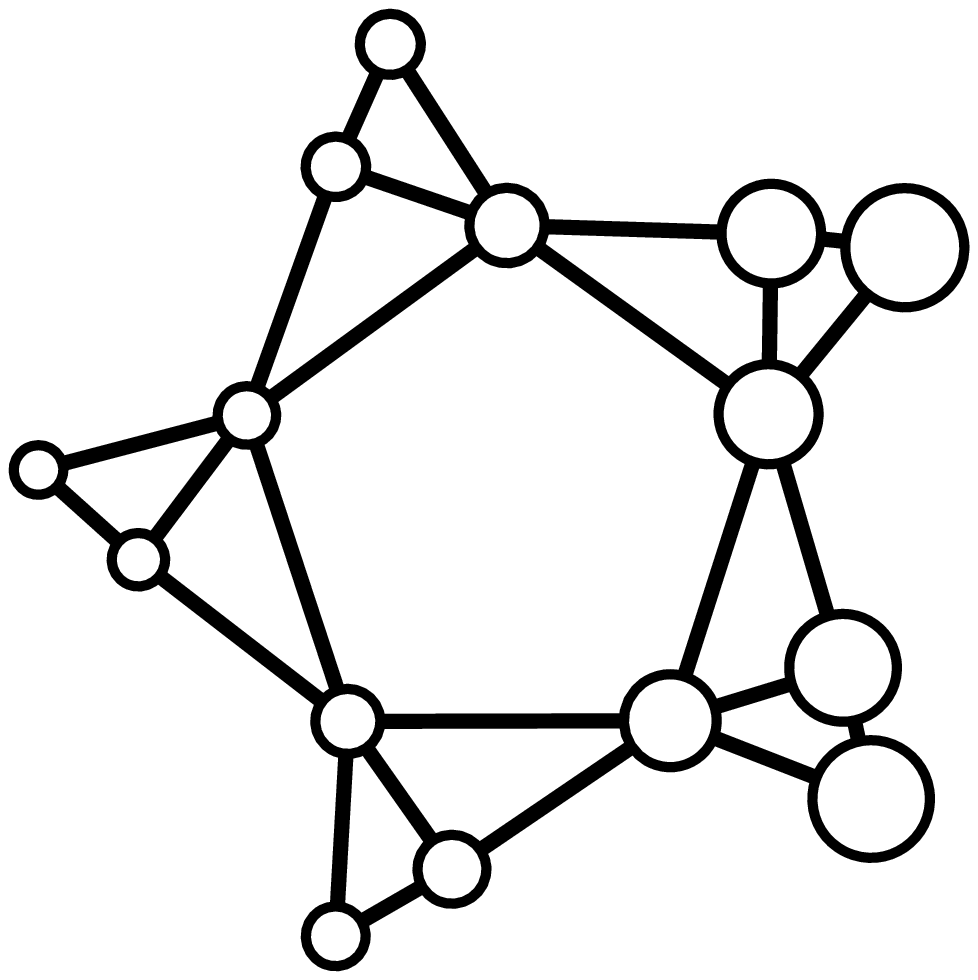} 
	       }
              }="v0c"
"O"
[r(6)] [u(0)]{  
               \rlap{
	       \kern 0.6cm
	       \raise .7cm
	       \hbox{$\scriptstyle \Gamma_{2}$}
	       }
               \scalebox{\psize}{
               \includegraphics{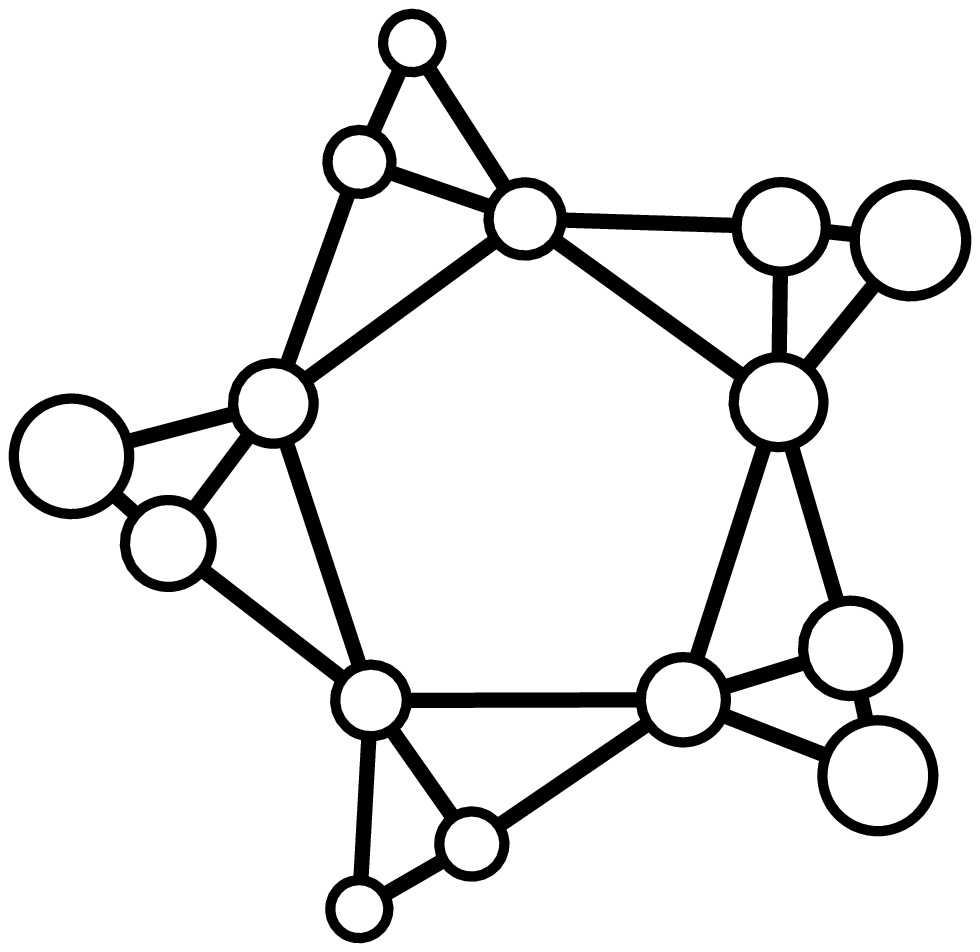} 
	       }
              }="v0d"
"O"
[r(5)] [u(3)]{
               \rlap{
	       \kern 0.6cm
	       \raise .7cm
	       \hbox{$\scriptstyle \Gamma_{1}$}
	       }
               \scalebox{\psize}{
               \includegraphics{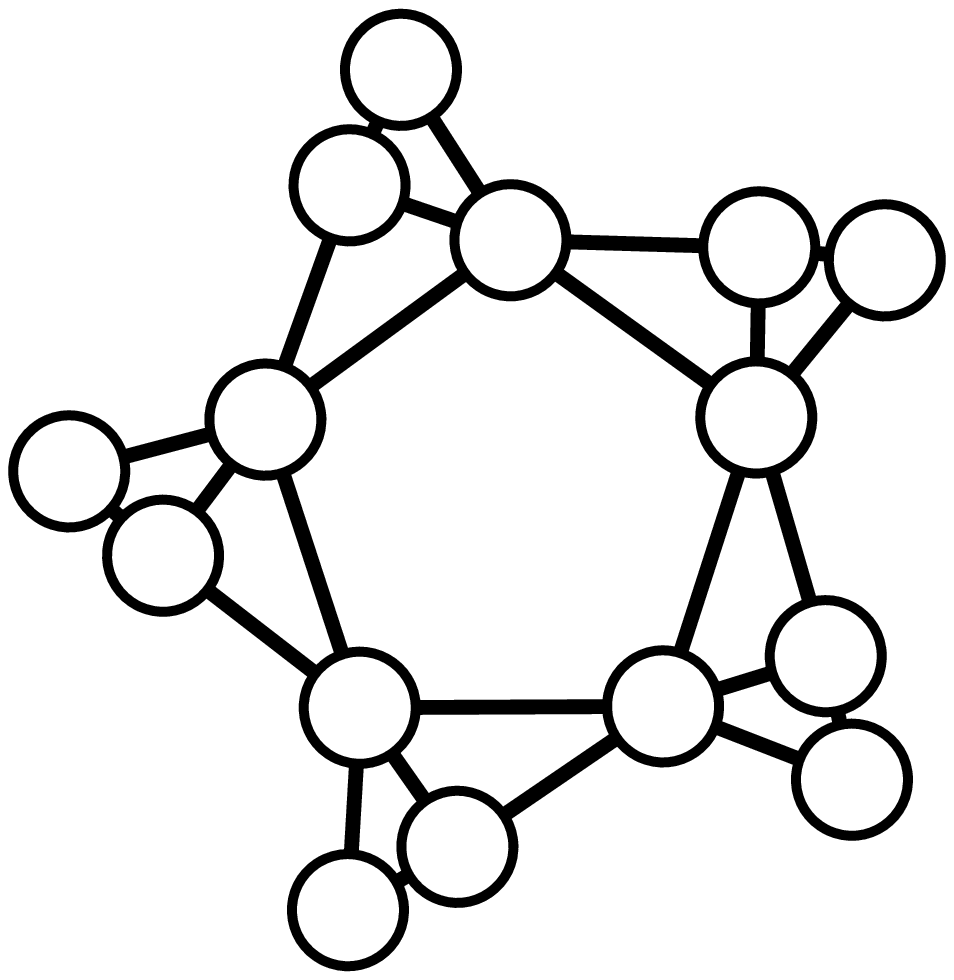} 
	       }
              }="v1"
"O"
[r(1)] [u(4)]{ 
               \rlap{
	       \kern 0.6cm
	       \raise .7cm
	       \hbox{$\scriptstyle \Gamma_{0}$}
	       }
               \scalebox{\psize}{
               \includegraphics{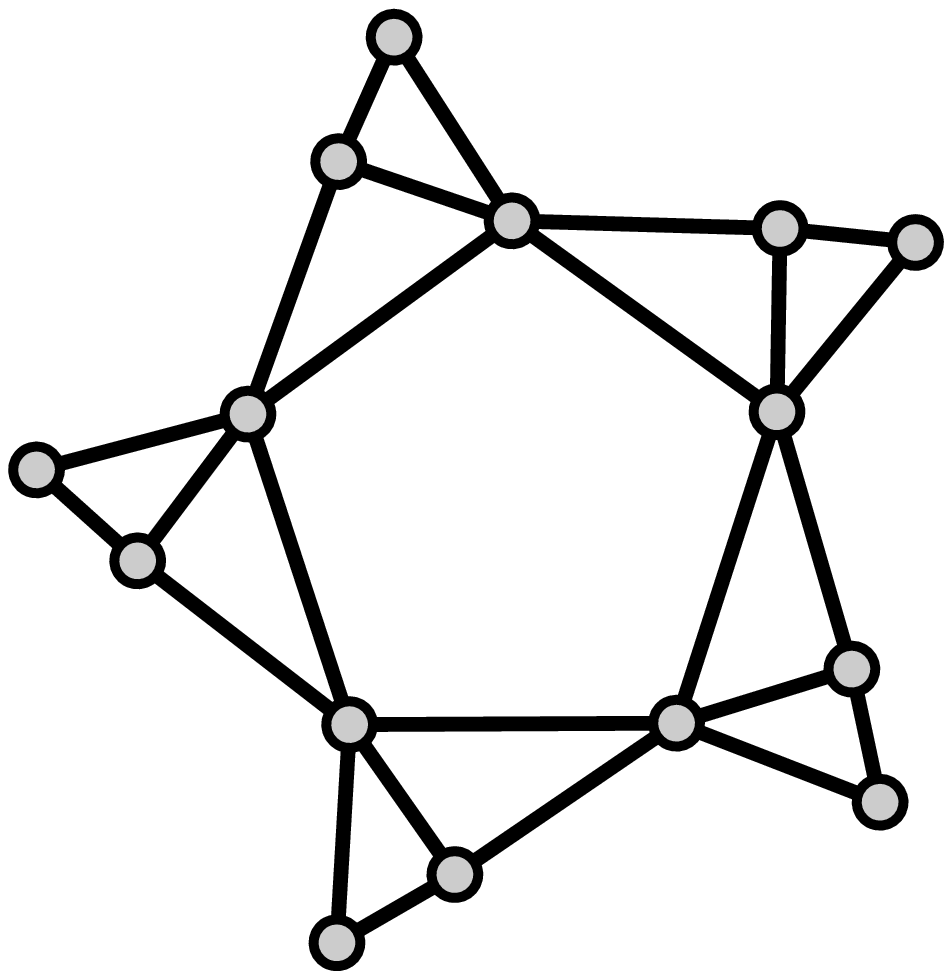} 
	       }
	      }="v2"
"v2":^1_{\Z_2}"v1"
"v1":@<-4pt>@{.>}_2^{\Z_5}"v0c"
"v1":@<4pt>@{.>}^3_{\Z_5}"v0d"
"v2":@<4pt>@{.>}
^<>(.7)
3
_<>(.7)
{\Z_{10}}"v0b"
"v2":@<-4pt>@{.>}
_<>(.7)
2
^<>(.7)
{\Z_{10}}
"v0a"
}
}
\caption
{
The irreducible representations of $\Z_{10}$ over $\R$,
the bifurcation digraph,
and all five bifurcation arrows for the Cayley graph of $\Z_5$
(Example~\ref{Cayley_Z5}).
The matrix of the rotation of $\R^2$ about the origin by $\theta$  
is denoted by $R(\theta)$. 
The bifurcation arrows are labeled with the irreducible representation $k$.
The arrow type and 
the group $\Gamma_i/\Gamma_{i, k}'$ of the bifurcation
are included in the arrow label
to facilitate comparison with the bifurcation digraph.
}
\label{Z5_digraph}
\end{figure}


In Figure~\ref{Z5_bifurcation_diagram} we show bifurcations corresponding
to the five bifurcation arrows.
For a gradient system,
a nondegenerate bifurcation with $\Z_{10}$ symmetry creates 20 daughter branches in two group orbits of size 10,
while a bifurcation with $\Z_5$ symmetry creates 10 daughter branches in two group orbits of size five.
These bifurcations are similar to bifurcations with $\D_{10}$ and $\D_5$ symmetry, respectively.
For example, a calculation shows that
the normal form \cite{GSS} for a gradient bifurcation with $\Z_5$ symmetry
is $g: \C \rightarrow \C$, $g(z) = \lam z \pm z |z|^2 + (a + i b) {\bar z}^4$,
where $\lam$, $a$ and $b$ are real.
The normal form for a bifurcation with $\D_5$ symmetry is the same,
but with $b = 0$ so that the real $z$
axis contains solutions to $g(z) = 0$.
The search directions in $E$ that lead to daughter solutions
are determined by the symmetry for $\D_5$, but
remain a matter of trial and error for $\Z_5$.
The situation is similar for the group $\Z_{10}$,
except that ${\bar z}^4$ is replaced by ${\bar z}^9$ in the normal form.

Our continuation solver had no trouble finding the bifurcating branches in this example.  To test this,
we ran the program that produced Figure~\ref{Z5_bifurcation_diagram} with
several different seeds in the random number generator used in Algorithm~\ref{find_daughters}.
In most cases, all branches would have been found with $f_{nc}(2) = 5$.
In all cases, the default $f_{nc}(2) = 21$ was large enough to find all of the branches.

The keen observer will notice that the $\Z_5$ bifurcation with $k = 2$ in Figure~\ref{Z5_bifurcation_diagram}
appears to have one branch bifurcating to the left
and one branch bifurcating to the right, in contradiction to the fact that nondegenerate $\Z_5$ bifurcations
have both of the daughter
branches bifurcating the same direction. 
However, an extreme blow-up of the figure shows
that both branches bifurcate to the left, but one branch has a fold point very close to the bifurcation point
in addition to the visible fold point.
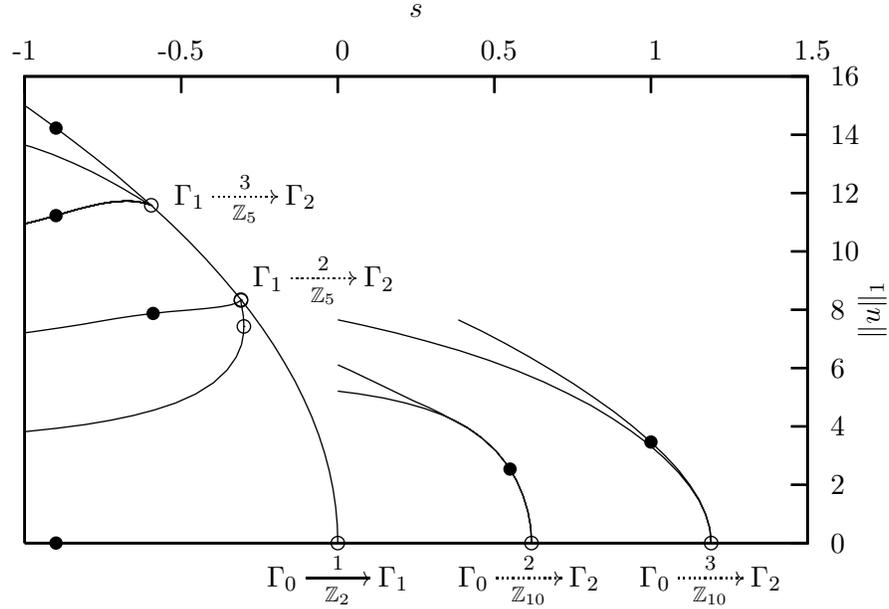
\begin{figure}
\def\are{$\xymatrix{\Gamma_0 \ar@{..>}[r]^3_{\Z_{10}} &\Gamma_2}$}
\def\ard{$\xymatrix{\Gamma_0 \ar@{..>}[r]^2_{\Z_{10}} &\Gamma_2}$}
\def\arc{$\xymatrix{\Gamma_1 \ar@{..>}[r]^2_{\Z_5} &\Gamma_2}$}
\def\arb{$\xymatrix{\Gamma_1 \ar@{..>}[r]^3_{\Z_5} &\Gamma_2}$}
\def\ara{$\xymatrix{\Gamma_0 \ar[r]^1_{\Z_{2}} &\Gamma_1}$}
\input Z5_bif.tex
\bigskip
\caption
{
Bifurcation diagram showing some of the solutions to Equation~(\ref{pde})
for the Cayley graph of $\Z_5$.
The primary branches created at the first three distinct eigenvalues
of $L$ are shown, along with the daughter branches from
two bifurcations with $\Z_5$ symmetry on the constant branch.
The arrows and the black dots correspond to the six solutions shown
in Figure~\ref{Z5_digraph}.  
}
\label{Z5_bifurcation_diagram}
\end{figure}

\end{subsection}

\begin{subsection}{The Cayley graph of the quaternion group $Q$}
\label{Qexample}
We are interested in the quaternion group
$Q = \langle i, j, k \mid i^2 = j^2 = k^2 = i j k = -1\rangle = \{\pm1,\pm i,\pm j,\pm k\}$
because it is the smallest group with a quaternionic representation
(see Section~\ref{isoSub}.)
We find several examples of
bifurcation with $Q$ symmetry, which are interesting and complicated.
Figure~\ref{Q_graphs} shows the Cayley color graph of $Q$
with generating set $\{i,j\}$ and the corresponding decorated Cayley graph. 
\begin{figure} \begin{center} \scalebox{1}{ \xygraph{ !{0;<2.2cm,0cm>:} []="O"
"O" [r(0)] [u(0)]{(j)}="v0" "O" [r(1)] [u(2)]{(k)}="v1" "O" [r(2)]
[u(2)]{(\text{-}j)}="v2" "O" [r(1)] [u(0)]{(\text{-}k)} ="v3" "O" [r(0)]
[u(2)]{(i)}="v4" "O" [r(0)] [u(1)]{(\text{-}1)}="v5" "O" [r(2)]
[u(0)]{(\text{-}i)}="v6" "O" [r(2)] [u(1)]{(1)}="v7" "v0":@{..>}"v3"
"v0":@{-->}"v5" "v1":@{..>}"v0" "v1":@{-->}"v6" "v2":@{..>}"v1" "v2":@{-->}"v7"
"v3":@{..>}"v2" "v3":@{-->}"v4" "v4":@{-->}"v1" "v4":@{..>}"v5" "v5":@{-->}"v2"
"v5":@{..>}"v6" "v6":@{-->}"v3" "v6":@{..>}"v7" "v7":@{-->}"v0" "v7":@{..>}"v4"
} } \qquad
\scalebox{.48}{\includegraphics{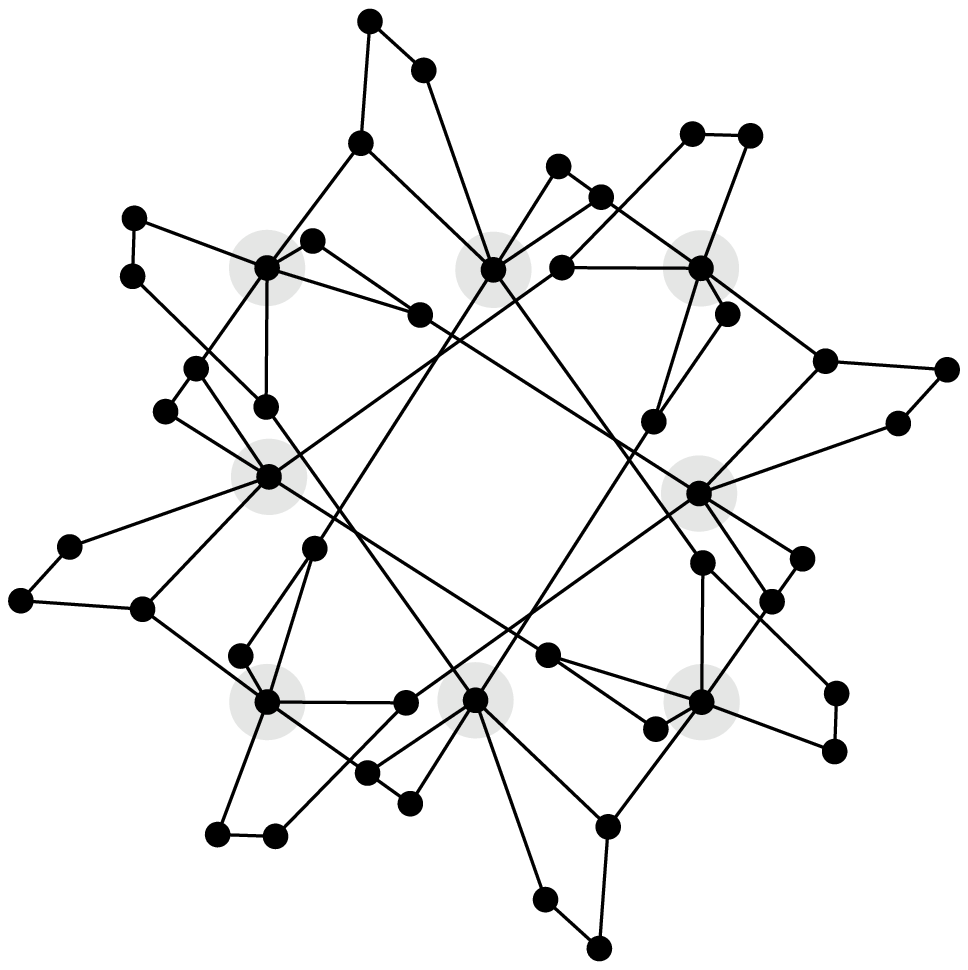}} \end{center} \medskip
\begin{center} \scalebox{1}{ \xygraph{ !{0;<1.0cm,0cm>:} []="O" "O" [r(0)]
[u(0)]{g}="v0" "O" [r(2)] [u(0)]{gi}="v1" "O" [r(6)] [u(0)]{g}="v2" "O" [r(7)]
[u(0)]{\bullet}="v3" "O" [r(8)] [u(0)]{gi}="v4" "O" [r(7.5)]
[u(0.7)]{\bullet}="v5" "O" [r(0)] [u(1.5)]{g}="v10" "O" [r(2)]
[u(1.5)]{gj}="v11" "O" [r(6)] [u(1.5)]{g}="v12" "O" [r(7)]
[u(1.5)]{\bullet}="v13" "O" [r(8)] [u(1.5)]{gj}="v14" "O" [r(7)]
[u(2.2)]{\bullet}="v15" "O" [r(8)] [u(2.2)]{\bullet}="v16" "v0":@{..>}^i"v1"
"v2"-"v3" "v3"-"v4" "v3"-"v5" "v4"-"v5" "v10":@{-->}^j"v11" "v12"-"v13"
"v13"-"v14" "v13"-"v15" "v14"-"v16" "v15"-"v16" } } \qquad \end{center}
\caption { The Cayley graphs of $Q$ (Example~\ref{Qexample}).  The graph on the
top left is the Cayley color digraph $\text{Cay}_{\{i,j\}}Q$.  The graph on the
top right is a decorated Cayley graph, which has 48 vertices and 72 edges.
The eight highlighted vertices in the decorated Cayley graph corresponding to the group elements $g \in Q$.
The bottom pictures show how the colored, directed edges are replaced with
undirected edges.  } \label{Q_graphs} \end{figure}

Figure~\ref{Q_digraph} shows the condensed bifurcation digraph computed by GAP.
There are 14 symmetries in $\mathcal G$, and each of the symmetry
types is a singleton:   $S_i = \{ \Gam_i \}$ for all $i$.
Since $-1 \in Q$, we need to use the notation
$(g, \pm 1)$, with $g \in Q$ to denote elements of $Q \times \Z_2$.
The symmetries isomorphic to $Q$ are
$$
\Gam_1 = \langle (i, 1), (j, 1) \rangle, \ 
\Gam_2 = \langle (i, -1), (j, 1) \rangle, \ 
\Gam_3 = \langle (i, 1), (j, -1) \rangle, \ \text{and }
\Gam_4 = \langle (i, -1), (j, -1) \rangle. \ 
$$
These four subgroups of $\Gam_0 = Q \times \Z_2$ are not conjugate,
but they are related by outer automorphisms
and the four symmetry types are in the same condensation class,
as seen in Figure~\ref{Q_digraph}.
The six symmetries isomorphic to $\Z_4$ are
$$
\Gam_{5} = \langle (i,  1) \rangle, \
\Gam_{6} = \langle (i, -1) \rangle, \
\Gam_{7} = \langle (j,  1) \rangle, \
\Gam_{8} = \langle (j, -1) \rangle, \
\Gam_{9} = \langle (ij,  1) \rangle , \text{ and }
\Gam_{10} = \langle (ij, -1) \rangle .
$$
There are two symmetries isomorphic to $\Z_2$, namely
$$
\Gam_{11} = \langle (-1, -1) \rangle, \text{  and } \Gam_{12} = \langle (-1, 1) \rangle.
$$
The symmetry types $S_{11}$ and $S_{12}$ are in different condensation classes.
Finally, the trivial symmetry is $\Gam_{13}$.

\begin{figure} \scalebox{1}{ \xygraph{ !{0;<.6cm,0cm>:} 
[]="O" 
"O" [r(8)][u(0)]{\fbox{$S_{13}$}}="v0" 
"O" [r(16)] [u(2)]{ \fbox{$ S_{12} $} }="v1"
"O" [r(0)] [u(2)]{\fbox{$ S_{11} $} }="v2" 
"O" [r(16)] [u(6)]{\fbox{$ \begin{aligned} & S_5, S_6, S_7, \\ 
                                                & S_8, S_9, S_{10} 
                                \end{aligned} $}}="v3" 
"O" [r(8)] [u(8)]{ \fbox{$ \begin{aligned} &S_1, S_2, \\ 
                                                 &S_3, S_4
			         \end{aligned} $} }="v4" 
"O" [r(0)] [u(10)]{\fbox{$S_0$}}="v5"
"v1":^{\Z_2}"v0" 
"v2":_{\Z_2}"v0" 
"v3":@{..>}^{\Z_4}_<>(.9)6"v0" 
"v3":^<>(.5){\Z_2}_<>(.8)6"v1"
"v4":@{..>}_Q_<>(.9)4"v0" 
"v4":^<>(.5){\Z_2}_<>(.1)3_<>(.9)2"v3" 
"v5":@{..>}_Q"v2" 
"v5":^{\Z_2}_<>(.1)4"v4" } }
\caption { The condensed bifurcation digraph for the Cayley graph of $Q$.
Since the symmetry types are singletons,
it is quite easy to deduce the bifurcation arrows from this figure
using the description of the symmetries in terms of generators.
For example, 
there are four bifurcation arrows emanating
from $\Gam_2$:
$\Gam_2 \rightarrow \Gam_{6}$,
$\Gam_2 \rightarrow \Gam_{7}$,
$\Gam_2 \rightarrow \Gam_{10}$, and
$\Gam_2 \rightarrow \Gam_{13}$.
} \label{Q_digraph} \end{figure}

\begin{figure}
\begin{center} \begin{tabular}{cc}
\scalebox{.65}{\includegraphics{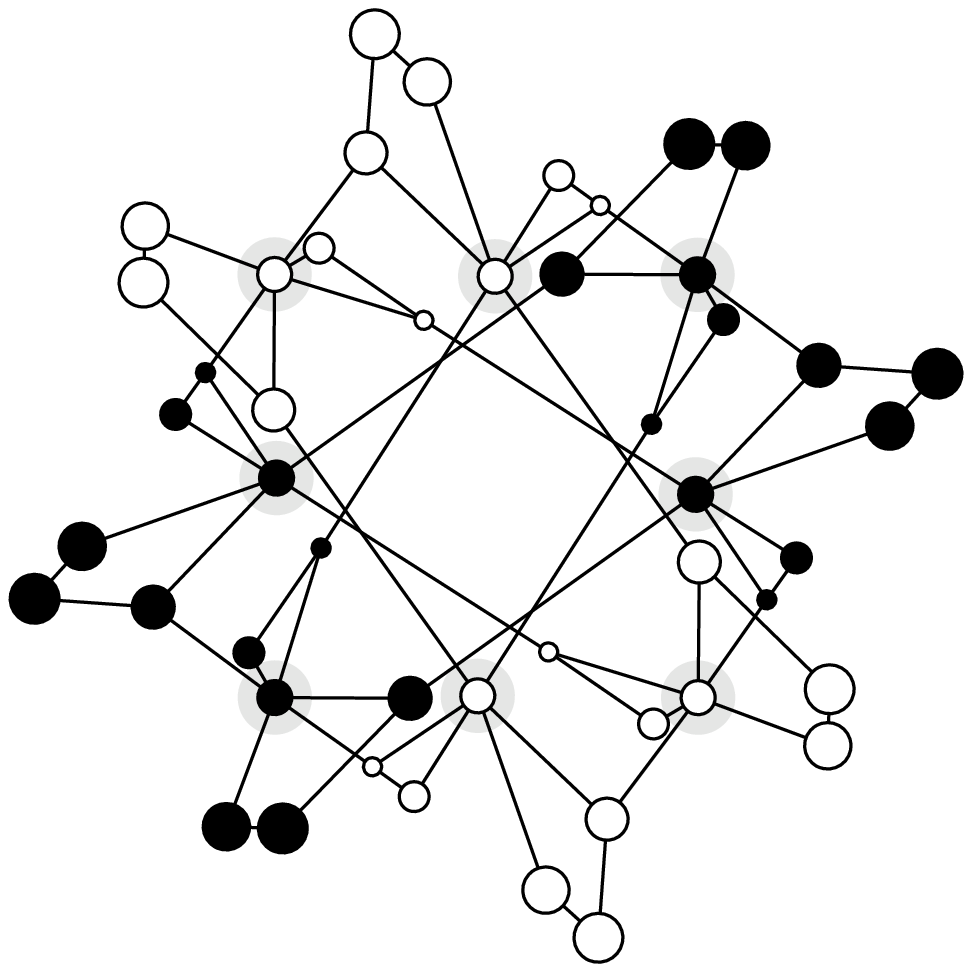}} &
\scalebox{.65}{\includegraphics{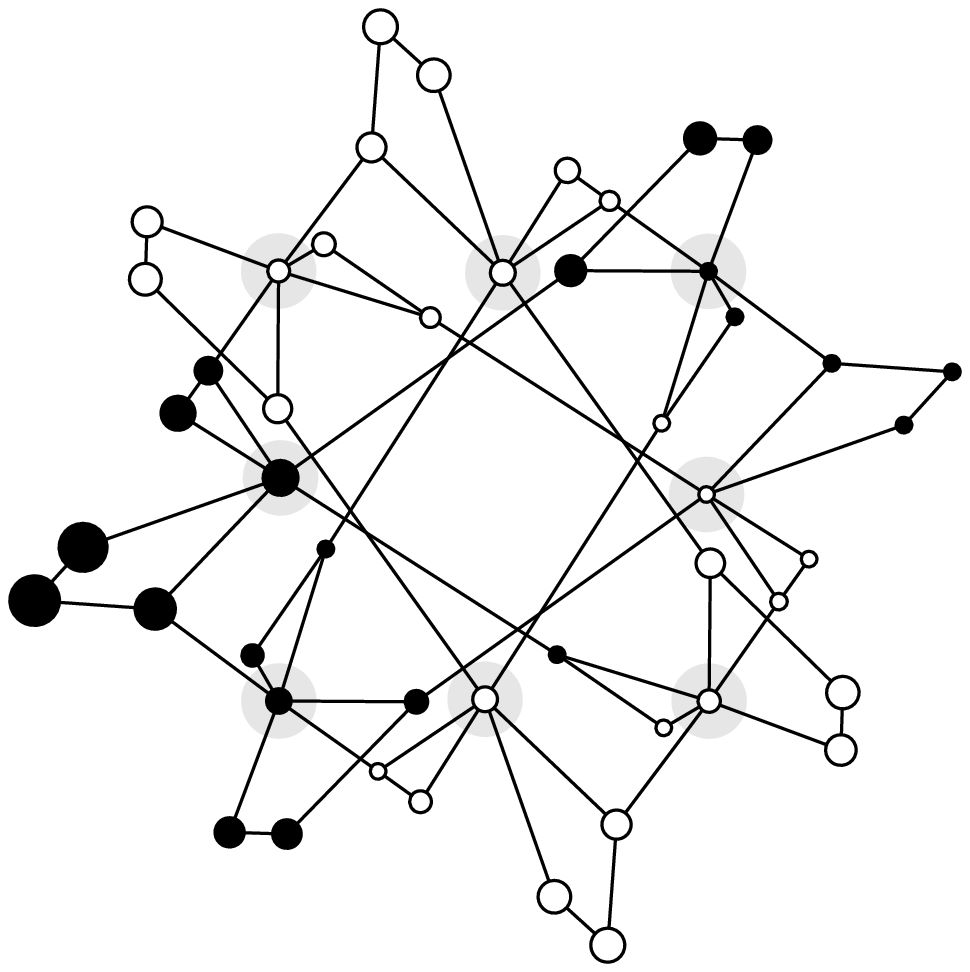}} \\
MI 2 for $s \in (s^*, \lam_2)$, MI 6 for $s < s^*$ &
MI 2 for $s < s^*$
\end{tabular}
\vspace*{.2in} \caption{
Contour maps of solutions on the decorated Cayley
graph of $Q$.
The solutions are the mother (left) and one of the 80 daughters (right)
of a bifurcation $\Gam_2 \rightarrow \Gam_{13}$ with $Q$ symmetry
at $s^* \approx 0.328$.
The mother solution lies on the primary branch born at $s = \lam_2 \approx 0.347$.
The symmetry of the solutions
can be determined from the value of $u$ at the 8 highlighted vertices.
For the mother solution, $|u|$ is the same at all highlighted vertices, and the action
of $i$, shown in Figure~\ref{Q_graphs}, switches the signs of $u$.
Similarly, the signs are unchanged under the action of $j$.
Hence the symmetry of the mother is $\Gam_2 = \langle (i, -1), (j, 1) \rangle$.
The daughter has trivial symmetry.
}
\label{bifWithQ}
\end{center}
\end{figure}

Bifurcations with $Q$ symmetry have a four-dimensional critical eigenspace $E$.
These bifurcations are ``highly non-EBL'',
since all points in the critical eigenspace,
except for the origin, have the same symmetry.
Our continuation solver found examples of each of the 5 bifurcations with $Q$ symmetry implied
by Figure~\ref{Q_digraph}:
$\Gam_0 \rightarrow \Gam_{11}$ and $\Gam_i \rightarrow \Gam_{13}$ for $i =$
1, 2, 3, and 4.
This is the first time, to the best of our knowledge,
that bifurcations with $Q$ symmetry have been observed.
Figure~\ref{bifWithQ} shows the mother and one of the daughter solutions
for a bifurcation $\Gam_2 \rightarrow \Gam_{13}$.
In a non-gradient system of differential equations, one would expect the bifurcation to create
periodic solutions, but
in our gradient system the daughters are solutions of $\nabla J_s = 0$
that come in group orbits of size 8.

We find the daughters by trial and error with repeated
calls to the {\tt cGNGA} function, as described in
Algorithm~\ref{find_daughters}.
The number of random guesses 
can be modified by the user if index theory of something else suggests that not all 
daughter branches are found.
In particular, the Poincar{\'e}-Hopf index theorem \cite{Arnold} implies that
\begin{equation}
\label{indexInvariant}
\sum_{(u,s^*-\eps)\in X} (-1)^{\text{MI}(u,s^*-\eps)} = 
\sum_{(u,s^*+\eps)\in X} (-1)^{\text{MI}(u,s^*+\eps)},
\end{equation}
where $\eps$ is chosen so that the only bifurcation with $s \in [s^* - \eps, s^*+\eps]$ is at $s = s^*$.
In practice, we only need to sum over the mother and daughter solutions of the bifurcation.
This invariant proved particularly useful at bifurcations with $Q$ symmetry.
The MI on the mother branch changes by four, and each daughter
has a group orbit of size 8,
so in fact we sum Equation~(\ref{indexInvariant}) over the set of nonconjugate daughters found by our continuation solver.
The normal form for this bifurcation has not been computed, to the best of our knowledge.
We have no theory predicting exactly how many daughters there are,
or where in the 4-dimensional critical eigenspace the daughter solutions lie.

Bifurcations with $Q$ symmetry are extremely complicated.
The daughters all go to the left in some examples,
but some go to the left and some go to the right in others.
We observed bifurcations with eight and ten nonconjugate daughter branches,
for a total of 64 and 80 daughters, respectively.
Some solutions had a very small basin of attraction in the cylinder, necessitating
a large number of random guesses.  The index invariant in Equation~(\ref{indexInvariant})
was key in recognizing that some solutions were initially missing.
However, index theory can never prove that all of the daughters have been found.

One such bifurcation with $Q$ symmetry occurred on the CCN branch.
We found
10 nonconjugate daughters, all branching to the left, with MI 2, 3, 3, 3, 4, 4, 4, 4, 5, and 5, respectively.
We had to increase {\tt num\_no\_changes} $= f_{nc}(4)$ to at least 6,500 to satisfy index theory;
an additional 360,000 calls to {\tt cGNGA} did not find any more daughters.
For brevity, only the bifurcating CCN solution is shown in Figure~\ref{bifWithQ}.

\end{subsection}

\begin{subsection}{The Petersen graph}
\label{petersensub}

\begin{figure}
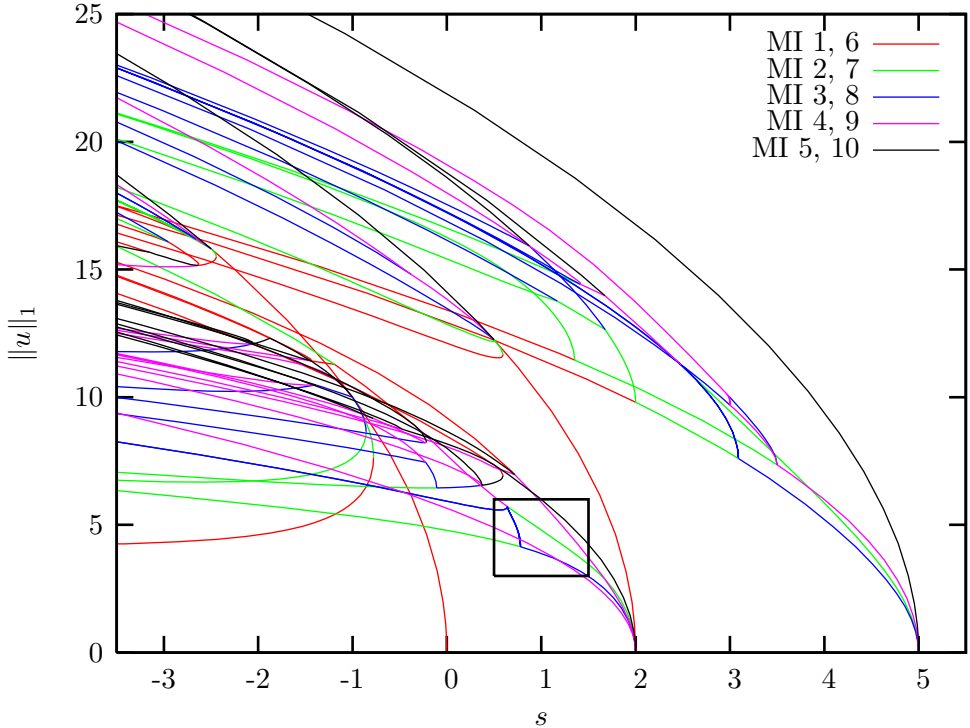

\begin{center}
\input Petersen.all.tex
\caption{Full bifurcation diagram for the Petersen graph (Example~\ref{petersensub}).  
There are high dimensional critical eigenspaces at $u=0$,  $s=\lambda_2=\cdots=\lambda_6=2$ and   $s=\lambda_7=\cdots=\lambda_{10}=5$,
as well at some secondary bifurcation points.  
Our code takes advantage of symmetry to search lower dimensional subspaces in order to efficiently find most solutions.
The green MI 2 and blue MI 3 branches bifurcating from the second eigenvalue (see inset)
contain a disconnected set of CCN solutions and are featured in Figure~\ref{ccn.Petersen}.
}
\label{Petersen.all}
\end{center}
\end{figure}

We considered the well known Petersen graph in our experiments.  The second eigenvalue for this graph Laplacian is of multiplicity 5 and there
are 210 symmetries (20 types) of possible solutions to look for, which
presents a challenge for our code.  We are fairly confident that we have accurately followed a representative from each equivalence class of primary
branches, and most if not all connected secondary branches.  Under greater magnification, 
Figure~\ref{Petersen.all} reveals that roughly 1300 points were used in following 75 branches, connected via 52 bifurcation points, 
with 3 reported (but not visually apparent) branch following failures.  
On our 3 GHz Linux workstation
it took about 3 seconds to perform
1726 calls to {\tt tGNGA}, with 4427 iterations at 2.565 iterations per call;
433 {\tt cGNGA} calls, with 2541 iterations at 5.868 iterations per call; and
  83 {\tt secant} calls, with 572 iterations at 6.892 iterations per call.

A particularly interesting feature of the bifurcation diagram regards solutions with the minimum $J$ value
among all sign-changing solutions for that $s$ parameter value, 
which are necessarily of MI 2 (see \cite{CCN}; for convenience we will call these CCN solutions here).
In \cite{CDN}, we proved that CCN solutions exist up to $\lambda_2$, and  
in \cite{Neu} we applied the GNGA to graphs and extended the CCN existence theorem (and related theorems) to graphs.
We have since conjectured that there should exist
a {\it connected branch} of CCN solutions for $s \in (-\infty, \lambda_2)$,
but this numerical experiment indicates that this is not true.
In Figure~\ref{ccn.Petersen} we show the
symmetry type $S_5$ and $S_{11}$ primary branches bifurcating from the multiplicity 5 second eigenvalue.
To the right of $s^*\approx0.694$, CCN solutions lie on the upper branch and have
symmetry type $S_5$,
whereas to the left they lie on the lower branch and have symmetry type $S_{11}$.
In Figure~\ref{ccnplots} we provide contour plots of these two solutions at the crossover point $s^*$,
where both are global minimizers of $J$ over the set of nontrivial sign-changing solutions.

\begin{figure}

\begin{center}

%
\input Petersen_ccn_bif.tex
\caption{
The CCN solutions for the Petersen graph, indicated by the thicker green lines, are not connected
(see the inset in Figure~\ref{Petersen.all}).
To the right of $s^*\approx 0.694$, CCN solutions have symmetry type $S_5$ and lie on the upper branch,
whereas to the left they lie on the lower branch and have symmetry type $S_{11}$.
The CCN solutions are global minimizers of $J$ over the set of sign-changing solutions,
and always have MI 2.
This is a numerical counterexample to our previous conjecture that a continuous branch of CCN solutions 
exists for $s<\lambda_2$.  
}
\label{ccn.Petersen}
\end{center}
\end{figure}

\begin{figure}
\begin{center}
$S_{5}$  \scalebox{.4}{\includegraphics{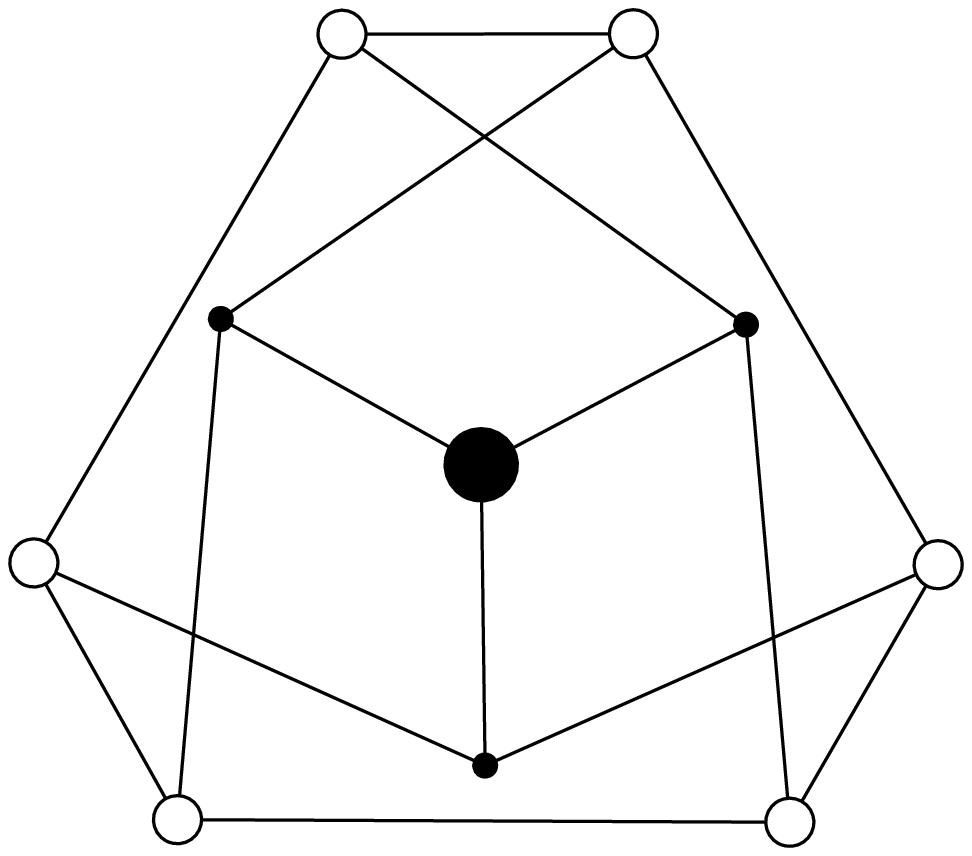}}  \ 
$S_{18}$ \scalebox{.4}{\includegraphics{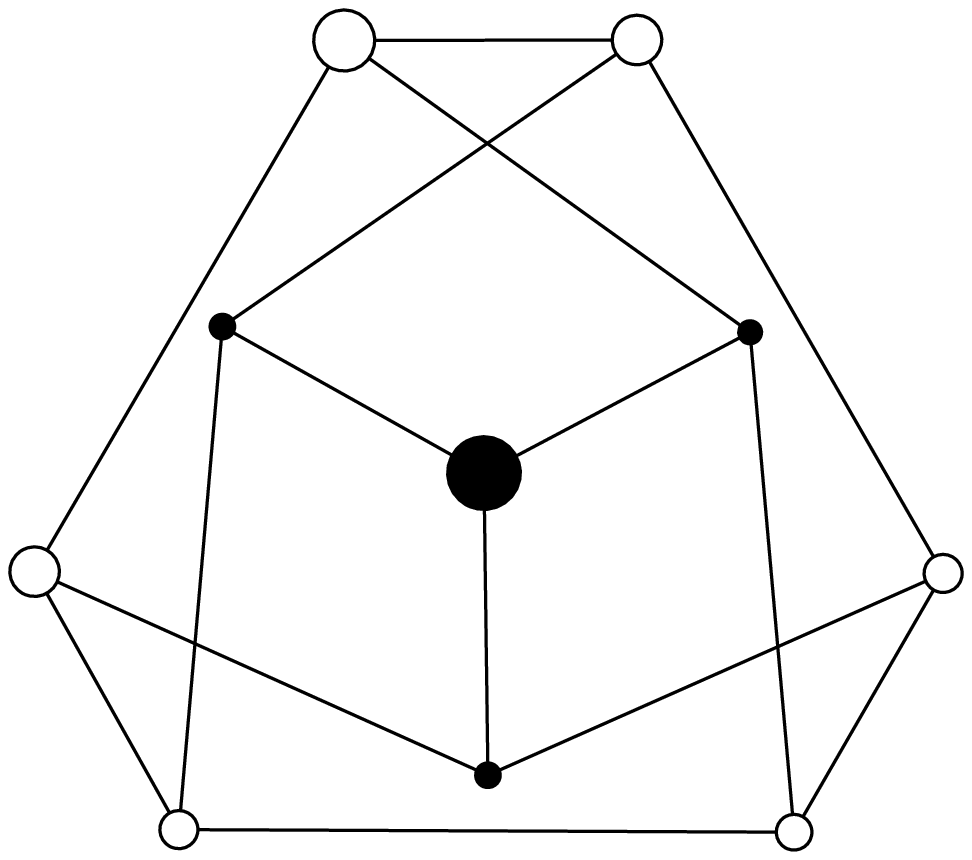}} \ 
$S_{11}$ \scalebox{.4}{\includegraphics{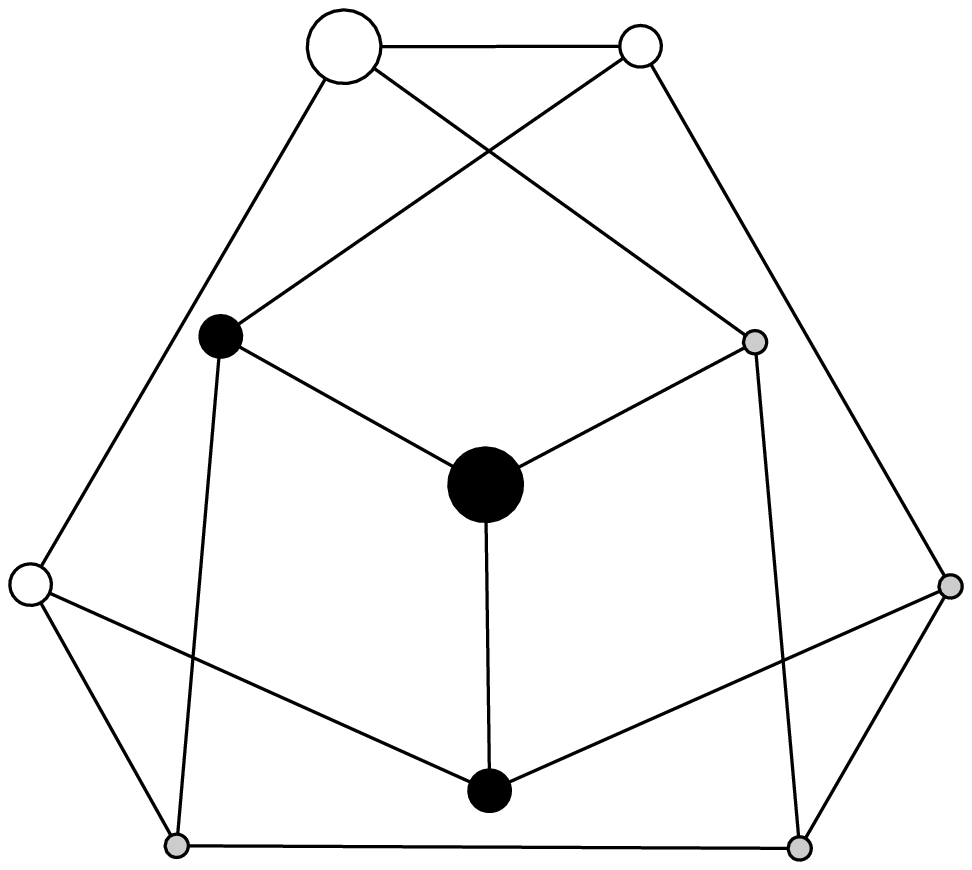}} 
\caption{
Contour plots of two simultaneous CCN solutions and a MI 3 solution from the connecting branch for the
Petersen graph,
corresponding to the dots in Figure~\ref{ccn.Petersen}.
The $S_5$ solution has 12 elements in its symmetry group, six of which are visible in this layout.
The $S_{11}$ and $S_{18}$ solutions have symmetry groups of size four and two, respectively, although
no nontrivial symmetries are visible in this layout.
This layout, as well as the traditional layout of the Petersen graph were found automatically
by our layout program.
The layout used here makes more of the symmetries of the $S_5$ solution visible.
}
\label{ccnplots}
\end{center}
\end{figure}

\end{subsection}

\begin{subsection}{Dodecahedron}
\label{dodeca_sub}
The space of functions on the dodecahedron graph admits 383 symmetries and 39 symmetry types.
The default layout found by our program is close to the
orthogonal projection of the 3-dimensional dodecahedron onto a plane parallel to a face (see Figure~\ref{type3Degeneracy}).

The dodecahedron features a Type 3 accidental degeneracy (see Definition~\ref{accidental}),
which are rare in the examples we studied. 
At this bifurcation point
the critical eigenspace is a direct sum of two 1-dimensional irreducible subspaces lying in
the same isotypic component.  That is,  ${\textsf K} = \{k\}$ is a singleton set and
$\dim(V_{\Gamma_i}^{(k)} \cap E) = 2$ whereas $d_{\Gamma_i}^{(k)} = 1$.

The accidental degeneracy in this example can be explained using AIS.
There is an AIS ${\mathcal A}_a$ for the dodecahedron comprising
of functions with $u_i = u_j$ if vertex $i$ and $j$ are antipodal.
It is well-known that the Petersen graph is the dodecahedron with antipodal points identified.
Therefore,
the nonlinear operator $\nabla J$ restricted to ${\mathcal A}_a$ is the same as $\nabla J$ acting
on functions on the Petersen graph.

The operator $\nabla J |_{{\mathcal A}_a}$ is equivariant under permutations on antipodal pairs of vertices.
There are 120 such permutations, since that is the size of the automorphism group of the Petersen graph.
However, only 60 of
these permutations are
symmetries of the dodecahedron. 
We say that the other 60 permutations are {\em anomalous symmetries} of $\nabla J |_{{\mathcal A}_a}$.

The numerical results for the Petersen graph can be used to understand solutions for the dodecahedron in ${\mathcal A}_a$.
At $s^* \approx 0.8727$ there is a bifurcation with $\Z_4$ symmetry for the Petersen graph,
at the same $s$ value, there is a 
degenerate bifurcation of Type 3 in the dodecahedron for which
the mother and daughters all lie in ${\mathcal A}_a$.
One might expect that the symmetry of the mother in the dodecahedron would be isomorphic to $\Z_4 \times \Z_2$, but in fact
it is $\Z_2 \times \Z_2$. Similarly, the symmetry of the bifurcation in the dodecahedron is not $\Z_4$, but rather
there are two simultaneous bifurcations with $\Z_2$ symmetry.

At the non-EBL bifurcation with $\Z_4$ symmetry in the Petersen graph,
8 daughters are created in two group orbits of size 4. 
One element in each group orbit is shown in Figure~\ref{type3Degeneracy}.
In the dodecahedron, the program also finds 8 daughters, but they lie in 4 group orbits of size two.
We only show two daughters, since the others are conjugate under anomalous symmetries.
\begin{figure}
\begin{center}
\begin{tabular}{ccc}
    \scalebox{.36}{\includegraphics{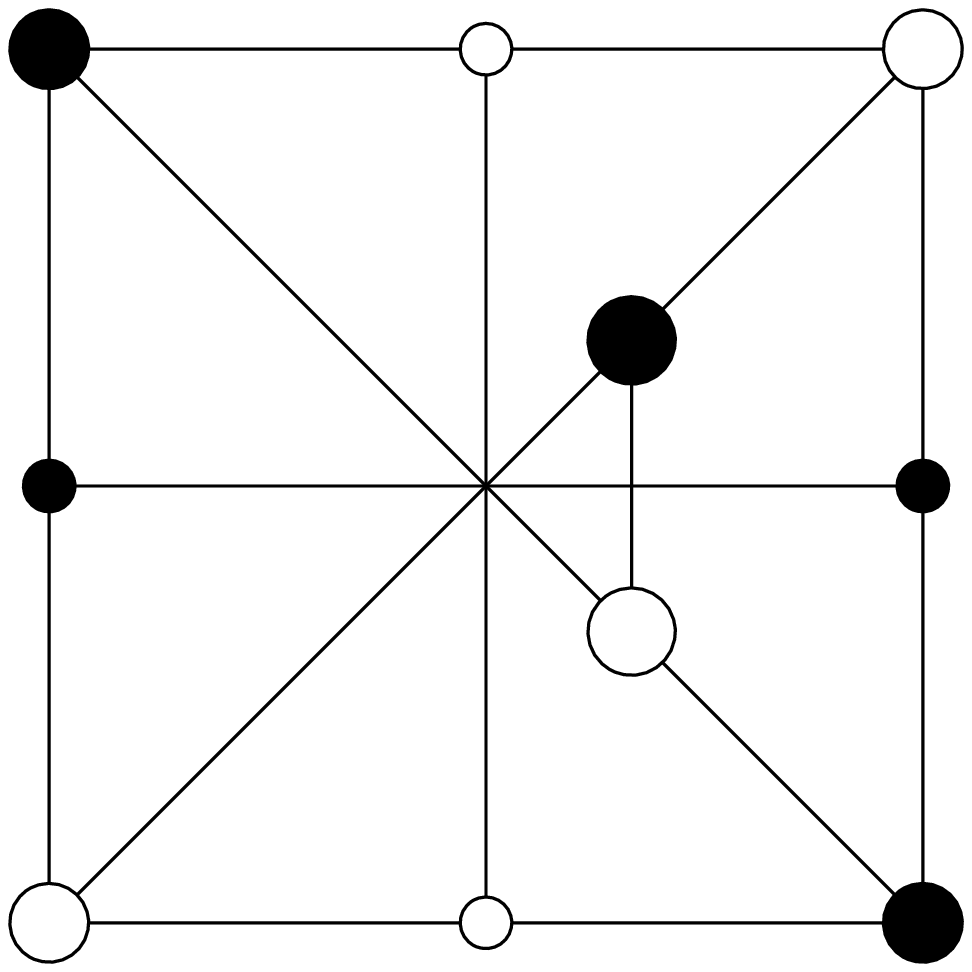}}
&
    \scalebox{.36}{\includegraphics{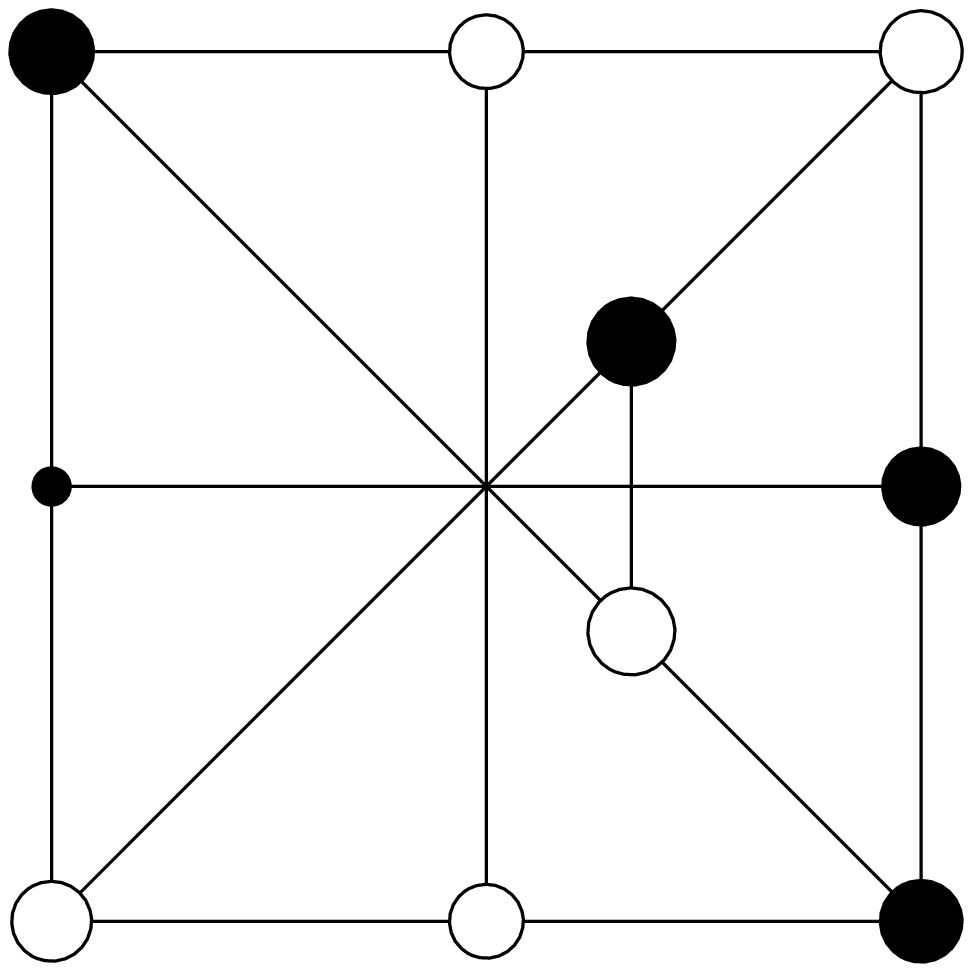}}
&
    \scalebox{.36}{\includegraphics{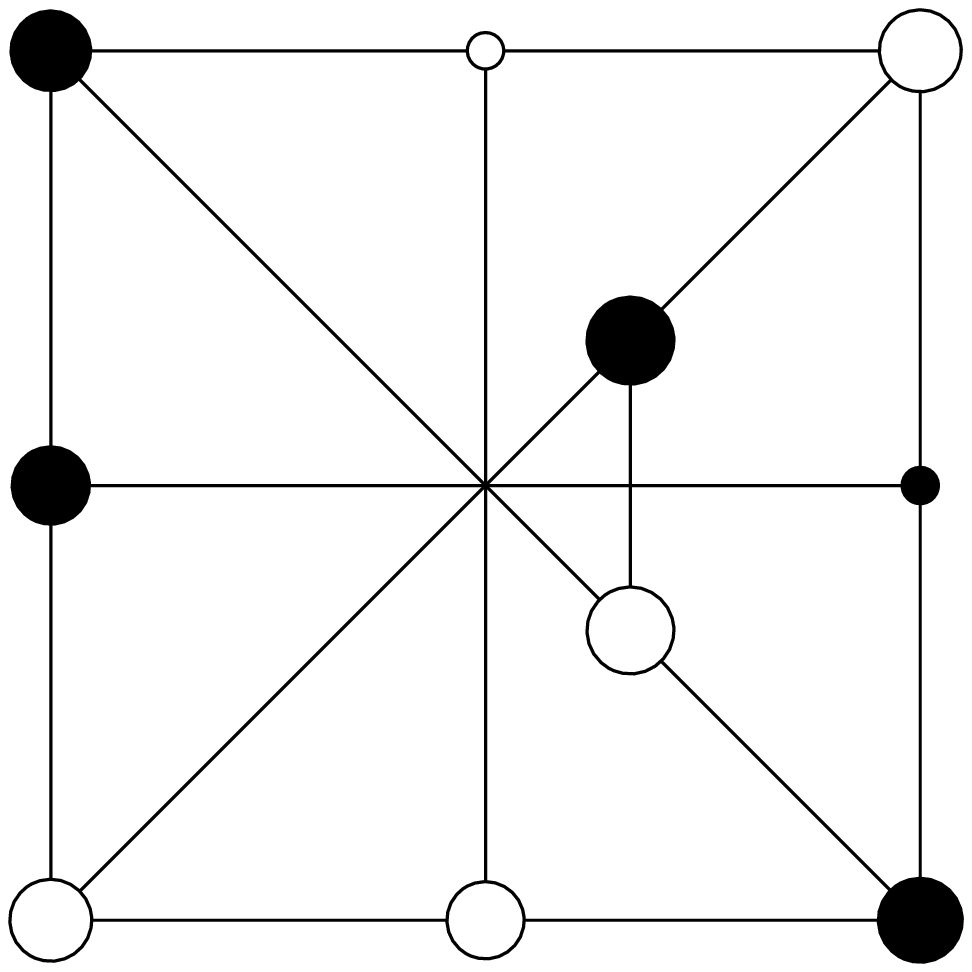}}
    \\
$\Gamma \iso \Z_4$ 
&
$\Gamma \iso \Z_1$
&
$\Gamma \iso \Z_1$\\

MI 10 for $s < s^*$ &
MI 9 for $s < s^*$ &
MI 8 for $s < s^*$ \\
MI 8 for $s > s^*$ & & \\
\\
    \scalebox{.45}{\includegraphics{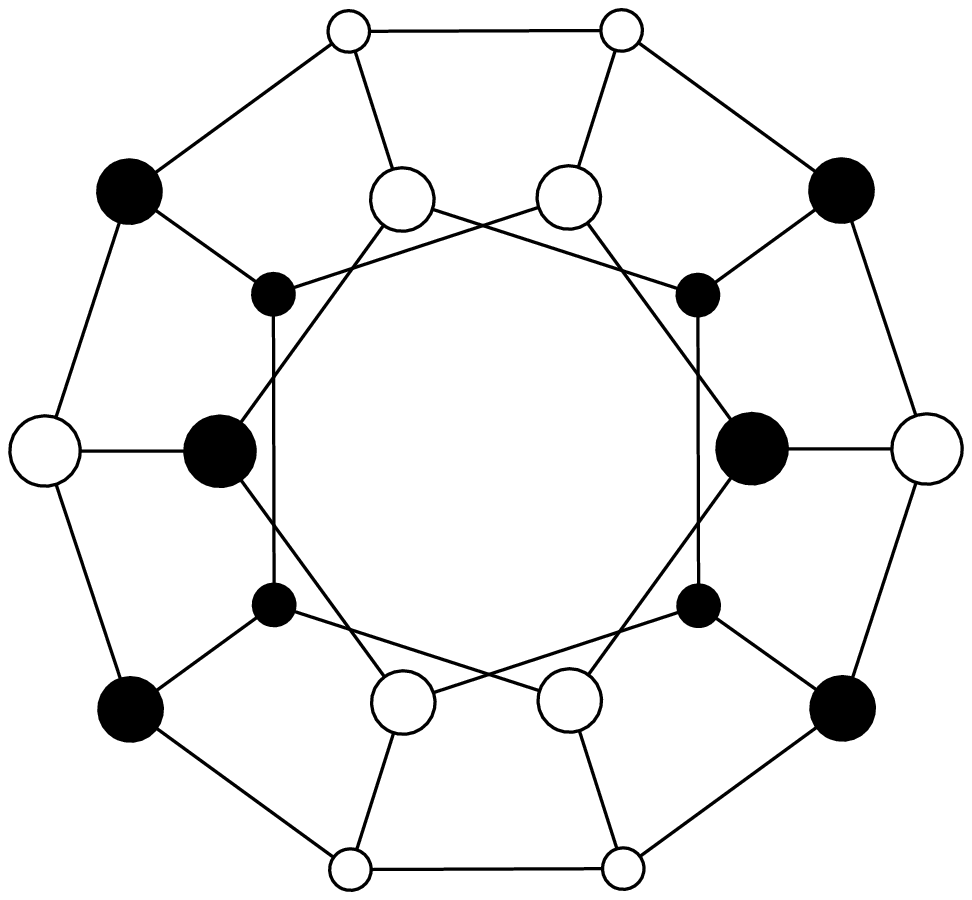}}
&
    \scalebox{.45}{\includegraphics{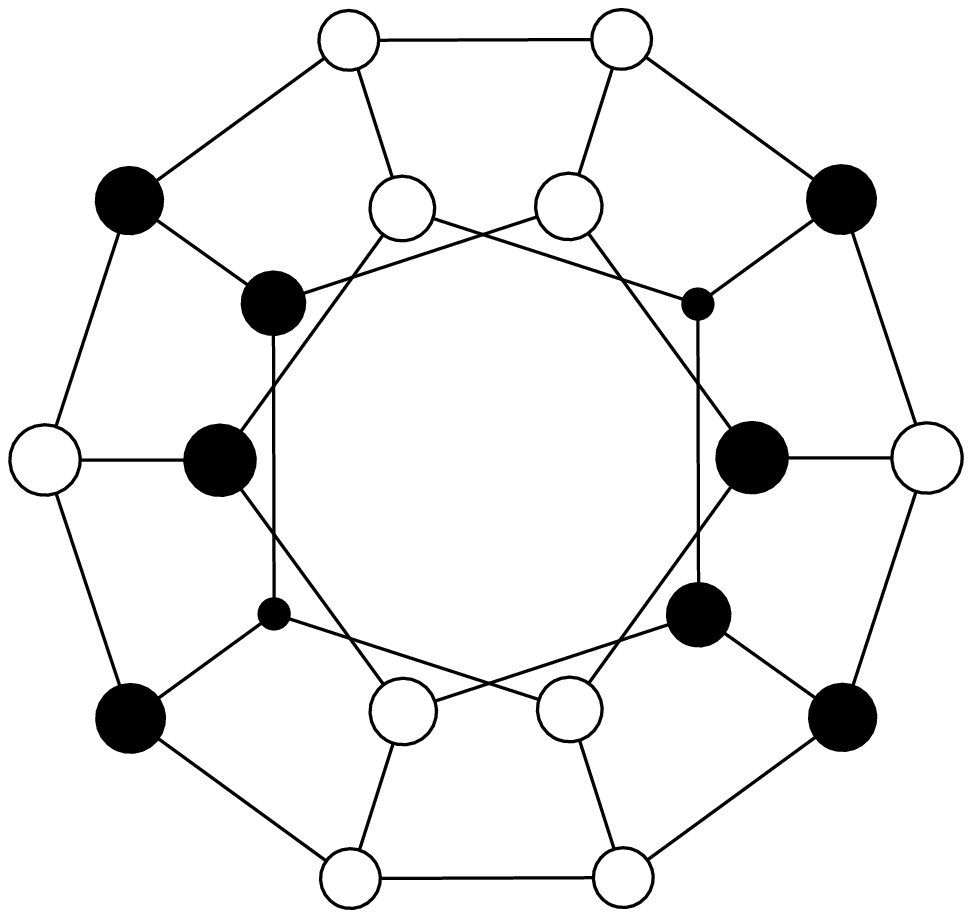}}
&
    \scalebox{.45}{\includegraphics{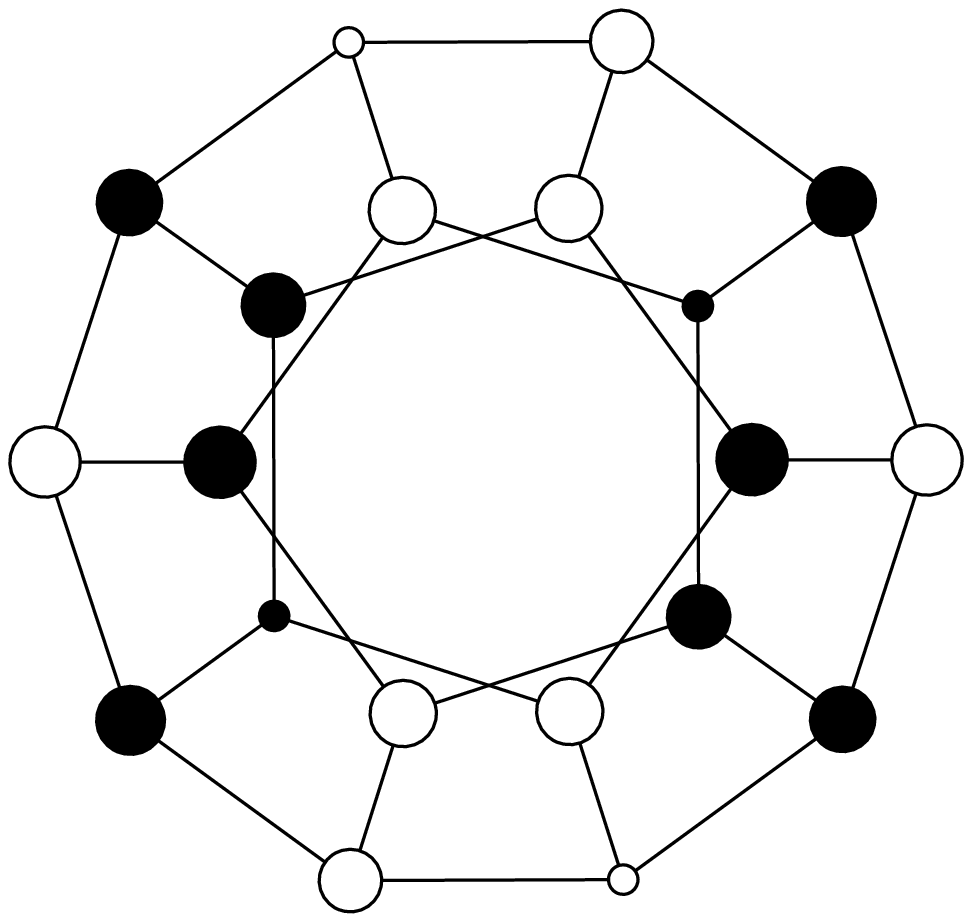}}
    \\
$\Gamma \iso \Z_2 \times \Z_2$ 
&
$\Gamma \iso \Z_2$
&
$\Gamma \iso \Z_2$\\

MI 18 for $s < s^*$ &
MI 17 for $s < s^*$ &
MI 16 for $s < s^*$ \\
MI 16 for $s > s^*$ & &
\end{tabular}
%
\caption{
Contour plots of solutions to Equation~(\ref{pde}) for the Petersen graph (top row) and dodecahedron (bottom row)
near the bifurcation point $s^* \approx 0.8727$.
The bifurcation has $\Z_4$ symmetry for the Petersen graph, but is a degenerate bifurcation of Type 3 for
the dodecahedron.
The MI information is valid on an interval $(s^* - \eps, s^* + \eps)$ local to the bifurcation.
The contour plots are obtained at $s = 1.5$ for the mother solutions (left column) and at $s = -2$ for the daughter solutions.
This layout of the Petersen graph shows the $\Z_4$ symmetry of
the mother solution.  The coordinates of the vertices in this layout were typed into a file rather than
automatically generated by our layout program.
Note that the four dots on the edge of the square are all the same size for the mother,
but the $\Z_4$ symmetry is broken for the two daughter solutions.
The bottom row shows corresponding solutions in the AIS of antipodal solutions ${\mathcal A}_a$
on the dodecahedron.
An interactive version of this figure, including three dimesional layouts of the Petersen graph
and the dodecahedron, can be found at
{\tt http://NAU.edu/Jim.Swift/PdE}.
}
\label{type3Degeneracy}
\end{center}
\end{figure}

\end{subsection}

\begin{subsection}{Truncated icosahedron (soccer ball)}
\label{soccerball_sub}

We include a final example with more vertices and a large number of
high multiplicity eigenvalues.
The truncated icosahedron, made famous via the Buckminsterfullerene molecule, 
has 60 vertices.
In Figure~\ref{soccerball_plots}
we display contour plots for 3 solutions from the second and third primary branches, near their bifurcations
from the trivial solution at $s=\lambda_2 = \lambda_3 = \lambda_4 \approx 0.2434$ and
$s = \lambda_5 = \cdots = \lambda_9 \approx 0.6972$, respectively.
This layout makes nearly all the symmetries of these solutions visible.
It was found by our graph layout code,
although it is not the layout with least complexity.
\begin{figure}
\begin{center}
\includegraphics[scale = 0.45]{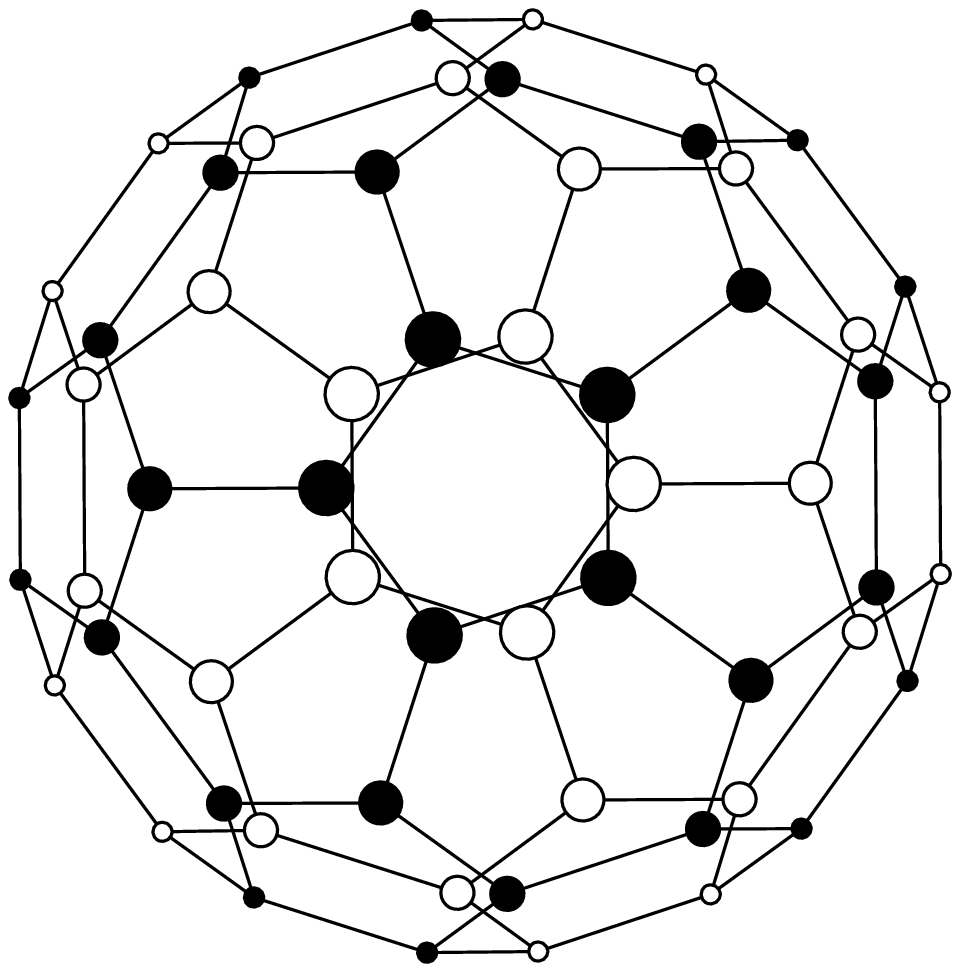}
\includegraphics[scale = 0.45]{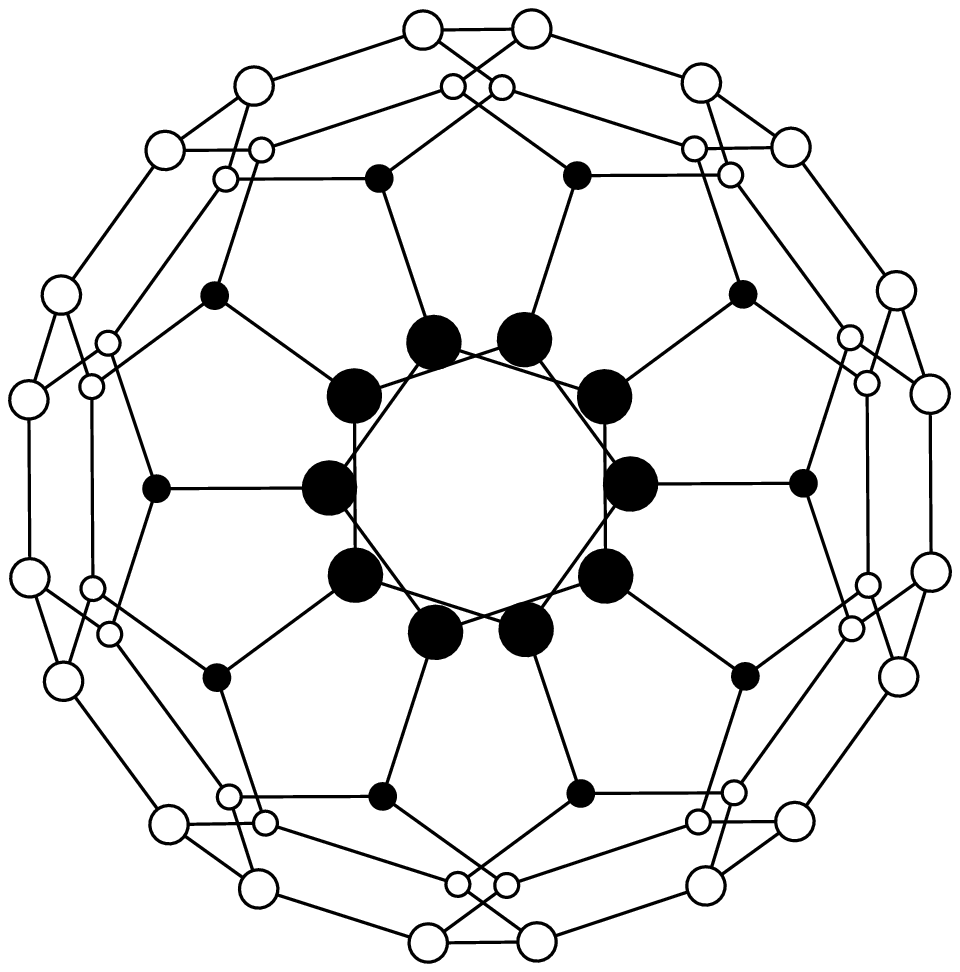}
\includegraphics[scale = 0.45]{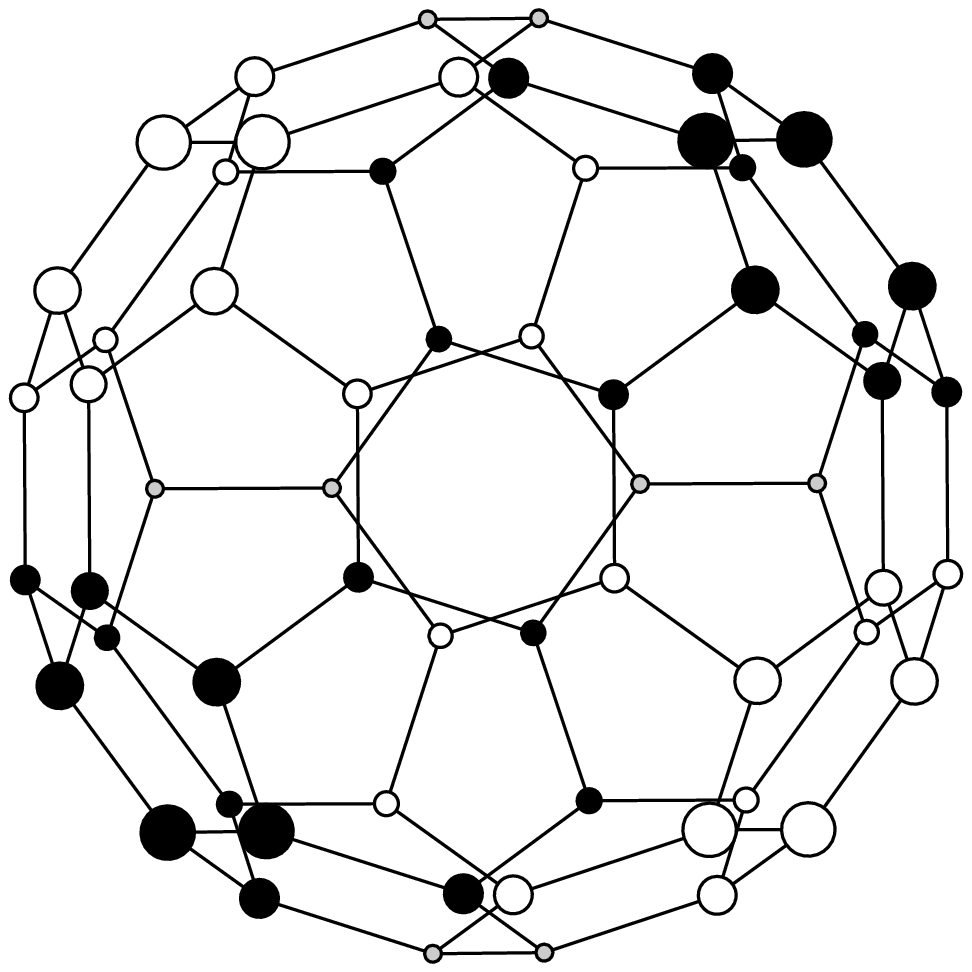}
\caption{
Contour plots of three solutions to Equation~(\ref{pde}) for the truncated icosahedron.
The MI 2 CCN solution (left) has 20 symmetries, all visible in this layout.
The front hemisphere is positive and the back is negative.
The MI 5 solution (center) has 20 of 20 visible symmetries as well.
It has a negative equatorial band separating front and back positive caps.
The MI 6 solution (right)
has 8 symmetries of which only 4 are visible.
The nodal structure with two positive and two negative components is clear.
All three solutions are very close to eigenvectors of $L$,
which in turn resemble  
eigenfunctions of the PDE Laplacian $-\Delta$ on the sphere.
Specifically, they have similar nodal structures as the spherical harmonics $Y_{1,0}$, $Y_{2,0}$, and ${\rm Re}\,(Y_{2,2})$, 
respectively \cite{boas}.
%
%
}
\label{soccerball_plots}
\end{center}
\end{figure}
\end{subsection}

\end{section}

\begin{section}{Future Research Directions.}
\label{conclusion}

Our suite of programs in their current state works well.
We achieved our goal of taking an edgelist as input and automatically generating
the symmetry information and solution data for PdE~(\ref{pde}) on the graph with that edgelist.
The figures in this paper required only minor formatting of the raw results.
We have tested our code on many other examples with a high degree of success.
We successfully automated, for general graphs, the symmetry analysis found in our
nonlinear snowflake code~\cite{NSS2}.  The results encoded in the bifurcation
digraph were used to follow most if not all bifurcating branches in these
gradient systems.

Bifurcation theory predicts that certain daughter branches must exist, but does
not rule out the existence of other branches.
Our continuation solver finds both types of daughter branches.
In Example~\ref{Qexample}
we used index theory as an indicator that not all of the daughters had yet been found.
It is an open problem to find more topological and variational theory to
predict in general how many branches bifurcate, and where they lie in the critical eigenspace.
Another open problem is to find a general theory of anomaly breaking bifurcations.

Our focus in the current project was on large groups, not large graphs.
For expedience we did not take advantage of all the methods used
in~\cite{NSS2} to speed up the calculations.
That code was very efficient, using hardcoded symmetry shortcuts to reduce
the number of integral calculations needed to define the linear systems required by Newton's method.
The Hessian matrix $h_s(u)$ is a block matrix,
with the number of zero blocks depending on the symmetry of $u$.
In the snowflake code these zero blocks were not computed, but in the current work 
it was tolerated to perform the calculation of each element of $h$.
We will implement this and other shortcuts
for solving PdE on significantly larger graphs,
specifically those obtained by discretizing PDE and using finite differences.
Our next project will start with the application of our automated branch following algorithms to
PDE on the square, which we first studied in \cite{NS}.

There are many PDE that merit an application of our PdE code.
One area of interest is PDE on fractal regions.
We propose to automatically generate 
large but finite pre-fractal graphs that in the limit converge to a fractal.  
Analytical and numerical research into the linear version of this problem has been done 
(see for example papers by R. S. Strichartz and A. Teplyaev and references therein),
but nonlinear research where bifurcation is considered is absent from the literature.  
One could also investigate large graphs embedded in 2D manifolds such as the torus and sphere
or 3D regions such as the cube.  
Generalizing $L$ to approximate the Laplace-Beltrami operator is a related idea for future investigation.
We are also interested in systems of PDE,
and have made initial demonstration programs for computing solutions to several systems.  Our code should
perform well in investigating these types of problems, depending in part on our understanding of the underlying theory for systems.  
Another area of future research
is to borrow from the established linear graph theory and our ghostpoint ideas from \cite{NSS} to accurately enforce 
alternate boundary conditions to our PdE (and hence PDE) code.

We are interested in the existence,
multiplicity, and nodal structure of solutions to nonlinear elliptic PDE.  
Thus, we will continue to perform experiments to support conjectures in the analytical theory for PDE,
seeking a better understanding of the underlying
variational structure. 
We found several interesting phenomena in the PdE examples and seek to determine if they persist for PDE.
For example,
in Section~\ref{petersensub} one sees an example where the CCN branches of solutions
for a PdE are disconnected, contrary to our conjecture that
a continuum of such solutions exists
for all $s < \lambda_2$
(assuming the standard subcritical/superlinear hypothesis found in~\cite{AR,CCN}).  
Furthermore, in Example~\ref{cycle_sub} we found strong numerical evidence that not all symmetry types
are present in the solution set $X$ for PdE.
It would be instructive to find examples of PDE with similar features.
Finally, we would like to find a PDE result analogous to the grouping by MI property commented on in Example~\ref{nosym_sub}. \\
\end{section}


\bibliography{nss3}

\def\cprime{$'$}
\providecommand{\bysame}{\leavevmode\hbox to3em{\hrulefill}\thinspace}
\providecommand{\MR}{\relax\ifhmode\unskip\space\fi MR }
\providecommand{\MRhref}[2]{%
  \href{http://www.ams.org/mathscinet-getitem?mr=#1}{#2}
}
\providecommand{\href}[2]{#2}
\begin{thebibliography}{10}

\bibitem{AR}
Antonio Ambrosetti and Paul~H. Rabinowitz, \emph{Dual variational methods in
  critical point theory and applications}, J. Functional Analysis \textbf{14}
  (1973), 349--381.

\bibitem{Arnold}
Vladimir~I. Arnol{\cprime}d, \emph{Ordinary differential equations}, Springer
  Textbook, Springer-Verlag, Berlin, 1992, Translated from the third Russian
  edition by Roger Cooke.

\bibitem{Bap}
R.~B. Bapat, \emph{The {L}aplacian matrix of a graph}, Math. Student
  \textbf{65} (1996), no.~1-4, 214--223.

\bibitem{Big}
Norman Biggs, \emph{Algebraic graph theory}, second ed., Cambridge Mathematical
  Library, Cambridge University Press, Cambridge, 1993.

\bibitem{boas}
Mary~L. Boas, \emph{Mathematical methods in the physical sciences}, Wiley, July
  2005.

\bibitem{CCN}
Alfonso Castro, Jorge Cossio, and John~M. Neuberger, \emph{A sign-changing
  solution for a superlinear {D}irichlet problem}, Rocky Mountain J. Math.
  \textbf{27} (1997), no.~4, 1041--1053.

\bibitem{CDN}
Alfonso Castro, Pavel Dr{\'a}bek, and John~M. Neuberger, \emph{A sign-changing
  solution for a superlinear {D}irichlet problem. {II}}, Proceedings of the
  Fifth Mississippi State Conference on Differential Equations and
  Computational Simulations (Mississippi State, MS, 2001) (San Marcos, TX),
  Electron. J. Differ. Equ. Conf., vol.~10, Southwest Texas State Univ., 2003,
  pp.~101--107 (electronic).

\bibitem{Chu}
Fan R.~K. Chung, \emph{Spectral graph theory}, CBMS Regional Conference Series
  in Mathematics, vol.~92, Published for the Conference Board of the
  Mathematical Sciences, Washington, DC, 1997.

\bibitem{Dornhoff}
Larry Dornhoff, \emph{Group representation theory. {P}art {A}: {O}rdinary
  representation theory}, Marcel Dekker Inc., New York, 1971, Pure and Applied
  Mathematics, 7.

\bibitem{FieldRichardson}
M.~J. Field and R.~W. Richardson, \emph{Symmetry breaking and the maximal
  isotropy subgroup conjecture for reflection groups}, Arch. Rational Mech.
  Anal. \textbf{105} (1989), no.~1, 61--94.

\bibitem{GS}
Victor~A. Galaktionov and Sergey~R. Svirshchevskii, \emph{Exact solutions and
  invariant subspaces of nonlinear partial differential equations in mechanics
  and physics}, Chapman \& Hall/CRC Applied Mathematics and Nonlinear Science
  Series, Chapman \& Hall/CRC, Boca Raton, FL, 2007.

\bibitem{GSS}
Martin Golubitsky, Ian Stewart, and David~G. Schaeffer, \emph{Singularities and
  groups in bifurcation theory. {V}ol. {II}}, Applied Mathematical Sciences,
  vol.~69, Springer-Verlag, New York, 1988.

\bibitem{GAP}
The~GAP Group, \emph{{GAP -- Groups, Algorithms, and Programming, Version
  4.4.9}}, { \verb+http://www.gap-system.org+}, 2006.

\bibitem{Lee}
Jason Lee, \emph{Existence of asymptotic solutions to semilinear partial
  difference equations on graphs}, {A Research Experience of Undergraduates
  Report, Northern Arizona University,\\
  \verb+http://math.nau.edu/researchInterests/students+}, 2007.

\bibitem{nauty}
Brendan~D. McKay, \emph{Practical graph isomorphism.}, {Numerical mathematics
  and computing, Proc. 10th Manitoba Conf., Winnipeg/Manitoba 1980, Congr.
  Numerantium 30, 45-87 (1981).}, 1981.

\bibitem{Neu2}
John~M. Neuberger, \emph{G{NGA}: recent progress and open problems for
  semilinear elliptic {PDE}}, Variational methods: open problems, recent
  progress, and numerical algorithms, Contemp. Math., vol. 357, Amer. Math.
  Soc., Providence, RI, 2004, pp.~201--237.

\bibitem{Neu}
\bysame, \emph{Nonlinear elliptic partial difference equations on graphs},
  Experiment. Math. \textbf{15} (2006), no.~1, 91--107.

\bibitem{NSS}
John~M. Neuberger, N{\'a}ndor Sieben, and James~W. Swift, \emph{Computing
  eigenfunctions on the {K}och snowflake: a new grid and symmetry}, J. Comput.
  Appl. Math. \textbf{191} (2006), no.~1, 126--142.

\bibitem{NSS2}
\bysame, \emph{Symmetry and automated branch following for a semilinear
  elliptic {PDE} on a fractal region}, SIAM J. Appl. Dyn. Syst. \textbf{5}
  (2006), no.~3, 476--507 (electronic).

\bibitem{NS}
John~M. Neuberger and James~W. Swift, \emph{Newton's method and {M}orse index
  for semilinear elliptic {PDE}s}, Internat. J. Bifur. Chaos Appl. Sci. Engrg.
  \textbf{11} (2001), no.~3, 801--820.

\bibitem{Rabinowitz}
Paul~H. Rabinowitz, \emph{Minimax methods in critical point theory with
  applications to differential equations}, CBMS Regional Conference Series in
  Mathematics, vol.~65, Published for the Conference Board of the Mathematical
  Sciences, Washington, DC, 1986.

\bibitem{WangZhou1}
Zhi-Qiang Wang and Jianxin Zhou, \emph{A local minimax-{N}ewton method for
  finding multiple saddle points with symmetries}, SIAM J. Numer. Anal.
  \textbf{42} (2004), no.~4, 1745--1759 (electronic).

\bibitem{WangZhou2}
\bysame, \emph{An efficient and stable method for computing multiple saddle
  points with symmetries}, SIAM J. Numer. Anal. \textbf{43} (2005), no.~2,
  891--907 (electronic).

\bibitem{White}
Arthur~T. White, \emph{Graphs, groups and surfaces}, second ed., North-Holland
  Mathematics Studies, vol.~8, North-Holland Publishing Co., Amsterdam, 1984.

\end{thebibliography}

\end{document}